%% file: main.tex
\documentclass[11pt]{article}
\usepackage[a4paper,margin=22mm]{geometry}
\usepackage{amsmath,amssymb,amsthm,mathtools,booktabs,array,tikz,microtype}
\usepackage[hidelinks,hypertexnames=false]{hyperref}
\hypersetup{pdftitle={Eigenvalues of doubly stochastic matrices in large orders},
pdfauthor={Brecht Verbeken and Vincent Ginis},
pdfsubject={The Perfect--Mirsky eigenvalue region}}
\newtheorem{theorem}{Theorem}
\newtheorem{lemma}[theorem]{Lemma}
\newtheorem{proposition}[theorem]{Proposition}

\newcommand{\cx}{\preceq_{\rm cx}}
\newcommand{\conv}{\operatorname{conv}}
\title{Eigenvalues of doubly stochastic matrices in large orders}
\author{Brecht Verbeken \and Vincent Ginis}
\date{}
\begin{document}
\maketitle
\begin{abstract}
For $1\le k\le N$, every point in the convex hull of the $k$th roots
of unity is an eigenvalue of some doubly stochastic matrix of order $N$.
Perfect and Mirsky asked whether any other eigenvalues occur. We prove
that the union of these polygons is the whole region for every $N\ge450$.
Our pencil-and-paper proof for $N\ge10^{30}$ starts with an eigenvector
$v$ of Euclidean norm one. Place one unit of mass at each coordinate
$v_j$. If $Av=zv$, row $i$ draws mass $A_{ij}$ from each mark $v_j$,
so the mass it receives has centre $zv_i$. The column sums ensure that each mark
supplies one unit in total. It follows that the sum of the row variances is
$1-|z|^2$. Since the regular $N$-gon contains the disk of radius
$\cos(\pi/N)$, an eigenvalue outside the polygons would force this sum
below $\sin^2(\pi/N)$. We use this bound to compare permutations of the
outer marks and to follow a randomly chosen mark through a quarter turn.
To reach $450$, we count marks between concentric circles and verify
a finite set of inequalities by computer. We also determine the regions
in orders five and six by computer-assisted proofs. At order five the
full region is obtained by adding to the polygons the segments from
the origin to two conjugate arcs traced by roots of quartic polynomials.
At order six the square, pentagon and hexagon give the whole region.
\end{abstract}

\begingroup
\setcounter{tocdepth}{1}
\tableofcontents
\endgroup

\section{What the polygon remembers}\label{sec:object}
A real square matrix is \emph{row-stochastic} if its entries are
nonnegative and each row sums to one. It is \emph{doubly stochastic}
if each column also sums to one. Start with a row-stochastic matrix $A$.
If $Av=zv$
with $v\ne0$,
the $i$th row expresses $zv_i$ as an average of the complex numbers
$v_1,\dots,v_N$. Draw their convex hull $K$. Multiplication by $z$
takes each marked point into $K$, and therefore $zK\subseteq K$.
Conversely, if $zK\subseteq K$, we can choose one such average for each
point. The coefficients give a stochastic matrix with eigenvector $v$.
Karpelevi\v c determined the eigenvalue region for row-stochastic
matrices~\cite{karp}. The invariant-polygon proof of that theorem
uses this equivalence: the vertices of $K$ already suffice to build
a realizing matrix of their own order~\cite{vg}.

Mirsky asked for the eigenvalue region when the columns also sum to
one~\cite{mirsky63}. What does that extra condition ask of the polygon?
Place one unit of mass
at each $v_j$, even when several coordinates occupy the same point.
In forming row $i$, use the amount $A_{ij}$ from site $j$. That row
receives one unit and its centre of mass is $zv_i$. Allocating mass to one row limits what is available to the others. The column sums say exactly that, once
all $N$ rows have received their unit, the supply at every site has
been exhausted. We call the list of sites with these equal masses a
\emph{marked configuration}. We say that multiplication by $z$ leaves the
configuration invariant if the mass can be distributed in this way.
We always require $v\ne0$. For $z\ne1$ the configuration is centred,
since summing $Av=zv$ over the rows gives $(1-z)\sum_i v_i=0$.

\paragraph{Two marks inside a square.}
Take the vertices $1,i,-1,-i$ and put two further marks at zero.
Multiplication by $i$ simply permutes the six marks. Move the two central
marks to $a$ and $-a$, with $0<a<1$. The same square is still invariant
under a quarter turn. How much of this turn can the six masses realise?
If $it$ is a multiplier with $t\ge0$, the rows at $1$ and $a$ ask for
positive imaginary parts $t$ and $ta$. The only source with positive
imaginary part is $i$, which has just one unit. Hence $t(1+a)\le1$.
Equality is possible. Send $t$ and $ta$ units of $i$ to those rows,
and distribute $-i$ in the same proportions to the opposite rows.
In the row at $i$, take $(1-t)/2$ units from $1$ and $(1+t)/2$ from
$-1$, giving the mean $-t$. Reverse these weights in the row at $-i$.
Use the two untouched masses at $a,-a$ in equal proportions to complete
the four partially filled rows. Thus the furthest point on the positive
imaginary axis changes from $i$ to $i/(1+a)$.

\begin{figure}[htbp]
\centering
\begin{tikzpicture}[x=1.45cm,y=1.45cm,font=\small]
\begin{scope}
\node at (0,1.53) {The six marks, $a=1/2$};
\draw[fill=black!3,thick] (1,0)--(0,1)--(-1,0)--(0,-1)--cycle;
\draw[gray!55] (-1,0)--(1,0);
\foreach \p/\pos/\name in {(0,1)/above/i,(0,-1)/below/-i,(-1,0)/left/-1,(-.5,0)/below/-a}{
  \fill \p circle (1.9pt) node[\pos=2pt] {$\name$};}
\fill[teal!75!black] (.5,0) circle (2.5pt) node[below=3pt,black] {$a$};
\fill[teal!75!black] (1,0) circle (2.5pt) node[right=2pt,black] {$1$};
\end{scope}
\begin{scope}[xshift=6.9cm]
\node at (0,1.53) {Sums of two distinct marks};
\draw[fill=teal!6,thick] (1.5,0)--(1,1)--(-1,1)--(-1.5,0)--(-1,-1)--(1,-1)--cycle;
\draw[gray!60] (-1.5,0)--(1.5,0);
\draw[teal!75!black,line width=1.35pt] (-1,1)--(1,1);
\draw[->,thick] (0,0)--(0,.93);
\fill[teal!75!black] (0,1) circle (2.5pt) node[above=2pt,black] {$i$};
\fill[teal!75!black] (1.5,0) circle (2.5pt) node[right=2pt,black] {$1+a$};
\node[below left] at (0,0) {$0$};
\end{scope}
\end{tikzpicture}
\caption{The rows at the two highlighted marks on the left require
the sum $it(1+a)$. On the right, the convex hull of all pair sums
has top edge at height one. The required sum reaches that edge at
$t=2/3$.}\label{fig:square}
\end{figure}
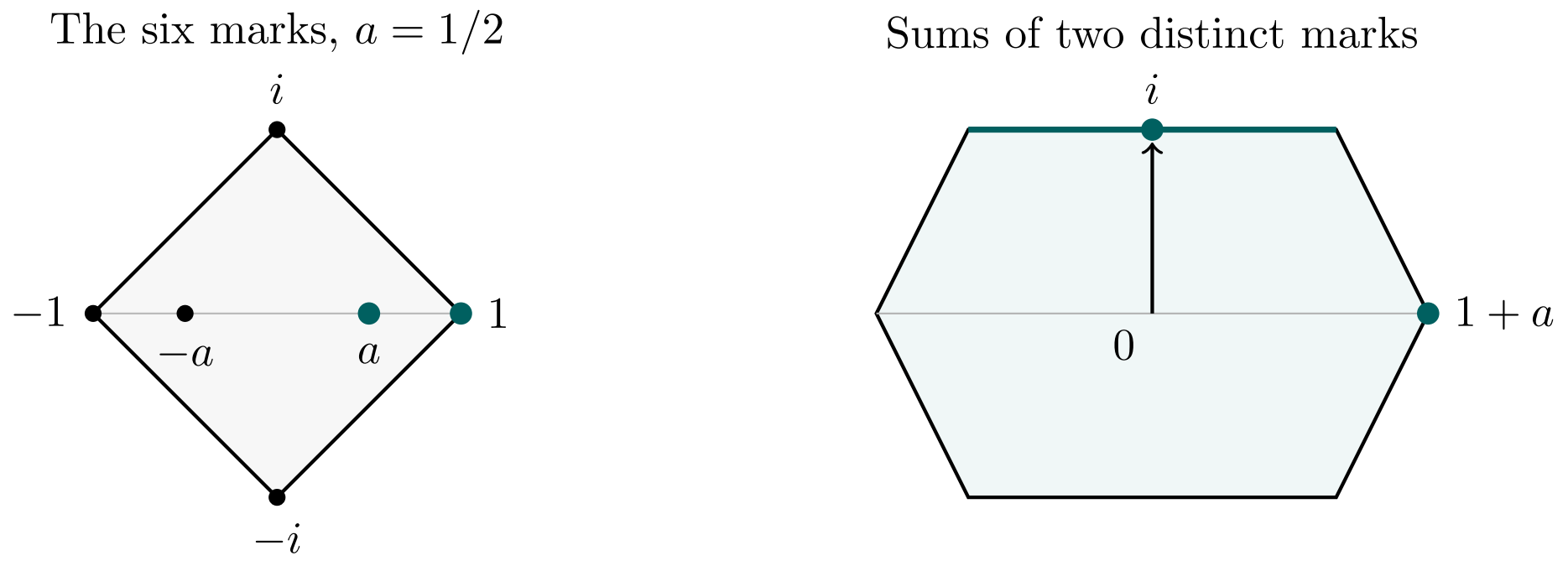

\paragraph{The region to be determined.}
Let $D_N$ be the set of complex numbers which occur as an eigenvalue
of at least one doubly stochastic matrix of order $N$. Put
\[
 \Pi_j=\conv\{e^{2\pi i k/j}:0\le k<j\},\qquad
 U_N=\bigcup_{j=1}^{N}\Pi_j .
\]
Here $\Pi_1=\{1\}$ and $\Pi_2=[-1,1]$. The polygons have a common vertex at $1$. They need not contain one
another. For instance $i$ belongs to $\Pi_4$ and lies outside $\Pi_3$,
and $e^{2\pi i/3}$ lies outside $\Pi_4$. Thus the union remembers
which roots can occur on a single cycle. Every point of $\Pi_j$ is
an eigenvalue. Indeed, put $\zeta=e^{2\pi i/j}$ and let $C$ shift
the entries of $(1,\zeta,\ldots,\zeta^{j-1})^T$ one place to the left,
moving the first entry to the end. This multiplies the vector by $\zeta$.
The powers $I,C,\ldots,C^{j-1}$ therefore give all the vertices on
the same vector. Take the convex combination giving the desired point and adjoin
$N-j$ fixed coordinates. This proves $U_N\subseteq D_N$.

\begin{figure}[htbp]
\centering
\begin{tikzpicture}[font=\small]
\begin{scope}[scale=1.42]
\node at (0,1.56) {$U_4=\Pi_3\cup\Pi_4$};
\fill[teal!9] (0:1.07)--(90:1.07)--(180:1.07)--(270:1.07)--cycle;
\fill[teal!9] (0:1.07)--(120:1.07)--(240:1.07)--cycle;
\draw[teal!75!black,thick] (0:1.07)--(90:1.07)--(180:1.07)--(270:1.07)--cycle;
\draw[black,semithick,dashed] (0:1.07)--(120:1.07)--(240:1.07)--cycle;
\fill (0:1.07) circle (1.8pt) node[right=2pt] {$1$};
\fill (90:1.07) circle (1.8pt) node[above=1pt] {$i$};
\fill (120:1.07) circle (1.8pt) node[left=2pt] {$e^{2\pi i/3}$};
\end{scope}
\begin{scope}[xshift=6.5cm,scale=1.42]
\node at (0,1.56) {The disk inside $\Pi_8$};
\draw[gray!55] (0,0) circle (1.07);
\draw[fill=teal!8,draw=teal!75!black,semithick] (0,0) circle[radius={1.07*cos(22.5)}];
\draw[thick] (0:1.07) \foreach \a in {45,90,...,315} {--(\a:1.07)} --cycle;
\draw[thick] (0,0)--(0:1.07);
\draw[thick,teal!75!black] (0,0)--(22.5:{1.07*cos(22.5)});
\fill (0,0) circle (1.4pt) node[below left] {$0$};
\fill (22.5:{1.07*cos(22.5)}) circle (1.8pt);
\draw (0:.29) arc(0:22.5:.29);
\node[anchor=west] at (.34,.11) {\scriptsize $\pi/8$};
\end{scope}
\end{tikzpicture}
\caption{The union on the left includes the two corners of the triangle
outside the square. On the right the radius to a side midpoint makes
an angle $\pi/8$ with the radius to a vertex. The right triangle gives
the inscribed disk radius $\cos(\pi/8)$, and the same calculation
works for $\Pi_N$.}\label{fig:region}
\end{figure}
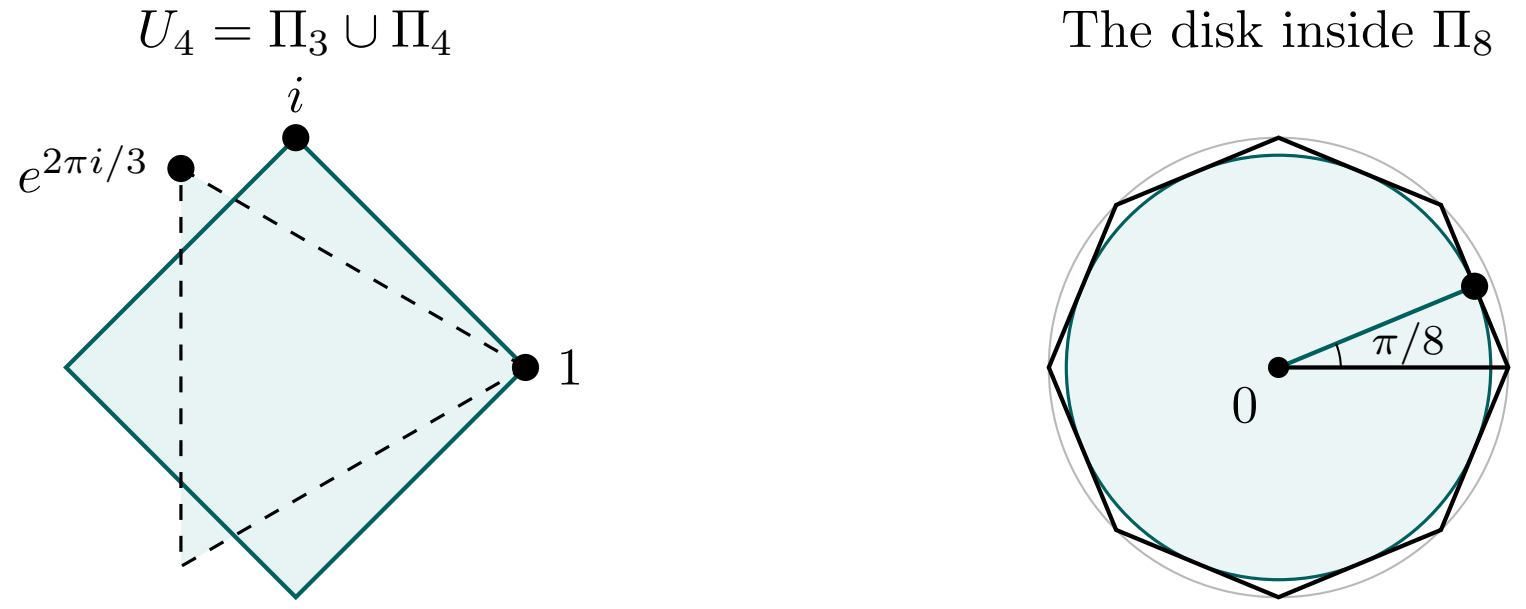

\paragraph{When do the polygons give the whole region?}
The inclusion $U_N\subseteq D_N$ goes back to an observation of
Kingman recorded by Mirsky~\cite{mirsky63}. Perfect and Mirsky gave
the cycle construction above, proved equality through order three and
asked whether equality holds in general~\cite{pm}. Order four already
requires more information than the ordinary coordinate hull supplies.
Levick, Pereira and Kribs also took the hull of the sums $v_i+v_j$
over distinct labels. Both polygons must be invariant, as the square
example illustrates. Together their constraints give $D_4=U_4$~\cite{four}.
At order five, Mashreghi and Rivard found a rational matrix with
an eigenvalue outside $U_5$~\cite{five}. The matrix at $t=1/2$ in
Appendix~\ref{app:loworders} gives their example with the labels in the
order $5,1,4,2,3$.
It is the average of two permutation matrices.

Is this a single exceptional point, or part of a curve?
Harlev, Johnson and Lim sampled the eigenvalues of $tP+(1-t)Q$,
with $P,Q$ permutation matrices, grouping together pairs that differ
only by a simultaneous relabelling~\cite{hjl}. Among their sampled
pairs at order five, one equivalence class produced eigenvalues outside $U_5$.
They conjectured that at most two permutation matrices suffice to attain
every boundary point of $D_N$. Their computations also led them to
conjecture $D_N=U_N$ for $6\le N\le11$. A mesh can reveal a curve.
To determine the region we must also exclude eigenvalues farther out,
including those of matrices using three or more permutations.

Two restrictions on the permutations explain part of what the plots
show. Jankowski, Johnson and Lim determined the region obtained from a
finite group of commuting matrices by taking all convex combinations of
its elements. Each root of unity that occurs as an eigenvalue in the
group supplies $\Pi_k$, where $k$ is the least positive integer for
which the $k$th power of the root is one. These polygons fill the
region~\cite{jjl}. For permutation matrices, each such root comes
from a cycle. In particular, if we take mixtures of powers of a single
permutation and then let that permutation vary, the region is exactly
$U_N$, as Casta\~neda Betancourt, Morneau-Gu\'erin and Singh also
showed directly~\cite{cms}. Johnson, Newman and Spitkovsky allowed permutations that
need not commute and prescribed their cycle lengths. They obtained an
outer bound by taking the convex hull of the corresponding polygons~\cite{jns}.
That hull fills gaps between polygons, precisely where an exception
like the one at order five can occur.

There is a different restriction one can put on the marks.
Bouthat, Mashreghi and Morneau-Gu\'erin announced that if no coordinate
lies in the convex hull of the remaining coordinates, then the
multiplier already lies in $\Pi_N$~\cite{bouthat}. An eigenvector
may also have repeated coordinates or marks strictly inside its hull.
We keep the mass at each of these sites. Theorem~\ref{large:main} uses
it to prove $D_N=U_N$ by hand for every $N\ge10^{30}$.
Appendix~\ref{sec:450} brings the threshold down to $450$ by computer.
Appendix~\ref{app:loworders} determines the regions at orders five and six.

\subsection{Three ways to use the column sums}
The square suggests a way to express its constraint without naming the
two competing rows. The function $f(w)=\max(0,\operatorname{Im}w)$
counts positive imaginary mass. In that example
$\sum_i f(v_i)=1$ and $\sum_i f(itv_i)=t(1+a)$. With column sums one, averaging each row can
only decrease the sum of a convex function. Could all such inequalities
together recover the column sums?

Give each label probability $1/N$. The resulting measure is
$\mu_v=N^{-1}\sum_i\delta_{v_i}$, where $\delta_w$ puts unit mass at $w$.
Thus $\int f\,d\mu_v=N^{-1}\sum_i f(v_i)$, with repeated sites counted
each time. For probability measures of finite support on the plane,
we say $\nu$ is below $\mu$ in \emph{convex order}, written $\nu\cx\mu$, if
$\int f\,d\nu\le\int f\,d\mu$ for every real-valued convex function
$f$. Multiplication by $z$ defines the map $T$. Moving each mass to
its image gives the \emph{push-forward} $T_\#\mu_v=N^{-1}\sum_i\delta_{zv_i}$.
If we keep the order of the sites,
the same assignment of masses gives a point in $\mathbb C^N$. The
natural convex set there is
\[
 \mathcal P(v)=\conv\{(v_{\sigma(1)},\dots,v_{\sigma(N)}):\sigma\in S_N\}.
\]
Here $S_N$ consists of all permutations of the $N$ labels, and $T$
acts on a vector by multiplying every entry by $z$.
The equivalence with convex order in the following proposition is the
equal-mass planar case of Sherman's theorem~\cite{sherman}. We give the
finite argument since it also explains what the column sums preserve.

\begin{proposition}\label{prop:assign}
Multiplication by $z$ leaves the marked configuration $v$ invariant
if and only if $T\mathcal P(v)\subseteq\mathcal P(v)$.
Equivalently, $T_\#\mu_v\cx\mu_v$.
\end{proposition}
\begin{proof}
A doubly stochastic matrix is a convex combination of permutations.
To check this, notice that any $k$ rows send mass $k$ into columns
of capacity one, and therefore meet at least $k$ columns. Choose as many
positive entries as possible with no row or column used twice. Such a
choice is a \emph{matching}. From an unmatched row, move to columns
along unchosen positive entries and back to rows along chosen entries.
Reaching an unmatched column would let us exchange chosen and unchosen
entries along the path and enlarge the matching. Otherwise each reached
column is paired with a reached row other than the starting row. The
reached rows then outnumber all their neighbouring columns, contrary to
the capacity bound. Hence the matching meets every row and gives a permutation.
Subtract the least of those $N$ entries along the permutation.
If any mass remains, renormalise and repeat. At each stage at least
one more entry vanishes. Conversely, a convex combination of permutations
has both sums one. Hence the vectors $Bv$, with $B$ doubly stochastic,
fill $\mathcal P(v)$. If $zv\in\mathcal P(v)$, permuting its coordinates
shows that every vertex image belongs to $\mathcal P(v)$. Convexity gives
the inclusion, and applying the inclusion to $v$ gives its converse.

If $y=Bv$, Jensen's inequality in each row gives
$\sum_i f(y_i)\le\sum_j f(v_j)$ for every convex $f$.
For the converse, suppose $y\notin\mathcal P(v)$. Let $p$ be its
nearest point in this polytope. The vector $y-p$ points outwards:
$(y-p)\cdot(x-p)\le0$ for every $x$ in the polytope, by minimizing
distance on the segment from $p$ to $x$. Its planar components $c_i$
therefore give a sum $\sum_i c_i\cdot y_i$ which exceeds
$\max_\sigma\sum_i c_i\cdot v_{\sigma(i)}$. Relabel the $v_j$ so that
the identity attains that maximum. Give the directed edge $i\to j$
weight $c_i\cdot(v_j-v_i)$. Every directed cycle has weight at most
zero, since a positive cycle would improve the assignment.
We want a convex function with a supporting plane of slope $c_i$
through $(v_i,b_i)$, where the heights $b_i$ are still to be chosen.
Such a plane stays below the height at $v_j$ provided
$b_j-b_i\ge c_i\cdot(v_j-v_i)$. The absence of a positive cycle
lets us satisfy all these inequalities at once. Let $b_j$ be the greatest weight of a path ending at $j$, allowing
an empty path. We may delete repeated cycles, so this is a finite
maximum. Appending an edge gives $b_j-b_i\ge c_i\cdot(v_j-v_i)$.
Set $a_i=c_i\cdot v_i-b_i$. The convex function
$f(w)=\max_i(c_i\cdot w-a_i)$ satisfies $f(v_j)=b_j$ for every $j$,
whereas $\sum_i f(y_i)\ge\sum_i(c_i\cdot y_i-a_i)>\sum_j b_j$.
Take $y=zv$.
\end{proof}

\paragraph{How to see the configuration in the plane.}
Choose complex coefficients $c_1,\ldots,c_N$. For each permutation
$\sigma$, plot the sum $\sum_i c_i v_{\sigma(i)}$ and take the convex
hull of the plotted points. Denote this polygon by $K_c(v)$. Its points
are exactly $\sum_i c_i(Bv)_i$ with $B$ doubly stochastic. If $Av=zv$,
multiplying any such sum by $z$ replaces $B$ by $BA$, so
$zK_c(v)\subseteq K_c(v)$. A single coefficient equal to one, with
all the others zero, recovers the ordinary coordinate polygon. Several
equal coefficients select sums of distinct coordinates, as in
Levick, Pereira and Kribs~\cite{four}. For example, coefficients one at the labels $1,a$ of
Figure~\ref{fig:square} and zero elsewhere give precisely the polygon
on the right of the figure. If we give two labels the coefficients $1$ and
$-1$, the points to plot are the differences $v_j-v_k$ with $j\ne k$.
The restriction on the two indices is exactly the column constraint.

Taken together, these planar polygons give an exact converse.
Indeed, if $zv$ lay outside $\mathcal P(v)$, a separating real linear
functional on $\mathbb C^N$ could be written as
$w\mapsto\Re\sum_i c_iw_i$. The point $\sum_i c_i v_i$ belongs to
$K_c(v)$, while its image under $z$ would lie strictly beyond a
supporting line of that polygon. Thus the familiar language of an
invariant polygon carries over once we allow all choices of $c$.

We can now replace the marks by another set of marks while preserving
the same multiplier. The analogue of cutting down an invariant polygon
is to choose a measure between the original measure and its image in convex order.
If $T_\#\mu_v\cx\nu\cx\mu_v$, then
$T_\#\nu\cx T_\#\mu_v\cx\nu$: composition of a convex function with
$T$ preserves convexity. Whenever $\nu\ne\delta_0$ consists of $N$ equal masses,
it gives another invariant configuration. The measure on the interior marks enters the two convex-order inequalities just as the mass on the boundary does.

\paragraph{What near equality gives us.}
Write $A$ as a weighted average of permutation matrices and choose a
random permutation $P$ with those same weights as probabilities.
Then $\mathbb E[v_{P(i)}]=zv_i$. In row $i$ the mean squared distance
to this barycentre is $\sum_j A_{ij}|v_j|^2-|zv_i|^2$.
On adding these quantities every $|v_j|^2$ in the first sum receives
total weight one. We obtain
\begin{equation}\label{eq:energy}
 \mathbb E\sum_i |v_{P(i)}-zv_i|^2
  =\sum_{ij} A_{ij}|v_j-zv_i|^2
  =(1-|z|^2)\sum_i|v_i|^2 .
\end{equation}
At $|z|=1$ each permutation of positive weight takes every $v_i$ to
$zv_i$. Following any nonzero mark around its cycle gives $z^j=1$
for some $j\le N$. How small is the budget if $z$ lies outside $U_N$?
The regular $N$-gon already contains the disk of radius $\cos(\pi/N)$
for $N\ge3$, as in Figure~\ref{fig:region}. Hence $1-|z|^2<\sin^2(\pi/N)\le\pi^2/N^2$.
A permutation can still move one site a long way, but that move has a cost.
Order the marks by distance from zero. If a narrow outer ring ends at
a large radial gap, a permutation crossing the gap pays for the crossing.
The permutations that stay in the ring have integer winding numbers.
Comparing their costs with the side between consecutive roots of unity
is the purpose of Sections~\ref{large:permutations} and~\ref{large:prefix}.
If no such gap occurs, many marks have to lie at intermediate radii.
Section~\ref{large:return} follows a randomly chosen mark through a
quarter turn. There are then too few marks in its thin annulus to
complete the turn at the energy allowed by the row averages.

A first reading can take Theorem~\ref{large:annulus} as given and go
directly to Section~\ref{large:return}. This proves the main theorem
from the ring estimate and explains why that estimate has precisely
the two alternatives in its statement. The proof of the estimate uses
finite Fourier series on a cycle. The path argument uses conditional
expectation for a finite Markov chain.

\input{large}
\clearpage
\appendix
\input{appendix450}
\clearpage
\input{smallorders}

\section*{Numerical verification}
Run \texttt{python verify450.py --show-endpoints} to reproduce the
last two columns of the endpoint tables with additional decimal places.
The verifier uses exact fractions where possible. For a quantity such
as a sine, it computes an interval containing the exact value. A strict
positive sign is certified only when the entire interval is above zero.
The README identifies the formulas and inputs checked by each module.
For the two counts in the worked example, add the option
\texttt{--explain-7r}. For Appendix~\ref{app:loworders}, the
electronic supplement contains the geometric arguments, the Python
programs for orders five and six, and the coefficient data they check.
From the \texttt{small-orders/} directory run
\texttt{python verify\_small\_orders.py 5} or replace $5$ by $6$.

\section*{Use of AI tools}
Discussions with ChatGPT Pro 6 helped the authors find inequalities
used in the pencil-and-paper proof for $N\ge10^{30}$.
GPT~6.2 Ultra (ChatGPT~6.2 Astra with ultra reasoning) wrote the
verification software and generated the certificates fully autonomously.
This includes the calculations for $N\ge450$ and the certificates for
orders five and six. The authors checked the software and the certificates.
OpenAI Codex also assisted with the exposition and assembly of the
electronic supplement. The authors take full responsibility for
the mathematical claims, proofs, computations, figures, and final text.

\end{document}

%% file: large.tex
\providecommand{\PME}{\mathbb E}
\providecommand{\PMP}{\mathbb P}
\providecommand{\PMav}[1]{\langle #1\rangle}
\providecommand{\PMarg}{\operatorname{Arg}}
\providecommand{\PMdist}{\operatorname{dist}}

\section{What two permutations tell us about an eigenvector}
\label{large:permutations}

Suppose that $Av=zv$ and that the coordinates of $v$ lie on a circle of positive radius about the origin. If a permutation of the coordinates multiplies each of them by the same number $\zeta$, then $\zeta^h=1$, where $h$ is the number of coordinates. Indeed, we can multiply the $h$ quotients and cancel all factors. We shall use this elementary argument even when the quotients differ. Their arguments still determine an integer. The difficulty is to compare the average of several such permutations with the segment between two consecutive roots of unity.

Where should we look for a circle when the moduli of the coordinates differ? We take the outermost coordinates first. A set $H$ of $h$ labels is a prefix if $|v_i|\ge |v_j|$ whenever $i\in H$ and $j\notin H$. Thus we can obtain a prefix by listing the labels in decreasing order of modulus and taking the first $h$. A gap after this prefix has a concrete consequence. If a permutation sends a label across the gap, its image has smaller squared modulus than the missing outer coordinate. We shall estimate how often this can happen.

\begin{theorem}\label{large:annulus}
Let $A$ be a real doubly stochastic matrix and let $Av=zv$, $v\ne0$. Suppose that $H$ is a prefix of $h\ge 10^{16}$ labels. Scale $v$ so that $h^{-1}\sum_{i\in H}|v_i|^2=1$, and suppose that
\[
                 1-\tfrac1{200}\le |v_i|\le 1+\tfrac1{200}
                 \qquad(i\in H).
\]
If $H$ contains all labels, or if
\[
       \min_{i\in H}|v_i|^2-\max_{j\notin H}|v_j|^2\ge \frac{5000}{h},
\]
then $z\in\Pi_h\cup\Pi_{2\lfloor h/2\rfloor}$.
\end{theorem}

The even polygon contains both $1$ and $-1$. For odd $h$ we shall use its sides only to estimate the distance from $z$ to these two real vertices. All comparisons of permutations take place on $h$ labels. We first work with a general number $n$ of labels, in order to see exactly where the size of the prefix enters.

\subsection{The integer attached to a permutation}

Write $\PMav{f_i}=n^{-1}\sum_{i=1}^n f_i$. Fix an angle $\theta$ and assume that
\begin{equation}\label{large:pairhyp}
 n\ge100,\qquad \PMav{|v_i|^2}=1,\qquad
 .995\le |v_i|\le1.005,\qquad d=|\PMav{v_i}|\le .01 .
\end{equation}
Put $g=2\pi/n$ and $a=|e^{ig}-1|^2=4\sin^2(\pi/n)$. These are the angle between consecutive $n$th roots of unity and the squared length of their common polygon side. We also need the complex second moment $M_2=\PMav{v_i^2}$. For instance, if $v_i=x_i+iy_i$, its real part is the difference between the averages of $x_i^2$ and $y_i^2$, and its imaginary part is the average of $2x_i y_i$. For a permutation $P$ of these labels let
\[
 \tau_{iP(i)}=\PMarg\frac{v_{P(i)}}{e^{i\theta}v_i}.
\]
Here $\PMarg$ is the angle in $(-\pi,\pi]$. Thus $\tau_{iP(i)}$ is the signed angular error after the prescribed rotation. We call $P$ good when $|\tau_{iP(i)}|\le10^{-4}$ at every label. Since the product of the normalized phase quotients equals $e^{-in\theta}$, we have $\exp(i\sum_i\tau_{iP(i)})=e^{-in\theta}$. Thus
\[
             k_P=\frac{n\theta+\sum_i\tau_{iP(i)}}{2\pi}
\]
is an integer. We call $k_P$ the winding integer of the draw. It counts the total turns made by all labels when each angular displacement is measured as $\theta+\tau_{iP(i)}$. In particular, $gk_P-\theta$ is the mean angular error.

Take two good permutations $P,Q$ and write $k=k_Q-k_P$. To compare their images, put
\begin{align*}
 E&=\PMav{|v_{Q(i)}-v_{P(i)}|^2},&
 B&=\PMav{(v_{Q(i)}-v_{P(i)})^2},\\
 T&=\PMav{|v_{Q(i)}-e^{igk}v_{P(i)}|^2}.
\end{align*}
The sum $E$ is the squared distance between the two images, averaged over all labels. In $T$ we first rotate one image through the difference of the two mean angles. The ordinary square in $B$ carries the same phase information as the second moment $M_2$. It will appear again when we square the row equations of $A$.

For example, place $n\ge32000$ labels at the $n$th roots of unity, let $P$ fix every label, and let $Q$ shift each label to the next root. For $\theta=g/2$ both draws are good, with $k_P=0$ and $k_Q=1$. Here $E=a$, $M_2=B=0$ and $T=0$. Give $Q$ probability $u$ and $P$ probability $1-u$. The mean sends each coordinate $v_i$ to $(1-u+ue^{ig})v_i$. The corresponding eigenvalue runs along the polygon side, with
\[
          1-|1-u+ue^{ig}|^2=au(1-u).
\]
Two independent draws differ with probability $2u(1-u)$. Half their mean squared distance is therefore exactly the loss of squared radius on the side. This is the comparison we want for general permutations. In Lemma~\ref{large:pair} we subtract $a|k|$ so that this regular example has $D=0$. The term containing $d^2$ accounts for a possible nonzero mean of the coordinates. When the cost is this small, the two lowest nonconstant Fourier modes have the same edge energy. The winding of the cycle will distinguish them.

\begin{lemma}\label{large:pair}
Assume \eqref{large:pairhyp}. Let $e=1$ if $|k|=1$ and $E<(1+10^{-6})a$, and put $e=0$ otherwise. With
\[
                     D=E-a|k|+10ea d^2,\qquad M=10^7,
\]
we have
\begin{align}
 D&\ge0,& ea&\le1.02E,& a(k^2-|k|)&\le3D, \label{large:pair1}\\
 |B-eaM_2|&\le MD,&
 T&\le MD+50ea|M_2|^2.&& \label{large:pair2}
\end{align}
\end{lemma}

\begin{proof}
Let $x_i=\tau_{iQ(i)}-\tau_{iP(i)}$. Its absolute value is at most $2\cdot10^{-4}$ and its sum is $2\pi k$. For points of moduli $r_1,r_2$ and angular difference $x_i$, the squared chord length is $(r_1-r_2)^2+4r_1r_2\sin^2(x_i/2)$. The edge from $v_{P(i)}$ to $v_{Q(i)}$ therefore costs at least $4(.995)^2\sin^2(x_i/2)\ge .99x_i^2$. Cauchy--Schwarz gives
\begin{equation}\label{large:basicpair}
           E\ge .99g^2k^2\ge .99ak^2.
\end{equation}
Always $|B|\le E$ and $T\le 2E+2ak^2$, since $|e^{igk}-1|\le |k|\sqrt a$.
For $k=0$ we have $T=D=E$. If $|k|\ge2$, then \eqref{large:basicpair} gives $D\ge .49E$. If $|k|=1$ and $e=0$, then $E\le1\,000\,001D$. These inequalities imply both assertions for $e=0$.

Suppose now that $e=1$. Conjugating the coordinates and changing $\theta$ to $-\theta$ if necessary, we may take $k=1$. We follow the permutation of the target labels $P(i)\mapsto Q(i)$. An edge with source label $P(i)$ has angular increment $x_i$. Following any cycle brings us back to its starting site, so the increments on that cycle sum to $2\pi w$ for an integer $w$. Since $|x_i|<\pi$, the straight edge between the two sites has the same angular increment and avoids the origin. Thus $w$ counts the signed turns of the closed polygon about zero. The integers from all cycles sum to $k=1$.

Choose a cycle with $w>0$ and write $\ell$ for its length. Spreading an angular change $2\pi w$ over $\ell$ edges requires $\sum_{\rm cycle}x_i^2\ge(2\pi w)^2/\ell$. Its contribution to $E$ is consequently at least $.99(2\pi w)^2/(n\ell)$. Dividing by $g^2=(2\pi/n)^2$ gives
\[
             .99w^2\frac{n}{\ell}\le E/g^2<1+10^{-6}.
\]
It follows that $w=1$ and $q=\ell/n\ge .98$. A short cycle could complete a turn only by using longer edges than the hypothesis on $E$ permits.

\smallskip
\noindent\emph{How much can lie outside this cycle?}
Write $s_0=1-q$, let $E_o$ be the contribution of the other cycles to $E$, and put $V_o=n^{-1}\sum_{i\text{ outside}}|v_i|^2\le1.011s_0$. Index the vertices on the chosen cycle in their cyclic order. We write
\[
 v_j=\sum_{b\bmod\ell}c_b e^{2\pi ibj/\ell},\qquad
 \lambda_b=4\sin^2(\pi b/\ell),\qquad a_\ell=\lambda_1 .
\]
Here $c_b=\ell^{-1}\sum_{j=0}^{\ell-1}v_j e^{-2\pi ibj/\ell}$. Moving from $j$ to $j+1$ multiplies mode $b$ by $e^{2\pi ib/\ell}$. The edge difference therefore multiplies it by $e^{2\pi ib/\ell}-1$, whose squared modulus is $\lambda_b$. The constant mode $c_0$ is the centre of the cycle and contributes no edge energy. The smallest positive value, $a_\ell$, occurs twice. Mode $1$ goes once around the circle counterclockwise and mode $-1$ goes once around it clockwise. Their edge energies are identical. We use a prime on a sum over $b$ to omit $0,1,-1$. The normalization and the bound on the mean give
\[
 q\sum_b|c_b|^2=1-V_o,\qquad q|c_0|\le d+1.005s_0 .
\]
Averaging a product of two Fourier modes, with one conjugated, gives zero unless their indices agree. This follows by summing the geometric series around the cycle. Thus the mean square of the list is the sum of the $|c_b|^2$, the finite Parseval identity. On applying this identity to the edge differences we get $E=E_o+q\sum_b\lambda_b|c_b|^2$. Subtracting $a$ therefore gives the exact expression
\begin{align*}
 D={}&E_o+q\sum_b'(\lambda_b-a_\ell)|c_b|^2\\
     &+a_\ell(1-V_o-q|c_0|^2)-a+10ad^2.
\end{align*}
For $0<x\le\pi/98$ the logarithmic derivative of $\sin^2x$ with respect to $\log x$ is $2x\cot x\in[1.9,2]$. Integration from $\pi/n$ to $\pi/\ell$ gives
\[
 1+1.9s_0\le q^{-1.9}\le a_\ell/a\le q^{-2}\le1.042,
 \qquad a_\ell/a-1\le2.1s_0 .
\]
Also $q|c_0|^2\le2.07(s_0^2+d^2)$. For $b\ne0,\pm1$, $\lambda_b\ge 2a_\ell$ since $\ell\ge98$. The expression for $D$ now yields
\begin{equation}\label{large:fouriercost}
 D\ge E_o+.49\sum_b'\lambda_b|c_b|^2+.8as_0+7ad^2 .
\end{equation}
To check the coefficients, the primed terms retain at least $q/2\ge.49$ of their squared edge energy. The term $a_\ell-a$ contributes at least $1.9as_0$. Subtracting $a_\ell V_o$ and $a_\ell q|c_0|^2$ leaves the coefficient of $as_0$ at least $1.9-1.042(1.011+2.07\cdot.02)>.8$ and that of $ad^2$ at least $10-1.042\cdot2.07>7$. Thus $t=D/a$ satisfies $0\le t<10^{-6}+10d^2<.00101$.

\smallskip
\noindent\emph{Which of the two first modes follows the winding?}
What distinguishes the two modes which cost the same edge energy? If $v_j=\alpha e^{2\pi ij/\ell}+\beta e^{-2\pi ij/\ell}$, then averaging $v_j^2$ gives $2\alpha\beta$, since the terms of frequency $2$ and $-2$ average to zero. Thus a small second moment forces a small product. We still have to decide which factor can be small. Figure~\ref{large:figmodes} shows the two directions and a sum in which both occur. We shall first use the sign of the winding to bound $|c_1|$ from below. Then $M_2$ will give an upper bound for $|c_{-1}|$.

\begin{figure}[htbp]
\centering
\begin{tikzpicture}[font=\small, line join=round]
\begin{scope}[scale=1.10]
\draw[gray!55] (0:1.35) \foreach \t in {45,90,...,315} { -- (\t:1.35) } -- cycle;
\foreach \t in {0,45,...,315} \fill[black!65] (\t:1.35) circle (1.65pt);
\draw[teal!70!black,very thick,->,shorten <=2pt,shorten >=2pt] (0:1.35)--(45:1.35);
\draw[black!65,thick,dashed,->,shorten <=2pt,shorten >=2pt] (0:1.35)--(-45:1.35);
\node[right=3pt] at (0:1.35) {$1$};
\node[above right] at (45:1.35) {$\zeta$};
\node[below right] at (-45:1.35) {$\zeta^{-1}$};
\node[below=6pt] at (0,-1.35) {Opposite first modes};
\end{scope}
\begin{scope}[xshift=7.0cm,scale=1.10]
\draw[gray!55] (-1.85,0)--(1.85,0) (0,-1.12)--(0,1.12);
\draw[teal!70!black,thick] (1.8,0) \foreach \t in {45,90,...,315}
  { -- ({1.8*cos(\t)},{.9*sin(\t)}) } --cycle;
\foreach \t in {0,45,...,315}
  \fill[teal!70!black] ({1.8*cos(\t)},{.9*sin(\t)}) circle (1.65pt);
\draw[teal!70!black,very thick,->,shorten <=3pt,shorten >=3pt]
  (1.8,0)--({1.8*cos(45)},{.9*sin(45)});
\node[right] at (1.8,0) {$\frac43$};
\node[above right] at (0,.90) {$\frac23 i$};
\node[below=6pt] at (0,-1.35) {$v_j=\zeta^j+\frac13\zeta^{-j}$};
\end{scope}
\end{tikzpicture}
\caption{Take $\zeta=e^{2\pi i/8}$. On the left the first edges of $\zeta^j$ and $\zeta^{-j}$ start at the same point and have equal length. On the right their indicated sum has real and imaginary coordinates $\frac43\cos(2\pi j/8)$ and $\frac23\sin(2\pi j/8)$. Its mean square without conjugation is $2/3$.}
\label{large:figmodes}
\end{figure}
 By $\sin x\ge2x/\pi$ on $[0,\pi/2]$,
\[
 \sum_b'\lambda_b^{-1}\le\ell^2/4,\qquad
 \max_j\left|\sum_b'c_b e^{2\pi ibj/\ell}\right|
       \le \tfrac12\ell\sqrt{D/.49}<.15 .
\]
For the first inequality pair $b$ with $-b$ and use $\sum_{j\ge1}j^{-2}<2$. The last strict bound uses $\ell\sqrt a\le 2\pi$. We also have $|c_0|<.031$ and $\sum_b'|c_b|^2\le2.05t$, by \eqref{large:fouriercost}. If $|c_1|<1/4$, Parseval and $|v_j|\ge .995$ give $|c_{-1}|>.9$. At every vertex, and hence at every point of corresponding linearly parametrized edges, our polygon lies within $.031+1/4+.15$ of the clockwise polygon with vertices $c_{-1}e^{-2\pi ij/\ell}$. The latter has distance at least $.9\cos(\pi/\ell)$ from the origin. The straight homotopy avoids the origin, so their winding numbers would coincide. The chosen cycle has winding one. Hence $|c_1|\ge1/4$.

In the full expansion of the ordinary square we have
\[
 \frac1\ell\sum_{j=0}^{\ell-1} v_j^2
  =\sum_{b,t}c_bc_t\left(\frac1\ell\sum_{j=0}^{\ell-1}e^{2\pi i(b+t)j/\ell}\right)
  =\sum_b c_b c_{-b}.
\]
Only opposite frequencies survive this average. The sum of the squares from the remaining labels, divided by $n$, has modulus at most $V_o$. More precisely,
\[
 |M_2-2qc_1c_{-1}|\le V_o+q|c_0|^2+q\sum_b'|c_b|^2\le4t .
\]
Here $V_o\le1.264t$ and $q|c_0|^2\le2.07(.02s_0+d^2)\le.35t$ by \eqref{large:fouriercost}. Since $2q|c_1|\ge .49$, we obtain
\begin{equation}\label{large:reverse}
 |c_{-1}|^2\le(2/.98)^2(2|M_2|^2+32t^2)\le 9|M_2|^2+t .
\end{equation}

To evaluate $B$, square the edge differences before averaging. The contribution of the cycle is $q\sum_b\lambda_b c_b c_{-b}$. In $B-aM_2$ the modes $1,-1$ therefore cost at most $2.1as_0$, the constant mode at most $.35D$, and the other modes at most $D/.49$. For this last estimate use $|\lambda_b-a|\le\lambda_b$ and pair $b$ with $-b$. The complementary cycles cost at most $E_o+aV_o$. By \eqref{large:fouriercost} the total is less than $9D$ when $D>0$, and zero when $D=0$.

For $T$ we compare the next point of the cycle with the present point rotated by $g$. On mode $b$ their difference is $(e^{2\pi ib/\ell}-e^{ig})c_b e^{2\pi ibj/\ell}$. In $T$ the multiplier of $|c_b|^2$ is therefore $|e^{2\pi ib/\ell}-e^{ig}|^2$. For $b=1,-1,0$ it is at most $1.1as_0^2$, $5a$, $a$, respectively. Indeed $g^2/a\le1.001$ and $q\ge .98$. On the other modes it is at most $4\lambda_b$ by the triangle inequality. The complementary cycles cost at most $2E_o+2aV_o$. Combining these bounds with \eqref{large:fouriercost} and \eqref{large:reverse} gives $T\le25D+45a|M_2|^2$. Finally \eqref{large:basicpair} gives $a\le1.02E$. This proves \eqref{large:pair1} and \eqref{large:pair2}.
\end{proof}

\subsection{Averaging the rotations}

For any probability law supported on good permutations, let $P,Q$ be independent draws. Each image has a mean rotation $e^{igk_P}$. We need to compare its average with the average of each coordinate quotient. An estimate linear in the size of a coordinate error would be too large here. The product identity removes that linear term after we sum over the labels.

\begin{lemma}\label{large:log}
Under \eqref{large:pairhyp} put $b_G=\PME e^{igk_P}$ and $z_i^G=\PME[v_{P(i)}/v_i]$. The principal logarithms satisfy
\[
               \left|\PMav{\log(z_i^G/b_G)}\right|\le10\PME T.
\]
\end{lemma}
\begin{proof}
For each draw set
\[
 h_{Pi}=\log\frac{|v_{P(i)}|}{|v_i|}
       +i\left(\tau_{iP(i)}-\PMav{\tau_{jP(j)}}\right).
\]
The index $j$ is averaged in the last term. Then $\PMav{h_{Pi}}=0$, $|h_{Pi}|<.0103$ and $e^{h_{Pi}}=v_{P(i)}/(e^{igk_P}v_i)$. The real parts sum to zero because $P$ permutes the moduli. The imaginary parts sum to zero because we have subtracted their mean. To see which average occurs after exponentiation, divide $z_i^G$ by $b_G$. This gives
\[
     z_i^G/b_G=\sum_P \nu_P e^{h_{Pi}}.
\]
The complex weights
$\nu_P=\PMP(P)e^{igk_P}/b_G$ sum to one and $|\nu_P|\le1.01\PMP(P)$, because $|b_G|\ge\cos(10^{-4})$.

Fix $i$, put $\mu=\PME h_{Pi}$, $x_P=h_{Pi}-\mu$, and write $w=\sum_P\nu_P(e^{x_P}-1)$. Thus $z_i^G/b_G=e^\mu(1+w)$. Since $|x_P|<.021$, the estimate $|e^x-1-x|\le e^{.021}|x|^2/2$ and the bound on $|\nu_P|$ give
\[
 |w-\textstyle\sum_P\nu_P x_P|\le\PME|x_P|^2,\qquad
 |w|\le1.04\PME|x_P|<.022 .
\]
Both logarithms are taken near one, so $\log(z_i^G/b_G)=\mu+\log(1+w)$. Using $|\log(1+w)-w|\le |w|^2$ gives
\[
        \left|\log(z_i^G/b_G)-\sum_P\nu_P h_{Pi}\right|
                 \le 4\PME|x_P|^2 .
\]
The variance on the right equals $\tfrac12\PME|h_{Pi}-h_{Qi}|^2$. For $|u|<.021$, $|e^u-1|\ge.98|u|$. Applying this to $u=h_{Qi}-h_{Pi}$, and using $|v_{P(i)}|\ge .995$, gives $|h_{Pi}-h_{Qi}|\le1.04|v_{Q(i)}-e^{igk}v_{P(i)}|$. Averaging over $i$ cancels the linear term and proves the lemma.
\end{proof}

\section{An eigenvalue from a prefix of nearly equal moduli}
\label{large:prefix}

We now prove Theorem~\ref{large:annulus}. Write $\PMav{f_i}=h^{-1}\sum_{i\in H}f_i$ throughout this section, and put
\[
 \delta=.005,\quad \delta_0=10^{-4},\quad K=5000,\quad h_0=10^{16},
 \qquad g=2\pi/h,\quad a=4\sin^2(\pi/h).
\]
Let $c=\PMav{v_i}$ and $M_2=\PMav{v_i^2}$. Suppose, to obtain a contradiction, that $z=re^{i\theta}$ lies outside both polygons of the theorem. We have $r\le1$ by stochasticity. Real $z$ belong to the even polygon. The incircle of $\Pi_h$ gives $r>\cos(\pi/h)>.99$, and consequently
\[
                             x=1-r^2<10/h^2 .
\]

\subsection{The loss when a permutation leaves the prefix}

On a side of $\Pi_m$ with outward normal at angle $\alpha_0$, the radius in direction $\alpha$ is
\begin{equation}\label{large:radial}
       \rho_m(\alpha)=\frac{\cos(\pi/m)}{\cos(\alpha-\alpha_0)} .
\end{equation}
The absolute derivatives of $\log\rho_m$ and $\rho_m^2$ are bounded by $\tan(\pi/m)$ and $2\tan(\pi/m)$, respectively. Each formula holds up to the side endpoints, so the same bound applies between points on different sides. Put $d_\theta=\PMdist(\theta,\pi\mathbb Z)$, the angular distance to the nearest real direction. The even polygon of the theorem has a vertex in each real direction. Since $z$ lies outside it, \eqref{large:radial} gives
\[
  x\le -2\log r\le2\tan(\pi/(h-1))d_\theta\le 8d_\theta/h .
\]
Also $|1-z^2|^2=x^2+4r^2\sin^2\theta$. Since $|\sin\theta|\ge 2d_\theta/\pi$, we obtain
\begin{equation}\label{large:anglebound}
              |1-z^2|\ge d_\theta\ge hx/8,
              \qquad |1-z|\ge hx/16 .
\end{equation}

Use the decomposition in the proof of Proposition~\ref{prop:assign} to represent $A$ by a random permutation $P$, with $(Pv)_i=v_{P(i)}$, and draw $Q$ independently with the same law. Put
\[
 s=\frac1h\PME\#\{i\in H:P(i)\notin H\},\qquad
 \Gamma=1-\PME\PMav{|v_{P(i)}|^2}.
\]
Each target outside $H$ takes the place of one missing target in $H$. On a proper prefix the difference of their squared moduli is at least $K/h$. On the full set both quantities vanish. Thus
\begin{equation}\label{large:leak}
                     s/h\le\Gamma/K .
\end{equation}
Define $E,B$ by the same expressions as in Lemma~\ref{large:pair}, now averaging over $i\in H$ for every pair of draws, and write $E_P^z=\PMav{|v_{P(i)}-zv_i|^2}$. For a single label $i$, the two independent images have the same mean $zv_i$. The cross term in $\PME|v_{Q(i)}-v_{P(i)}|^2$ is therefore $-2|zv_i|^2$. With $P$ held fixed, the same expansion around $zv_i$ leaves the variance of $v_{Q(i)}$ and the squared distance $|v_{P(i)}-zv_i|^2$. We write $\PME[E\mid P]$ for this average over $Q$ alone. Averaging over the labels in $H$ gives
\begin{align}
 x&=\Gamma+\tfrac12\PME E,&
 \PME E_P^z&=x-\Gamma,&
 \PME[E\mid P]&=E_P^z+x-\Gamma . \label{large:variance}
\end{align}
The first and second moments of the coordinates satisfy
\begin{equation}\label{large:means}
 |(1-z)c|\le3s,\qquad
 |(1-z^2)M_2-\tfrac12\PME B|\le3s .
\end{equation}
For the first, replace each target outside $H$ by a missing target in $H$ in the sum of the coordinates. One replacement costs at most $2(1+\delta)$, as all coordinates have modulus at most $1+\delta$. For the second estimate take ordinary squares in the target sum. The cost per replacement is at most $2(1+\delta)^2$, and $\tfrac12\PME B=\PME\PMav{v_{P(i)}^2}-z^2 M_2$.

If $x=0$, \eqref{large:variance} and \eqref{large:leak} give $s=0$ and $E_P^z=0$ almost surely. A permutation of positive probability preserves $H$ and multiplies every coordinate there by $z$. Multiplying the coordinate quotients gives $z^h=1$, a contradiction. Hence $x>0$. We have $|c|,|M_2|\le1$ by normalization, and
\begin{equation}\label{large:centroid}
             |c|\le 3s/|1-z|\le48/K=.0096 .
\end{equation}
This supplies the bound on $d$ needed in Lemma~\ref{large:pair}.

\subsection{The permutations omitted from the angular estimate}

Call a draw bad if it sends a label of $H$ outside $H$, or if it preserves $H$ and $|\tau_{iP(i)}|>\delta_0$ for some $i\in H$. Write $p$ for its probability, $F$ for the event that at least one of the two draws is bad, and $U=\PME[\mathbf1_F E]$. The indicator $\mathbf1_F$ is one on $F$ and zero otherwise, so $U$ counts the contribution of such pairs to the full mean. If $i,j\in H$, then
\[
       |v_j-zv_i|^2\ge 2r|v_iv_j|(1-\cos\tau_{ij})\ge \tau_{ij}^2/3 .
\]
Here $|\tau_{ij}|\le\pi$ and $r>.99$. A bad draw which preserves $H$ therefore has $E_P^z\ge10^{-9}/h$. A draw failing to preserve $H$ has at least one target outside. From \eqref{large:variance}, $\PME[\mathbf1_{\{P\text{ bad}\}} E_P^z]\le U$, and this expectation is also at most $x$. It follows that
\begin{equation}\label{large:bad}
 p\le hs+10^9hU,\qquad
 p\le h^2\Gamma/K+10^9 h x\le 10/K+10^{10}/h<1/4 .
\end{equation}

For a good draw use its winding integer $k_P$ on $H$. For a bad draw set $k_P=k_*$, where $k_*$ is a nearest integer to $\theta/g$. This convention assigns every omitted permutation the same rotation. Denote expectation conditional on one good draw by $\PME_G$. The definition of $k_P$ and Cauchy--Schwarz give
$g^2|k_P-\theta/g|^2\le 3E_P^z$ on good draws. Hence
\[
 \PME_G|k_P-k_*|^2
   \le \frac{6x}{g^2(1-p)}+\tfrac12<4 .
\]
For $k=k_Q-k_P$, a pair of bad draws has $k=0$. The two possibilities for exactly one bad draw contribute at most $8p$ to $\PME[\mathbf1_F k^2]$. Recalling $a\le4\pi^2/h^2$, we obtain
\begin{equation}\label{large:badk}
 a\PME[\mathbf1_F |k|]\le a\PME[\mathbf1_F k^2]\le 320p/h^2,
 \qquad p/h^2\le s/h+10^9 U/h .
\end{equation}

\subsection{The squared radius available on a side of the polygon}

Let $j=\lfloor\PME k_P\rfloor$, $u=\PME k_P-j$, and put
\begin{equation}\label{large:chord}
      \eta=e^{igj}(1-u+ue^{ig}),\qquad b=\PME e^{igk_P}.
\end{equation}
The point $b$ is the mean of the rotations of the draws. To obtain $\eta$, replace their winding integers by a variable taking only the values $j,j+1$, with the same mean as before. The probabilities are $1-u,u$, so $\eta$ belongs to the side joining $e^{igj}$ and $e^{ig(j+1)}$. In general its angle is different from $\theta$. We first compare its squared radius with $r^2$ and then use the excess to bound this difference of angles.

Let $p_e$ be the probability that both draws are good and the indicator $e$ of Lemma~\ref{large:pair} is one. On $F^c$ let $D$ be the quantity in that lemma with $d=|c|$, and put
\[
 J=\PME[\mathbf1_{F^c}D]+U,\qquad A_c=10ap_e|c|^2 .
\]
Since $ap_e\le1.02\PME E\le3x$, the compensation for the mean satisfies
\begin{equation}\label{large:ac}
          A_c\le30x\frac{3s}{|1-z|}\le1440s/h .
\end{equation}
We used $|c|^2\le |c|$ and \eqref{large:means}. Define
\[
             L=|\eta|^2-r^2,\qquad
             W=a\bigl(\tfrac12\PME|k|-u(1-u)\bigr).
\]

We claim $W\ge0$. Let $Y$ be an integer random variable with mean $\mu=j+u$, and let $Y'$ be an independent copy. Given $Y$, the function $t\mapsto|Y-t|$ is convex and $Y'$ has mean $\mu$. Jensen's inequality gives $\PME|Y-Y'|\ge\PME|Y-\mu|$. To bound the latter average, draw the graph of $|y-\mu|$ and the line through its values at $j$ and $j+1$. Every integer is at one of these two points or outside the interval between them. At every integer the graph thus lies on or above the line. Its equation gives
\[
 \PME|Y-Y'|\ \ge\ \PME|Y-\mu|
       \ \ge\ \PME\{u+(1-2u)(Y-j)\}=2u(1-u).
\]
The two-point law at $j,j+1$ attains equality. More generally, for two atoms at positive distances $A_0,B_0$ from their mean, half the mean difference is $A_0B_0/(A_0+B_0)$. One can also read this quantity from the cumulative probabilities $F_t=\PMP(Y\le t)$, since $\frac12\PME|Y-Y'|=\sum_{t\in\mathbb Z}F_t(1-F_t)$ counts each threshold between the two independent values. Taking $Y=k_P$ proves the claim.

We can now compare the loss on the polygon side with the loss of the eigenvalue. On $F^c$ the pair energy is $D+a|k|-10ea|c|^2$. On $F$ its unconditional contribution is $U$. Thus $\PME E=J+a\PME|k|-A_c-a\PME[\mathbf1_F|k|]$. By \eqref{large:variance}, half of this expectation, together with $\Gamma$, is $x$. Subtracting $au(1-u)=1-|\eta|^2$ gives
\[
 L=\Gamma+\tfrac12 J+W-\tfrac12\bigl(A_c+a\PME[\mathbf1_F|k|]\bigr).
\]
By \eqref{large:ac} and \eqref{large:badk} the last half-sum is at most $880s/h+1.6\cdot10^{11}U/h\le\Gamma/2+J/6$. Consequently
\begin{equation}\label{large:margin}
               L\ge W+\Gamma/2+J/3,\qquad 320p/h^2\le L .
\end{equation}
For the latter inequality use $p/h^2\le2L/K+3\cdot10^9 L/h$ and $h\ge h_0$.

We also need the distance of $b$ from this same side point. For any integer $y$, linear interpolation of $f(y)=e^{igy}$ at $j,j+1$ gives an error at most $a(y-j)(y-j-1)/2$. Indeed, if $\Delta f(t)=f(t+1)-f(t)$, then $\Delta^2 f(t)=e^{igt}(e^{ig}-1)^2$ has modulus $a$. At $y=j+m$ with $m\ge2$, the error of the linear interpolant is
\[
       \sum_{t=0}^{m-2}(m-1-t)\Delta^2 f(j+t),
\]
whose coefficients sum to $m(m-1)/2$. Summing from the other endpoint gives the same bound for $y<j$. Taking expectations of the interpolant gives $\eta$. Since $\operatorname{Var}(k_P)=\PME k^2/2$, the error bound gives
\[
 |b-\eta|\le W/2+(a/4)\PME(k^2-|k|)\le 3L.
\]
Indeed the good pairs contribute at most $3J$ to $a\PME(k^2-|k|)$ by \eqref{large:pair1}, and \eqref{large:badk} and \eqref{large:margin} bound the remainder by $L$.

\subsection{Comparing the two arguments}

The second equation in \eqref{large:means} is the reason for keeping the ordinary square $B$. On good pairs we can subtract its leading term $eaM_2$ by \eqref{large:pair2}. The other pairs have $|B|\le E$. With $M=10^7$ we obtain
\[
   \big|\big(1-z^2-ap_e/2\big)M_2\big|\le MJ+3s .
\]
By \eqref{large:anglebound}, $|1-z^2|\ge hx/8$, whereas $ap_e/2\le 3x/2$. The reverse triangle inequality bounds $|1-z^2-ap_e/2|$ below by $hx/8-3x/2\ge hx/16$. Divide by this bound, then multiply by $ap_e|M_2|\le3x$. We obtain
\[
           ap_e|M_2|^2\le(48/h)(MJ+3s).
\]
The factor $h^{-1}$ comes from the even polygon. Use this in the estimate for $T$ in \eqref{large:pair2}. With \eqref{large:leak} and \eqref{large:margin} it gives
\begin{equation}\label{large:mt}
 \PME[\mathbf1_{F^c}T]\le M(1+2400/h)J+7200s/h\le 7ML .
\end{equation}

Write $b_G=\PME_Ge^{igk_P}$ and $z_i^G=\PME_G[v_{P(i)}/v_i]$. Conditional on $F^c$ the two draws remain independent with their common conditional law. Since $(1-p)^2>1/2$, Lemma~\ref{large:log} and \eqref{large:mt} imply
\[
            |\PMav{\log(z_i^G/b_G)}|\le 140ML .
\]
After rotating by $e^{-i\theta}$, the real parts of $z_i^G,b_G,b,z$ all exceed $.99$. For the first of these use $(199/201)\cos\delta_0>.99$, for the two rotation means use $|gk_P-\theta|\le\delta_0$ on good draws and $|gk_*-\theta|\le\pi/h$ on the others. The principal logarithm in the rotated right half-plane is $2$-Lipschitz on these points, by integrating its derivative along a line segment. In particular differences of their arguments may be taken in one branch. The eigenvector equations on the omitted draws give
\[
 \PMav{|z_i^G-z|}
 \le \frac{\PME[\mathbf1_{\{P\text{ bad}\}}\PMav{|v_{P(i)}-zv_i|}]}{(1-p)(1-\delta)}
 \le 2\sqrt{pU}.
\]
Here we used Cauchy--Schwarz and $\PME[\mathbf1_{\{P\text{ bad}\}}E_P^z]\le U$. Moreover $|b-b_G|\le pg\PME_G|k_P-k_*|\le13p/h$. For two complex numbers of moduli at least $.99$, their circular angular distance is at most $2$ times their Euclidean distance, by the chord inequality. Apply this also to $b$ and $\eta$, since $|\eta|\ge\cos(\pi/h)$. Passing from $z$ to the mean argument of the $z_i^G$, then to $b_G,b,\eta$, yields
\[
       \PMdist(\arg z,\arg\eta)\le 150ML+26p/h+4\sqrt{pU}.
\]
Here and below the distance between arguments is taken modulo $2\pi$.

Recall that $\eta$ belongs to the boundary. To conclude, it suffices to show that changing its argument to $\theta$ can reduce the squared radius of that boundary by at most $L=|\eta|^2-r^2$. By \eqref{large:radial} the squared radial function of $\Pi_h$ has Lipschitz constant at most $8/h$. Hence the loss in squared radius between $\eta$ and the direction of $z$ is at most
\[
                      (1200M/h)L+208p/h^2+32\sqrt{pU}/h .
\]
We keep the two sources of the last term separate. By \eqref{large:badk} and \eqref{large:leak},
\[
 \frac{\sqrt{pU}}h\le\sqrt{\Gamma U/K}+\sqrt{10^9/h}\,U .
\]
From \eqref{large:margin}, $\Gamma/2+U/3\le L$ gives $\sqrt{\Gamma U}\le\sqrt6\,L/2$. Using the same inequalities for $p/h^2$, we bound the loss by
\[
 \left(\frac{1200M}{h}+\frac{416}{K}+\frac{624\cdot10^9}{h}
        +16\sqrt{6/K}+96\sqrt{10^9/h}\right)L .
\]
At $K=5000$ and $h\ge10^{16}$ the coefficient is less than $.67$. In particular $\rho_h(\theta)^2\ge|\eta|^2-L=r^2$, also when $L=0$. We have proved $z\in\Pi_h$, contradicting our assumption and completing the proof of Theorem~\ref{large:annulus}.

\section{An angular return and a bound on the order}
\label{large:return}

The prefix theorem applies when the moduli are nearly equal and a sufficiently large gap follows the last label. What if there is no such gap? Then many distinct labels lie in each of a number of thin annuli. We shall follow a stationary path through a quarter turn. With high probability it stays in the annulus in which it starts. To turn through this angle using only the labels in that annulus costs a definite amount of squared edge length. The eigenvalue equations give an upper bound for the same cost.

\begin{figure}[htbp]
\centering
\begin{tikzpicture}[font=\small,line join=round]
\def\visited{0/0,1/1,2/2,3/3,4/5,5/7,6/9}
\pgfmathsetmacro{\modulus}{cos(5)}
\begin{scope}
  \draw[gray!60] (0,0)--(3.55,0) (0,0)--(0,3.50);
  \draw[gray!50] (0:3.14) arc(0:90:3.14);
  \foreach \a in {0,10,...,90} \fill[gray!65] (\a:3.14) circle (1.4pt);
  \foreach \k/\j in \visited {\coordinate (v\k) at ({10*\j}:3.14);}
  \draw[black!50,dotted,thick] (v0)--(v6);
  \draw[gray!55] (30:1.30)--(v3) (45:1.30)--(45:{3.14*\modulus^3});
  \draw[gray!70] (30:1.30) arc (30:45:1.30);
  \node[fill=white,inner sep=1pt] at (37.5:1.01) {$15^\circ$};
  \draw[teal!75!black,very thick] (v0)--(v1)--(v2)--(v3)--(v4)--(v5)--(v6);
  \foreach \k/\j in \visited \fill[teal!75!black] (v\k) circle (2.05pt);
  \coordinate (prediction) at (45:{3.14*\modulus^3});
  \draw[gray,dashed,thick] (v3)--(prediction);
  \draw[black!70,fill=white,thick] (prediction) circle (2pt);
  \draw[black!50] (prediction)--(2.13,2.65)--(2.52,2.88);
  \node[right,inner sep=1.5pt] at (2.52,2.88) {$\lambda^3v_0$};
  \node[below left] at (0,0) {$0$};
  \node[below] at (v0) {$k=0$};
  \node[right=4pt] at (v3) {$k=3$};
  \node[left=2pt] at (v6) {$k=6$};
\end{scope}
\begin{scope}[xshift=6.45cm,yshift=.05cm,x=.74cm,y=10cm]
  \draw[->,gray!70] (0,0)--(6.62,0) node[below right,black] {$k$};
  \draw[->,gray!70] (0,0)--(0,.316);
  \node[above right,black] at (0,.317) {$|\lambda^{-k}v_{X_k}-v_0|$};
  \foreach \y in {.1,.2,.3} {
    \draw[gray!23] (0,\y)--(6,\y);
    \node[left,inner sep=2pt] at (0,\y) {$\y$};
  }
  \foreach \k/\j in \visited {
    \pgfmathsetmacro{\a}{\modulus^(-\k)}
    \pgfmathsetmacro{\err}{sqrt((\a-1)^2+4*\a*sin((10*\j-15*\k)/2)^2)}
    \coordinate (e\k) at (\k,\err);
    \draw[gray!65] (\k,0)--(\k,-.004);
    \node[below,inner sep=2pt] at (\k,0) {$\k$};
  }
  \draw[teal!75!black,thick] (e0)--(e1)--(e2)--(e3)--(e4)--(e5)--(e6);
  \foreach \k/\j in \visited \fill[teal!75!black] (e\k) circle[radius=1.9pt];
\end{scope}
\end{tikzpicture}
\caption{The visited labels are $0,1,2,3,5,7,9$. The dark path joins their coordinates. The pale points mark the remaining sites in this quadrant, and the dotted chord joins the two endpoints. The open point gives the expected position after three steps. The right panel uses the same seven labels to plot their distances from the starting point after division by $\lambda^k$.}
\label{large:figpath}
\end{figure}
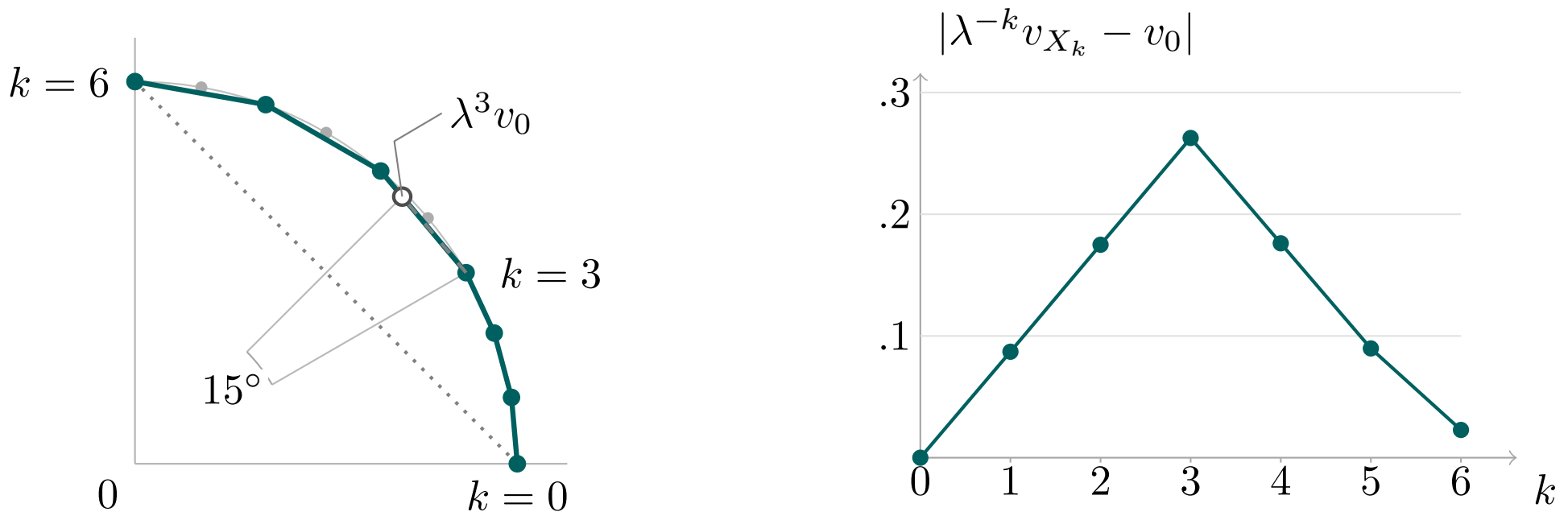

\paragraph{One path on a circle.}
To see which part of the motion to measure, place $36$ states at $v_j=\zeta^j$ with $\zeta=e^{2\pi i/36}$ and $j$ taken modulo $36$. From $j$ move to $j+1$ or $j+2$ with probability $1/2$ each. If $X_k$ denotes the state after $k$ moves, the mean next coordinate given $X_k=j$ is
\[
   \tfrac12(\zeta^{j+1}+\zeta^{j+2})=\lambda v_j,\qquad
   \lambda=R e^{i\pi/12},\quad R=\cos(\pi/36).
\]
After six moves the expected coordinate has turned through a right angle. Consider the path that starts at $0$, takes three single steps and then three double steps. At its third move its coordinate has argument $30^\circ$, while $\lambda^3v_0$ has argument $45^\circ$. At its sixth move these arguments agree again. We can undo the expected rotation and contraction at time $k$ by dividing the visited coordinate by $\lambda^k$. Its distance from the starting point is then $|\lambda^{-k}v_{X_k}-v_0|$. Figure~\ref{large:figpath} gives this distance at each of the seven times. The largest value along a path can be much bigger than the value at its endpoint. We shall therefore need a bound for the maximum over all times.

\begin{theorem}\label{large:main}
For every integer $N\ge 10^{30}$ we have $D_N=U_N$.
\end{theorem}

\begin{proof}
We have already realised every point of $U_N$ as an eigenvalue by taking powers of a cyclic permutation. It remains to prove $D_N\subseteq U_N$.

Take $\delta=1/200$, $K=5000$ and $h_0=10^{16}$ as in Theorem~\ref{large:annulus}. Set
\begin{equation}\label{large:constants}
 \gamma=\min\{\delta/(1+\delta),1/200\},\quad
 C=1.52\max\{1,K\},\quad J=\lceil 40C/\gamma\rceil,\quad w=.999\gamma/J .
\end{equation}
The width $w$ will be used for an annulus and $J$ counts how many of them fit in the chosen range of moduli. We prove the inclusion whenever
\begin{equation}\label{large:cutoff}
                 N\ge\max\left\{4h_0,\frac{4250C}{\gamma w^2}\right\}.
\end{equation}
In particular $N\ge10^6$, $N\ge 50000/\gamma$, $N\ge4000K$, $N\ge1000C/w$ and $N\ge20000/w$. These inequalities follow directly from \eqref{large:constants} and $0<w<\gamma\le .005$.

\smallskip
\noindent\emph{Rotating through a small angle.}
Suppose $Av=zv$ with $z=re^{i\theta}\notin U_N$, and multiply $v$ by a scalar so that $\max_i|v_i|=1$ and one coordinate is $1$. On the circle the $N+1$ angles $0,\theta,\ldots,N\theta$ include a pair at distance at most $2\pi/(N+1)$. Let $d$, $1\le d\le N$, be the difference of their indices and let $\varphi=\PMdist(d\theta,2\pi\mathbb Z)$. Put $m=d\lfloor N/d\rfloor$. Then $N/2<m\le N$ and every $d$th root of unity is a vertex of $\Pi_m$. If $\varphi=0$, the radial segment through $z$ ends at such a vertex and $z\in\Pi_m$. Hence $0<\varphi\le2\pi/N$. Conjugating if necessary, we can take $B=A^d$ and $\lambda=z^d=Re^{i\varphi}$ with $Bv=\lambda v$.

There is a vertex of $\Pi_m$ at angular distance $\varphi/d$ from $e^{i\theta}$. Applying \eqref{large:radial} between that vertex and the direction of $z$, and using $z\notin\Pi_m$, gives
\begin{align}
 -\log R&\le 6.29\varphi/N,\nonumber\\
 1-R\cos\varphi&\le 6.29\varphi/N+\varphi^2/2\le9.44\varphi/N,
 &R\sin\varphi&\ge .999\varphi .\label{large:returnbounds}
\end{align}
Indeed $\tan(\pi/m)\le\tan(2\pi/N)\le6.29/N$ for $N\ge10^6$, and $\sin\varphi/\varphi\ge1-\varphi^2/6$. The row and column sums of $B$ are one. Upon expanding the squares, we obtain the energy
\begin{equation}\label{large:returnenergy}
 \mathcal E=\sum_{i,j}B_{ij}|v_j-v_i|^2
       =2(1-R\cos\varphi)\sum_i|v_i|^2\le18.88\varphi .
\end{equation}

\smallskip
\noindent\emph{How many coordinates have modulus near one?}
Write $\mathcal H=\operatorname{conv}\{v_i\}$. Since $\lambda\mathcal H\subseteq\mathcal H$ and $1\in\mathcal H$, the hull contains $\lambda^k$ for $0\le k\le\lfloor2\pi/\varphi\rfloor$. The angular gaps between these points are at most $\varphi$, and by \eqref{large:returnbounds} their moduli are at least $e^{-12.58\pi/N}$. In every direction the support function of $\mathcal H$ is therefore at least $e^{-12.58\pi/N}\cos(\varphi/2)\ge1-40/N$. In the last inequality one can use $e^{-t}\ge1-t$, $\cos t\ge1-t^2/2$ and $N\ge10^6$. To see the disk inclusion, take any unit direction. Some $\lambda^k$ lies within angle $\varphi/2$ of it, so the orthogonal projection of that point on the direction is at least $1-40/N$. Thus every supporting line of $\mathcal H$ is at least $1-40/N$ from zero, with the origin on the inner side. The closed disk of this radius lies in $\mathcal H$, and every hull vertex has at least this modulus.

We also need to count the vertices. For comparison, a regular $b$-gon on the unit circle has inradius $\cos(\pi/b)\ge1-\pi^2/(2b^2)$. Taking $b$ of order $\sqrt N$ already gives an inscribed disk of radius $1-40/N$ if the multiplicative constant is large enough. We want a number of vertices proportional to $N$. The small angle of $\lambda$ supplies a second condition on the same polygon. Let $u_1,\ldots,u_b$ be its distinct vertices. Every $\lambda u_i$ is a convex combination of them. The weights of these combinations supply a stochastic matrix $S$ satisfying $Su=\lambda u$. A stationary probability vector $\pi$ satisfies $\pi^TS=\pi^T$. To obtain one, average the first $m$ iterates of any initial probability row under $S$ and take a convergent subsequence. Multiplying the average by $S$ changes only its two end terms, each with weight $1/m$, so the limit is stationary. The weights $\pi_i S_{ij}$ on ordered pairs have equal inflow and outflow at each label. They are a sum of nonnegative multiples of directed simple cycles. Indeed, follow positive edges until a label repeats, subtract the smallest weight on the cycle, and repeat. At least one edge disappears at each subtraction.

On a cyclic list $u_{i_0},\ldots,u_{i_{l-1}}$ write $y_j=u_{i_j}$ and $y_l=y_0$. We have
\[
   \sum_j\operatorname{Im}(\overline y_j y_{j+1})
       \le \tfrac12\cot(\pi/b)\sum_j|y_{j+1}-y_j|^2 .
\]
To verify this, take the discrete Fourier expansion of the list. A nonconstant mode of frequency $k$, $1\le k<l$, gives the ratio
\[
 \frac{\sin(2\pi k/l)}{4\sin^2(\pi k/l)}=\tfrac12\cot(\pi k/l)\le\tfrac12\cot(\pi/b),
\]
since $l\le b$. Constant modes and loops contribute zero. Summing over the cycle decomposition, the left side becomes $R\sin\varphi\sum_i\pi_i|u_i|^2$. The total squared edge length on the right is $2(1-R\cos\varphi)\sum_i\pi_i|u_i|^2$ by stationarity. Every vertex has positive modulus, so the common factor cancels. In view of \eqref{large:returnbounds},
\[
  .999\varphi\le R\sin\varphi
       \le\cot(\pi/b)(1-R\cos\varphi)\le\frac{b}{\pi}\frac{9.44\varphi}{N},
       \qquad b>.332N .
\]
Each of these vertices occurs at a distinct coordinate label.

\smallskip
\noindent\emph{The gaps between consecutive moduli.}
Order the coordinate moduli as $1=b_1\ge\cdots\ge b_N$. In a prefix with $b_h\ge1-\gamma$, the root mean square lies in $[1-\gamma,1]$. Dividing by it puts all the prefix moduli in $[1-\gamma,(1-\gamma)^{-1}]\subseteq[1-\delta,1+\delta]$. Thus if $b_N\ge1-\gamma$, the full set satisfies Theorem~\ref{large:annulus}, which contradicts $z\notin U_N$. In the remaining case the range of the moduli covers
\[
                          I=[1-\gamma,1-.001\gamma] .
\]
Every hull vertex has modulus above $I$, because $40/N<.001\gamma$. Consider an open gap $(b_{h+1},b_h)$ meeting the interior of $I$. We have $b_h\ge1-\gamma$ and $b_{h+1}$ is smaller than every hull-vertex modulus. The upper prefix contains all the hull vertices, so $h>.332N\ge h_0$. Normalizing it can only increase its squared gap, as the root mean square is at most one. Theorem~\ref{large:annulus} would now place $z$ in $U_N$ if this squared gap reached $K/h$. Hence
\[
              b_h-b_{h+1}\le \frac{K}{h(b_h+b_{h+1})}.
\]
First $b_h+b_{h+1}\ge .995$ gives $b_h-b_{h+1}\le 3.03K/N<.001$. Substitution back gives $b_h+b_{h+1}\ge2(1-\gamma)-.001\ge1.989$. Since $1/(.332\cdot1.989)<1.52$, every such gap has width at most $C/N$. It follows that a subinterval of $I$ of length $t$ contains at least $tN/C-1$ labels. Indeed, $k$ labels divide it into at most $k+1$ positive gaps or parts of gaps, each of length at most $C/N$.

\smallskip
\noindent\emph{Following a path through a quarter turn.}
Partition $I$ into $J$ intervals of length $w$, assigning a shared endpoint to the interval on its right and keeping the final right endpoint. Let $n_j$ count the labels in interval $j$. Its closed central interval of length $.6w$ has at least $.6wN/C-1$ labels. A uniformly chosen label thus starts in this central interval with probability at least
\[
                           m_0=.599w/C
\]
because $N\ge 1000C/w$.

Start at a uniformly chosen label $X_0$. From label $i$ choose label $j$ with probability $B_{ij}$, independently of the past given $i$. This defines the Markov chain $(X_k)_{k\ge0}$. Its distribution stays uniform since $\sum_i B_{ij}/N=1/N$ for every $j$. Put $T=\lfloor\pi/(2\varphi)\rfloor$. We have $\pi/2-\varphi<T\varphi\le\pi/2$ and $1-R^T\le6.29\pi/(2N)<10/N$. As in Figure~\ref{large:figpath}, we divide out the expected factor at each time. Given the path through time $k$, the next coordinate has conditional mean $(Bv)_{X_k}=\lambda v_{X_k}$. Consequently $Y_k=\lambda^{-k}v_{X_k}$ has conditional mean $Y_{k-1}$ given the path through time $k-1$. This is the martingale property. In particular $\PME[Y_T\mid X_0]=Y_0$. By stationarity $v_{X_T}$ and $v_{X_0}$ have the same second moment, so
\[
 \PME|Y_T-Y_0|^2
 =(R^{-2T}-1)\PME|Y_0|^2
 \le e^{6.29\pi/N}-1\le 20/N .
\]
Indeed $\PME(\overline{Y_0}Y_T)=\PME|Y_0|^2$ by conditioning on $X_0$. For $|Y_k-Y_0|^2$, the conditional mean of the next value is at least the present value by convexity of the squared modulus. This gives a bound uniform in time
\begin{equation}\label{large:maxpath}
 \PMP\left(\max_{k\le T}|Y_k-Y_0|>.199w\right)
 \le\frac{20}{.199^2Nw^2}\le\frac{m_0J}{5}.
\end{equation}
For the first inequality, put $Z_k=|Y_k-Y_0|^2$. On the event that the first crossing occurs at time $j$, conditional Jensen gives $\PME[Z_T\mid X_0,\ldots,X_j]\ge Z_j>(.199w)^2$. Summing the disjoint first-crossing events bounds their total probability by $\PME Z_T/(.199w)^2$. The second follows from \eqref{large:cutoff}, $m_0J=.599\cdot.999\gamma/C$ and $20/(.199^2\cdot4250)<.599\cdot.999/5$.
The same path also has a bound on its total squared edge length. At every time $k$ the probability of the directed edge $i\to j$ is $B_{ij}/N$, since $X_k$ is uniform. Its expected cost is therefore $\mathcal E/N$, independently of $k$. Summing over the $T$ steps and applying \eqref{large:returnenergy} gives
\begin{equation}\label{large:pathupper}
            \PME\sum_{k<T}|v_{X_{k+1}}-v_{X_k}|^2
                =T\mathcal E/N\le18.88\pi/(2N)<29.7/N .
\end{equation}

Suppose a path starts in the central part of interval $j$ and $\max_{k\le T}|Y_k-Y_0|\le.199w$. For every $k\le T$,
\[
 \big||v_{X_k}|-|v_{X_0}|\big|\le .199w+1-R^T<.2w
\]
by $N\ge20000/w$. Thus all the visited labels have moduli strictly inside that same interval. The first and last coordinates, however, are far apart. Since $\cos(T\varphi)\le\varphi$,
\[
 |1-\lambda^T|^2\ge2-(20+4\pi)/N>1.414^2,
 \quad |v_{X_T}-v_{X_0}|\ge1.414(1-\gamma)-.199w>\sqrt{1.97}.
\]
Delete each closed piece between repeated labels until the path has no repetition. Every retained edge occurred in the original path. Its endpoints are still $X_0,X_T$, which are distinct by the last bound. The shortened path visits each label at most once and therefore has $s\le n_j-1$ edges. If their complex differences are $d_1,\ldots,d_s$, they add to the same endpoint difference. Hence
\[
  1.97 < |v_{X_T}-v_{X_0}|^2
       =\left|\sum_{\ell=1}^s d_\ell\right|^2
       \le s\sum_{\ell=1}^s|d_\ell|^2
       \le n_j\sum_{k<T}|v_{X_{k+1}}-v_{X_k}|^2.
\]
This gives the promised cost $1.97/n_j$ for a successful path in interval $j$.

Let $f_j$ be the probability of starting in the central part of interval $j$ and violating the bound $.199w$. The paths just considered have probability at least $q_j=\max\{m_0-f_j,0\}$ in interval $j$. Since their starting events are disjoint, \eqref{large:maxpath} implies
\[
           0\le q_j\le m_0,\qquad \sum_jq_j\ge 4m_0J/5 .
\]
The intervals can contain different numbers $n_j$ of labels. A larger count permits cheaper paths in that interval, but takes labels away from the others. Since $\sum_j n_j\le N$, Cauchy--Schwarz applied to $\sqrt{q_j/n_j}$ and $\sqrt{n_j}$ gives a lower bound which uses only the total number of labels. Each $n_j$ is positive by the central interval count. The expected path cost is at least
\begin{align*}
 1.97\sum_j\frac{q_j}{n_j}
 &\ge\frac{1.97}{N}\left(\sum_j\sqrt{q_j}\right)^2
 \ge\frac{1.97(4/5)^2 m_0J^2}{N}\\
 &\ge \frac{1.97(4/5)^2(.599)(.999)(40)}{N}>\frac{30}{N}.
\end{align*}
Here $\sqrt{q_j}\ge q_j/\sqrt{m_0}$ and $\gamma J/C\ge40$. This contradicts \eqref{large:pathupper} and proves the inclusion under \eqref{large:cutoff}.

\smallskip
\noindent\emph{The numerical cutoff.}
Substitution in \eqref{large:constants} gives
\[
 \gamma=1/201,\qquad C=7600,\qquad J=61\,104\,000,
 \qquad w=\frac{999}{1000\cdot201\cdot61\,104\,000}.
\]
To check $10^{30}$ by hand, use $4250C<3.24\cdot10^7$,
$201^3<8.2\cdot10^6$, and $J/.999<6.12\cdot10^7$. Since
$6.12^2<37.5$ and $3.24\cdot8.2\cdot3.75=99.63$, the second
entry in \eqref{large:cutoff} is less than $99.63\cdot10^{28}$.
The first entry is $4\cdot10^{16}$. For comparison, the exact ceiling of the second is
\[
 \left\lceil\frac{4250\cdot7600\cdot201^3\cdot(61\,104\,000)^2}
                         {(999/1000)^2}\right\rceil
 =981\,293\,676\,405\,877\,185\,293\,401\,509\,618<10^{30}.
\]
This establishes Theorem~\ref{large:main}.
\end{proof}

%% file: appendix450.tex
\section{The bound from order 450}\label{sec:450}
\begin{theorem}\label{main}
For every integer $N\ge450$, one has $D_N=U_N$.
\end{theorem}
The annulus argument in the main text allows a radial gap of order
$5000/h$ at the $h$th mark. At $h=450$ this permits a gap larger
than the entire range of squared radii. To retain useful information
at these orders we have to look at which gaps a permutation can cross
together. Order the marks from the outside in. A permutation preserves
the cut at $h$ when it permutes the first $h$ labels among themselves.
If it also preserves the cut at $h+s$, it permutes the $s$ labels
between the cuts. On this smaller block every permutation eigenvalue
is a root of unity of order at most $s$. A rotation close to a root
of larger order then costs squared distance. This is the source of the
estimates for several consecutive cuts.

The reduced fractions in $[0,1]$ with denominator at most $N$, listed
in increasing order, form the Farey sequence of order $N$. A fraction
$f$ gives the vertex $e^{2\pi i f}$. Together these vertices are exactly
the points of $U_N$ on the unit circle. The direction of an exterior
eigenvalue lies strictly between two consecutive vertex directions. Let $d$ be the smaller denominator. We group consecutive cuts into packets of $d$, allowing
a shorter last packet. For $d\ge18$, adding the bounds on successive radial gaps
gives less than the squared span the marks would have to cover.
For $2\le d\le17$, the same estimates give a count of points in
chosen radial intervals. We use those points to construct a convex
function $f$ with $\sum_i f(z^d v_i)>\sum_i f(v_i)$. This would
contradict Proposition~\ref{prop:assign} applied to $A^d$.

A row of a table describes a range of possible configurations. We
follow the same marks through two operations. Ordering them by radius
tells us how many lie between two circles. Raising $z$ to a small
integer power tells us how far those marks turn. For example, in row $7$R the direction of $z$ is near a primitive
seventh root, so $z^7$ returns near $1$. The matrix can have \emph{any}
order $N\ge450$.
The label $7$ records a denominator. The matrix order stays $N$.
The entries come from inequalities that hold throughout each parameter
interval. When the number of marks occurs in an error term such as
$1/N$, using $1/450$ bounds that error from above in every order under
consideration. For example, (L5) below spells out this substitution for
a count of marks. The choice $450$ is a sufficient cutoff.
A decimal used as a bound denotes the indicated
rational number. The two columns showing the endpoints $p_-,p_+$
in the final tables are for orientation. Their defining formulas,
before rounding, determine the intervals checked in the proof.
\begin{center}
\begin{tabular}{@{}lll@{}}
\toprule
Quantity & Meaning & Where it enters\\
\midrule
$N,\ n=450$ & actual order and its fixed lower bound & all estimates\\
$d\le e$ & denominators of the two Farey endpoints & the return $z^d$\\
$h,\ \tau=h/N$ & number and fraction of outer labels & gaps at a radial cut\\
$u_h$ & probability a permutation preserves that cut & cost of a packet\\
$\mu,E_p$ in (L5) & cost per cut and terminal excess & counts of interior marks\\
$\Delta(p)$ & lower bound for the total convex gain & final contradiction\\
\bottomrule
\end{tabular}
\end{center}

\subsection{The first mark met by a supporting plane}
The configuration has a planar test that can be drawn with a straightedge. Give the mark $v_i=r_ie^{i\theta_i}$ height $F(r_i)$ for a convex nondecreasing function $F$. Look at the lifted marks from the side. At $v_i\ne0$, choose a supporting slope $s_i\ge0$ for the radial graph. Thus the line of this slope through $(r_i,F(r_i))$ lies below $F$ at every radius. Hold the radial slope fixed and increase the slope in the counterclockwise tangent direction until it meets a mark in the forward half-plane. Let $T_i$ be this last slope. Writing $\alpha_{ij}=\theta_j-\theta_i$ modulo $2\pi$, one has
\[
 T_i=\min_{0<\alpha_{ij}<\pi,\ r_j>0}
 \left\{s_i\tan(\alpha_{ij}/2)+
 \frac{F(r_j)-F(r_i)-s_i(r_j-r_i)}{r_j\sin\alpha_{ij}}\right\}.
\]
The second term measures the clearance between the radial supporting line and the other height. This term vanishes when both radii lie on a linear piece of the radial graph with slope $s_i$. It can keep a mark at a different radius from blocking the tilt.

\begin{lemma}\label{tilt}
If $Av=wv$ for a doubly stochastic $A$ and $\operatorname{Im}w>0$, then
\[
 \sum_{r_i>0}r_iT_i\le\frac{1-\operatorname{Re}w}{\operatorname{Im}w}
                         \sum_{r_i>0}r_i s_i .
\]
\end{lemma}
\begin{proof}
The imaginary part of row $i$, after rotation by $-\theta_i$, gives a destination in the forward half-plane, so the minimum exists. Write $g_i=(s_i+iT_i)e^{i\theta_i}$. The plane condition is $g_i\cdot(v_j-v_i)\le F(r_j)-F(r_i)$; solving it for forward $j$ gives the formula. At every other mark radial convexity gives the inequality, since its tangential projection is nonpositive. At zero take $g_i=0$. Average over destinations with probabilities $A_{ij}$ and sum over sources. Every height occurs with equal total weight on the two sides. The remaining inequality is the displayed one.
\end{proof}

For the paraboloid $F(r)=r^2/2$, the plane of gradient $c$ through the lifted mark at $v_i$ lies below every lifted mark precisely when $|c-v_i|\le |c-v_j|$ for all $j$. Thus the tilt can be seen in the original plane. Start with the centre $c=v_i$ and slide it on the counterclockwise tangent ray. The disk with centre $c$ through $v_i$ grows until a mark in the forward half-plane reaches its boundary. At that instant $T_i=|c-v_i|$.
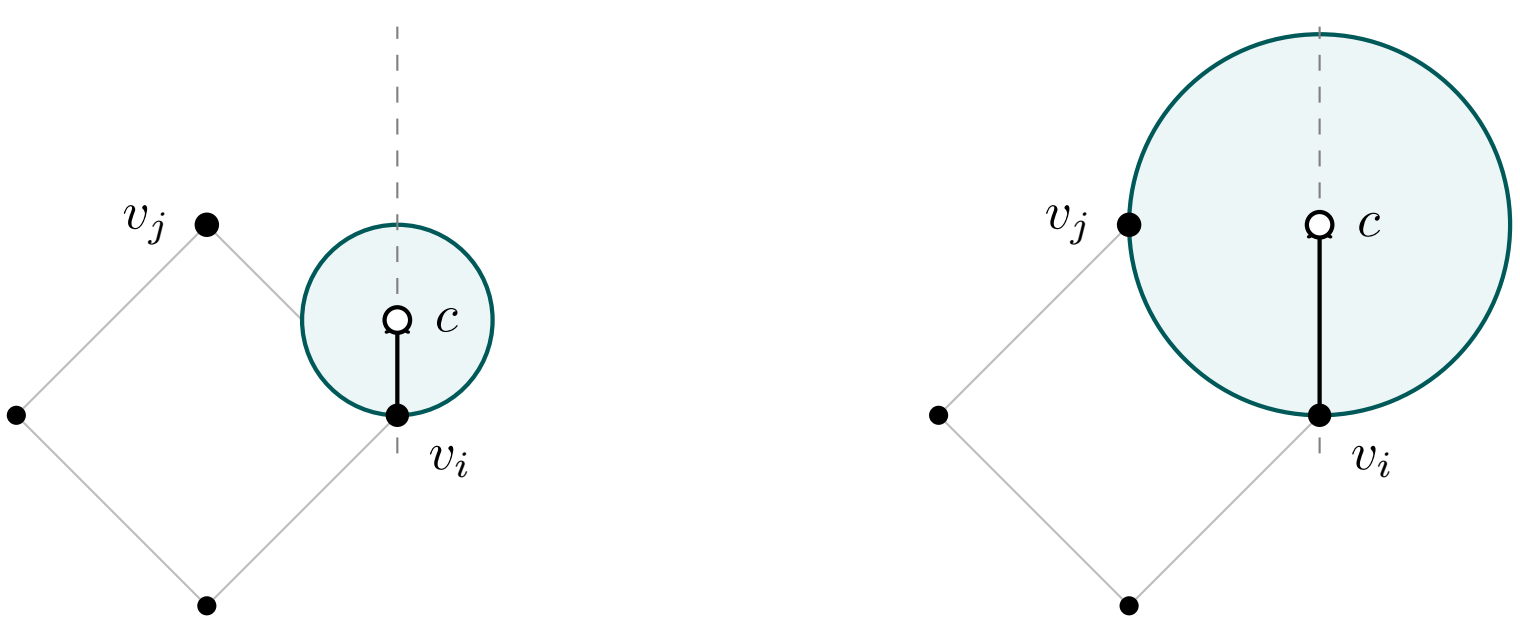
\begin{figure}[htbp]
\centering
\begin{tikzpicture}[font=\small,x=1.26cm,y=1.26cm]
\foreach \shiftx/\t/\head in {0/.5/{Before contact},6.1/1/{First contact}}{
 \begin{scope}[xshift=\shiftx cm]
 \node at (.45,2.49) {\head, $t=\pgfmathprintnumber{\t}$};
 \draw[thin,gray!52] (0,1)--(1,0)--(0,-1)--(-1,0)--cycle;
 \filldraw[draw=teal!70!black,fill=teal!7,thick] (1,\t) circle[radius=\t];
 \draw[dashed,gray] (1,-.2)--(1,2.04);
 \draw[thick,->] (1,0)--(1,\t);
 \foreach \pt in {(-1,0),(0,-1)} \fill \pt circle (1.8pt);
 \fill (0,1) circle (2.3pt) node[left=3pt] {$v_j$};
 \fill (1,0) circle (2.2pt) node[below right=2pt] {$v_i$};
 \draw[fill=white,thick] (1,\t) circle (2.4pt) node[right=3pt] {$c$};
 \end{scope}}
\end{tikzpicture}
\caption{Slide $c=(1,t)$ up the tangent at $v_i=(1,0)$, always
drawing the disk through $v_i$. For the four vertices of this square,
$v_j=(0,1)$ is the first forward mark on its boundary.
Indeed $|c-v_j|^2-|c-v_i|^2=2-2t$. Here the largest tangential
slope is exactly $T_i=1$.}\label{fig:first-contact}
\end{figure}

For instance, suppose $Av=wv$ with $\operatorname{Im}w>0$, and suppose its $m$ coordinates are distinct points on a circle of positive radius about zero. Take $F(r)=r$. The first contact is the next mark in angular order. All angular gaps $g_j$ are less than $\pi$, by the forward-neighbor condition. Jensen's inequality yields $m^{-1}\sum_j\tan(g_j/2)\ge\tan(\pi/m)$, since $\sum g_j=2\pi$. This recovers the familiar polygon sector inequality directly from the tilted planes. In the general proof we use a height with several slopes and average the position of its outer breakpoint. The first contact can then detect the radial distribution as well as the circular order.

\subsection{Farey returns and radial ranks}\label{sec:rank}
Suppose $N\ge n=450$, $Av=zv$, $A$ is doubly stochastic,
and $z=re^{i\theta}$ lies outside $U_N=\bigcup_{j=1}^N\Pi_j$,
where $\Pi_j$ is the convex hull of the $j$th roots of unity.
Normalize and label the coordinates so that $1=r_1\ge\cdots\ge r_N$. The radial rank $h$ is the position in this list. The cut at $h<N$ separates its first $h$ labels from the others.
Write $x=1-r^2$, $X=N^2x$, and let $b$ be the number of vertices
of their coordinate hull, counting distinct vertices once.
Write $a_h=h(r_h^2-r_{h+1}^2)$, $w_h=N(r_h^2-r_{h+1}^2)$ and
$\tau=h/N$. For the top $h$ coordinates put
$V=\max_{i\le h}r_i^2/[h^{-1}\sum_{i\le h}r_i^2]$.
We first estimate the angular return and its effect on the coordinate hull. Section~\ref{sec:prefix} uses these geometric estimates to exclude the full prefix and bound the individual gap at suitable radii. The following rank sections sum those gaps. The unnormalized mean square of a prefix is at most one, so a normalized squared gap smaller than $K/h$ will imply $a_h<K$.
\paragraph{One interval on the circle.}
For order eight the fractions $1/3$ and $3/8$ are neighbours.
Their difference is $1/24$. At a fraction $t$ of the way from
$1/3$ towards $3/8$, the angle in turns is $1/3+t/24$.
Three steps advance by $1+t/8$ turns. The full turn disappears and
the small angle $2\pi t/8$ remains. If two cuts at positions $h$ and
$h+2$ are preserved by the same permutation, its restriction to the two
labels between them has eigenvalues among $1$ and $-1$. Neither is the
third root $e^{2\pi i/3}$. The proof below estimates this mismatch
uniformly as the denominators and the order vary.

\paragraph{The two rational endpoints.}
The argument $\theta/(2\pi)$ lies between consecutive fractions of
Farey order $N$, since at a root-of-unity angle the unit disk on that
ray lies in the corresponding polygon (at angle zero use $\Pi_2$).
Denote their denominators by $d\le e$. Their separation is $1/(de)$
and $d+e>N$. To see both assertions, start with $0/1,1/1$. If $p/d<q/e$ satisfy $dq-ep=1$, then any reduced fraction $k/h$ strictly between them has $qh-ek>0$ and $dk-ph>0$. Expanding $(h,k)=(qh-ek)(d,p)+(dk-ph)(e,q)$ gives $h\ge d+e$, with equality only for the mediant $(p+q)/(d+e)$. Inserting that mediant when the allowed order reaches $d+e$ preserves determinant one on both sides. Induction proves the assertions. Put $y=e/N$, $s=d/N$, $m=s\lfloor 1/s\rfloor$.
Both $y,m$ exceed $1/2$. Let $t\in(0,1)$ be the fraction of the
interval measured from the endpoint of denominator $d$, so
$\phi=2\pi t/e$ is the absolute angular return of $z^d$.
The case $d=1$ is excluded by the stochastic sector inequality at $1$.
Give the directed edge $i\to j$ weight $A_{ij}/N$. At every vertex as much weight enters as leaves. Follow positive edges until a vertex repeats, subtract the smallest edge weight on that cycle from all its edges, and repeat. Each subtraction removes an edge and preserves balance. This writes the flow as a nonnegative sum of cycles. On each cycle
of length $j\le N$, expand the visited coordinates $w_k$ in Fourier modes as in Lemma~\ref{large:pair}. We obtain
\[
 \sum_k\operatorname{Im}(\bar w_k w_{k+1})
 \le \tfrac12\cot(\pi/N)\sum_k|w_{k+1}-w_k|^2.
\]
The ratio for mode $a\ne0$ is $\frac12\cot(\pi a/j)$.
For the uniform stationary distribution of $A$ the inequality gives
$|\operatorname{Im}z|\le\cot(\pi/N)(1-\operatorname{Re}z)$.
For $|\arg z|\le2\pi/N$ this places $z$ in $\Pi_N$.

A regular polygon of order $j$ has logarithmic radial slope at most
$\tan(\pi/j)$ in absolute value. Integrate from either endpoint, using
orders $Nm$ and $e$. Since $\tan u/u<\xi=1.00020$ for
$0\le u\le2\pi/n$, it follows that
\[
 \frac{N^3(-\log r)}{2\pi^2}
 \le\frac{\xi}{sy}\min\left(\frac{t}{m},\frac{1-t}{y}\right).
 \qquad dX\le\frac{4\pi^2\xi}{(m+y)y}<79.                 \tag{R1}
\]
The second bound maximizes over $t$ and uses $1-r^2\le2(-\log r)$.
For $d\ge10$ this gives $r>.99998$. At $2\le d\le9$ the bound (R3) below gives the same conclusion. The distance to the nearest integer turn
from the endpoint of denominator $d$ is at least $1/d$. For $d\ge3$
the same holds after doubling that reduced endpoint. The chord inequality
$\sin(\pi a)\ge2a$ for $0\le a\le1/2$ therefore gives
\[
 |1-z|\ge2\sqrt r\sin(\pi(1-1/e)/d)>3.99/d,\qquad
 |1-z^2|\ge2r\sin(\pi(1-2/e)/d)>3.8/d\quad(d\ge3).       \tag{R2}
\]
For the first strict bound at $d=2$ use $e>N-2$; at $d\ge3$ use the chord $\sin(\pi a)\ge(3\sqrt3/2)a$ up to $a=1/3$. Here $e>N/2$ and $hx/|1-z|<20/N$ follows at every rank.

\paragraph{A disk common to all denominators.}
For $d/N\le1/2$ choose the power $D=d$, and write
$z^D=Re^{\pm i\varphi}$. In this range $m\ge2/3$, so (R1) gives
$N(-\log R)/(\pi\varphi)\le3\xi/2$. For $s>1/2$ choose $D=e-d$.
The signed returns at the two endpoints are opposite, so
$N\varphi/(2\pi)=t/y+(1-t)/s$. In (R1) one may replace $m$ by $s$.
The ratio
\[
 \frac{(y-s)\min(t/s,(1-t)/y)}{st+y(1-t)}
 \le \frac{y-s}{s^2+y^2}<1
\]
increases up to $t=s/(s+y)$ and decreases thereafter. Hence the
same bound on $N(-\log R)/(\pi\varphi)$ holds. Both choices have
$\varphi<4\pi/N$. The iterates of any unit-modulus coordinate,
through $\lfloor2\pi/\varphi\rfloor$, are in the coordinate hull.
Their arguments have gaps at most $\varphi$, so the hull contains
a disk of radius at least
$\exp[-3\pi^2\xi/N]\cos(2\pi/N)>1-30/N$.
For example $1-e^{-u}\cos v\le u+v^2/2$ proves the last inequality.
Every hull vertex is therefore at radius greater than $1-30/N$.

We shall also need a vertex count when $s>1/4$. Write
$F=N\varphi/(2\pi)$, $E_0=N(-\log R)/(\pi\xi\varphi)$.
For $1/4\le s\le1/2$ the same return $D=d$ gives
$E_0\le\min(1/m,(1-t)/(yt))$, $F=t/y$.
The sum on the first branch increases in $t$ and the sum on the second
is convex. Its supremum is at most
$\max(1/y,1/m+m/[y(m+y)])$.
For $s\le1/3$ use $m\ge3/4,y\ge2/3$ and $m\le1$
to obtain $4/3+.9<7/3$. For $1/3<s\le1/2$, $m=2s$ and
$y>1-s$. The second term is bounded by
$1/(2s)+2s/(1-s^2)$, a convex function with endpoint values
$9/4,7/3$. In addition $(m+y)y\ge .75$ in both cases.
For the second case this follows from $(1+s)(1-s)\ge3/4$.

For $s>1/2$, take $D=e-d$ and $v=y/s\in(1,2)$.
The bound for $E_0$ above increases up to $t=s/(s+y)$ and then
decreases. When $F$ is added the first expression is convex in
$t$ and the second decreases. Its largest value, including both
endpoints, is at most
\[
 \frac1s\left[1+\frac{(v-1)^2}{v(1+v)(1+v^2)}\right]\le13/6 .
\]
Here $v-1\le v/2$ and $v-1\le (v+1)/3$.
In all these ranges $E_0+F\le7/3$. The cycle inequality on the
hull vertices (proved just below) and $R\sin\varphi>.99945\varphi$
now give $b/N>.99945/[(7/3)\xi]>.428$. For $s>1/2$ the sharper bound $13/6$ likewise gives $b/N>.461$.

\paragraph{The return at smaller denominators.}
For $d/N\le1/4$ take $D=d$, $R=r^d$, $\varphi=2\pi t/e$, and
let $u$ be any lower bound for $m,y$, such as $u=1-(d-1)/N$.
We use only choices $u\ge.75$. At distance $a\in[0,2\pi/j]$
from a vertex of $\Pi_j$, its logarithmic radial deficit is at most
$\frac12\sec^2(\pi/j)a(2\pi/j-a)$. This follows by integrating
the second derivative with its endpoint values zero. In our range
$\sec^2(\pi/j)<\xi$. Applied at both endpoints this gives
\[
 \frac{N(-\log R)}{\pi\varphi}
 \le\frac{\xi}{u}\min\left(1-\frac{ut}{d},
       \frac{1-t}{t}\left(1-\frac{1-t}{d}\right)\right).
\]
Set $v_0=\min(t,1-t)$ and $F=N\varphi/(2\pi)\le t/u$.
After adding $t$ to the first or second expression in parentheses,
the first increases to $t=1/2$ and the second decreases from that point.
The second derivative test is unnecessary since the derivative
of the second sum is $(1-1/d)(1-1/t^2)$. Using also the same
deficit bound in its original symmetric form yields
\begin{align}
 \frac{N(-\log R)}{\pi\varphi}+F
 &\le B_d=\xi\left(\frac{3}{2u}-\frac1{2d}\right),\nonumber\\
 dX&\le\frac{4\pi^2\xi}{u^2}v_0(1-u v_0/d)
 \le\frac{2\pi^2\xi}{u^2}(1-u/(2d)),\qquad
 -\log R\le\frac{\pi^2\xi}{u^2N^2}.                      \tag{R3}
\end{align}
For the second display, the deficit from the $d$ endpoint is bounded
by $4\pi^2\xi t(1-ut/d)/u^2$ in units of $dX$; the other endpoint
has the same bound with $t$ replaced by $1-t$.
The same proof applies to $2dN^2(-\log r)$, giving the bound on
$-\log R$. At $d=2$ we have $1-R=x$ and
$|1-z^2|\ge R\sin\varphi>.99945\varphi>.315Nx$.

Apply the stationary cycle inequality on the $b$ hull vertices to
this return. The matrix on the vertices is obtained by expressing
each return image as a convex combination of vertices. Any stationary
distribution has positive squared mass, since all vertices have positive
radius. Hence $R\sin\varphi\le\cot(\pi/b)(1-R\cos\varphi)$.
From $R\sin\varphi>.99945\varphi$ and (R3) we get
$b/N>.99945/B_d$. For the radius estimate use the power $d$. The $d$ initial
angles deviate to the same side of a regular $d$-grid by less
than $\varphi=2\pi t/e<2\pi/d$. Their order is preserved.
Starting from each of these $d$ points, take at most
$\lceil2\pi/(d\varphi)+1\rceil$ further rotations by $z^d$.
The original gaps are at most $2\pi/d+\varphi$, so the resulting
gaps are at most $\varphi$. Since $r^j=R_d^{j/d}$ for $j<d$,
the disk supported by these points has radius at least $R_d^{2\pi/(d\varphi)+3}\cos(\varphi/2)$,
where $R_d=r^d$. This argument places no upper bound on $d$.
Write the disk radius as $e^{-L_d}$. The logarithmic versions of
(R1) and (R3), together with $-\log\cos(\varphi/2)\le\varphi^2/4$,
give
\[
 L_d\le \frac{2\pi^2\xi}{u dN}
           +\frac{(3\xi+1)\pi^2}{u^2N^2}
 \quad (d/N\le.25,\ u=1-(d-1)/N),                       \tag{L1}
\]
and at all denominators
\[
 L_d\le \frac{4\pi^2\xi}{dN}
               +\frac{118.5+4\pi^2}{N^2}.
\]
For the latter inequality use $N(-\log R_d)/(\pi\varphi)\le
\xi/m\le2\xi$, $-\log R_d<39.5/N^2$, and $\varphi<4\pi/N$.
Thus $d/N>.25$ gives a hull radius greater than $1-317/N^2$.

A slightly shorter mesh will be useful. For each direction take the preceding point of the first $d$ powers, then advance by the nearest nonnegative integer number of $z^d$ rotations. The angular error is at most $\varphi/2$ and the total exponent, in units of $d$, is at most $2\pi/(d\varphi)+5/2$. Thus the disk has radius at least $R^{2\pi/(d\varphi)+5/2}\cos(\varphi/2)$. From (R3), $N(-\log R)/(\pi\varphi)\le\xi/u$ and $N^2(-\log R)\le\pi^2\xi/u^2$. Since $1-e^{-x}\cos y\le x+y^2/2$, the deficit times $dN$ is bounded by
\[
 2\pi^2\xi/u+(5\xi/2+1/2)\pi^2d/(Nu^2)<22\quad (d\le18).
\]
The last expression increases with $d$ and decreases with $N$.
Write $u_j=1-(j-1)/n$ and put
\[
 b_j=10^{-3}\left\lfloor\frac{999.45}{\xi(1.5/u_j-1/(2j))}\right\rfloor,
 \quad H_j=1-22/(jn),\quad
 X_j=\frac{2(9.86961)\xi}{j u_j^2}(1-u_j/(2j)).
\]
Up to $d=36$, $X_d$ decreases and $B_d$ increases. For the first assertion the derivative after multiplication by $d^2u^2/(2(9.86961)\xi)$ is $-1+2v+(u/d)(1-v/2)<0$, with $v=d/(n+1-d)<.087$. In $19\le d\le36$ (R3), (L1) give at least $.616N$ hull vertices above $.9968$ and $X<1.108$. Through $d=18$ we have $dX<22$, and through $36$ we have $dX<23.5$, giving the respective bounds $5.6,5.9$ on $Nhx/|1-z|$ by (R2). The two returns considered in the hull count satisfy $R\sin\varphi>.99945\varphi$, since $-\log R<39.5/N^2$ and $\varphi<4\pi/N$.

At $37\le d\le N/4$ take $u=.75$ in the hull count to get $b/N>.499$. In (L1) take $u=.90$ for $d\le N/10$, then $d>N/10$ and $u=.75$ otherwise. This gives hull radius greater than $.998$. The function $d(1-d/N)^2$ increases up to $N/4$, hence $X<2\pi^2\xi/[37(1-37/450)^2]<.641$. Equation (R3) gives $dX<35.1$, while (R2), $e>225$ and $\sin t\ge t-t^3/6$ give $d|1-z|>6.17$. Thus $Nhx/|1-z|<5.7$. Above $N/4$ the radius exceeds $1-317/N^2>.9984$, the hull count exceeds $.428N$ and (R1) gives $X<.470$. Here $(m+y)y\ge.75$ for $s\le.5$ and $(m+y)y>.5$ thereafter. Also $dX<79$ and $Nhx/|1-z|<13$.

\subsection{The outermost points}\label{sec:prefix}
A prefix consists of the first $h$ labels after sorting by decreasing modulus, with an arbitrary fixed order on ties. Temporarily divide the entire configuration by the root mean square on this prefix. In the next calculation a good permutation preserves the prefix and all its changes in angle lie near the argument of $z$. Every other permutation is called bad. A leaking permutation sends some prefix label outside it. Its squared radial loss will appear in the same estimate that controls movement across a polygon side.
Let $Av=zv$ with $A$ doubly stochastic. Normalize an outer prefix of $h$ labels by $\langle |v|^2\rangle=1$, where brackets are the prefix average, and suppose all its radii are at least $\alpha>0$. Let $P$ be a random permutation with mean $A$. Partition the draws into a good class, all of which preserve the prefix, and a bad class of probability $p$; put $q=1-p>0$. Let $v_P=(v_{P(i)})_{i\le h}$ be the vector of images on these labels. The conditional mean $\mu_G=\mathbb E_G v_P$ averages only the good draws, dividing their original weights by $q$. Write $\mu_B=\mathbb E_B v_P$ for the bad draws when $p>0$. On these vectors use the squared norm $\|w\|^2=\langle |w|^2\rangle$. Put
\[
 z_G=\langle\bar v\mu_G\rangle,\quad H=\|\mu_G-z_Gv\|^2,\quad
 V_G=\mathbb E_G\|v_P-\mu_G\|^2,\quad Z=\langle |v|^{-2}\rangle .
\]
Since $\|v\|=1$, $z_Gv$ is the orthogonal projection of $\mu_G$ onto the complex line through $v$. Thus $H$ measures how far the conditional mean is from any common scalar multiple of the marking. Define the radial loss $\Gamma=1-\mathbb E\|v_P\|^2$. The good images have norm one, so $\Gamma=p(1-\mathbb E_B\|v_P\|^2)$. We use the convention that terms from an empty bad class are zero.

For any side normal $e^{i\psi}$ of $\Pi_h$ put $Y=\mathbb E[\mathbf1_{\{P\ {\rm bad}\}}\|v_P-e^{i\psi}v\|^2]$. Projecting the row equation $zv=q\mu_G+p\mu_B$ onto $v^\perp$ gives
\[
 q^2 H\le pY .
\]
Writing $z_B=\langle\bar v\mu_B\rangle$, one also has the exact supporting-line identity
\[
 2p(\operatorname{Re}(e^{-i\psi}z_B)-\cos(\pi/h))
 =2p(1-\cos(\pi/h))-\Gamma-Y.
\]
In particular the loss at leaking draws is already measured in the direction of the very side one is trying to prove.

\begin{lemma}\label{pay}
Suppose $|z_G|>0$ and $\max_i|\mu_{G,i}/v_i-z_G|\le\epsilon<|z_G|$. Set $c_h=\cos(\pi/h)$, $S_h=2\tan(\pi/h)$,
\[
 C_1=\frac1{2\alpha^2(|z_G|-\epsilon)^2},\qquad S_h C_1\le1 .
\]
Take $\eta\in\partial\Pi_h$, with positive modulus. Suppose that the average of the arguments of $\mu_{G,i}/v_i$ differs from $\arg\eta$ by at most $d_0$, choosing representatives of all these angles in the same neighbourhood of $\arg z_G$. For a real $t$ write $[t]_+=\max(t,0)$. Define
\[
 E_0=q\,\big[S_h d_0-V_G+1-|\eta|^2\big]_+ .
\]
If $z$ violates the side with normal $e^{i\psi}$, then for every $0<\delta<1$,
\[
 \Gamma+\delta Y
 <\frac{E_0}{c_h}
 +p\left[\,2(1-c_h)+
 \frac{\tan^2(\pi/h)(Z-1)}{(1-\delta)c_h^2 |z_G|^2}\,\right].
 \tag{P0}\label{app:facet}
\]
\end{lemma}
\begin{proof}
Put $w=\mu_G-z_Gv$, so $\langle\bar v w\rangle=0$. The difference of the weighted and unweighted angular means only involves the nonconstant part of the reciprocal radii:
\[
 \big|\langle w/v\rangle\big|
 =\big|\langle (1/v-\bar v)w\rangle\big|\le\sqrt{Z-1}\sqrt H,
 \qquad \langle |w/v|^2\rangle\le H/\alpha^2.
\]
Taylor's formula for $\log(z_G+w_i/v_i)$, integrated on the segment from $z_G$, yields
\[
 \operatorname{dist}(\arg z_G,\arg\eta)
 \le d_0+|z_G|^{-1}\sqrt{Z-1}\sqrt H + C_1 H .
\]
The squared radial function $f_h$ of $\Pi_h$ has Lipschitz constant $S_h$, whereas $|z_G|^2=1-V_G-H$. Thus
\[
 2c_h q[\operatorname{Re}(e^{-i\psi}z_G)-c_h]_+
 \le q[|z_G|^2-f_h(\arg z_G)]_+
 \le E_0+ S_h q |z_G|^{-1}\sqrt{Z-1}\sqrt H .
\]
The first inequality follows from the radial comparison whenever the left side is positive, and otherwise from nonnegativity of the right side. In the second we used $S_hC_1\le1$. Now $z=qz_G+p z_B$. Combine the strictly violated side inequality with the exact bad supporting-line identity and with $q\sqrt H\le\sqrt{pY}$. Young's inequality absorbs the resulting square root into $(1-\delta)c_hY$ and the last term of~\eqref{app:facet}.
\end{proof}

If every leaking permutation loses squared norm at least $K/h^2$, then $\Gamma\ge Kp_{\rm leak}/h^2$. On a nonleaking bad draw with some angular error at least $\tau$ relative to $\arg z$, its contribution to $Y$ is at least $4\alpha^2\sin^2((\tau-\pi/h)/2)/h$ when $\psi$ is the normal of the side at $\arg z$ and $\tau>\pi/h$. The coefficient of $p/h^2$ from the square brackets of~\eqref{app:facet} tends to $\pi^2[1+(Z-1)/(1-\delta)]$ if $|z_G|\to1$. We next bound this error for the prefixes that occur in the proof.
\subsection{How large can a single radial gap be?}\label{app:individual}
Sorting the sites makes a gap visible. If a permutation takes a site
across the gap, the first $h$ radii lose a definite amount of squared
norm. If it preserves the first $h$ sites, we can follow the cycles
there and use their winding. We count both contributions against the
same supporting side of $\Pi_h$. This gives a bound which reflects
the spread of the actual radii through the reciprocal moment $Z$.

Write $t=h/N$. Use $K_0=10.15$ if $d=2$, $K_0=10.10$ if $d>N/4$, and $K_0=10$ otherwise. Rescale each prefix temporarily to mean square one, so that $Z=\langle |v|^{-2}\rangle$ is the product of the physical mean square and the physical mean reciprocal square. Suppose a prefix ends at physical radius $\rho$ or above, and includes at least $b_0N$ of the stated labels at physical radius $H_0$ or above. We shall prove that it is proper and that its normalized squared gap is strictly less than
\[
                 K/h,\qquad K=K_0+11.85(Z-1).                    \tag{P1}
\]
The inputs for the three smallest denominators are
\[
\begin{array}{l|rrrrrrrr}
\text{class}&2&3&4&2a&34a&2b&34b&3c\\\hline
\rho&.53&.53&.54&.51&.50&.50&.49&.48\\
b_0&.7975&.7459&.7217&.850&.850&.895&.895&.910\\
H_0&.9779&.9851&.9887&.970&.970&.965&.965&.965\\
X_0&7.443&5.539&4.382&7.443&5.539&7.443&5.539&5.539
\end{array}
\]

Here $L=5.6$ bounds $Nhx/|1-z|$. Classes $2,3,4$ use the hull, as shown in the Farey table for the final angular test. Classes $2a,34a,2b,34b,3c$ use the indicated additional high population, and $34$ covers both denominators with the larger deficit from order three. The other classes are
{\small\setlength{\tabcolsep}{3.2pt}\begin{center}\begin{tabular}{lrrrrr}
class&$\rho$&$b_0$&$H_0$&$X_0$&$L$\\\hline
5&.55&.706&.99022&3.622&5.6\\
6--7&.57&.689&.99185&3.088&5.6\\
8--18&.59&.652&.99388&2.39&5.6\\
5--8 R&.52&.79&.976&3.622&5.6\\
19--36&.61&.616&.9968&1.108&5.9\\
$37\le d\le N/4$&.7&.499&.998&.641&5.7\\
$d>N/4$&.75&.428&.9984&.47&13\\
\end{tabular}\end{center}}

The row $5$--$8$ R will be used when $.79N$ labels are at or above $.976$. For the ordinary classes from $5$ through $18$, use $b_0=b_j$, $H_0=H_i$ rounded down, $X_0=X_i$ rounded up on a group $i\le d\le j$. Later, a rank calculation on a single $d$ may use the corresponding more precise inputs from the preceding section. A containing prefix has the full specified high population. We use $b,H,X_*$ for $b_0,H_0,X_0$ in the estimates that follow.

In an interval $A\le h/N\le T$ let $h_0=\lceil nA\rceil$, $\lambda=b/T$, and
\[
 V_1=\frac1{\rho^2+\lambda(H^2-\rho^2)},\qquad
 Z_1=1+D_\rho\,w(\lambda)+D_H/4,\quad
 \]
\[
 D_s=(s-1/s)^2,\qquad
 w(\lambda)=\begin{cases}\lambda(1-\lambda)&\lambda\ge.5,\\1/4&\lambda<.5.\end{cases} \tag{P2}
\]
The normalized maximum squared radius $V$ is at most $V_1$. To see that $Z\le Z_1$, write the moment product as the sum of contributions of pairs of squared radii. It is separately convex in each, so each high radius can be placed at $H$ or $1$ and each remaining radius at $\rho$ or $1$. If the proportions at $\rho,H$ are $p,q$, respectively, the excess over one is
\[
 p(1-p)D_\rho+q(1-q)D_H-pq[D_\rho+D_H-(\rho/H-H/\rho)^2].
\]
The bracket is nonnegative and $p\le1-\lambda$, proving (P2).

\paragraph{A supporting side bounds the probability of a bad draw.}
For the moment return to the physical normalization $1=r_1\ge\cdots\ge r_h=m$. Write $L_0=\sum_{i\le h} r_i$, $Q_0=\sum_{i\le h}1/r_i$ and
\[
 H(s)=\max_{\sum t_i\equiv s\ (2\pi)}\sum_{i\le h}r_i\cos t_i.
\]
We call $s$ the twist. In a direction of angle $\psi$ the deviations $t_i$ add to $\sum_{i\le h}\arg v_i-h\psi$ modulo $2\pi$, because the coefficient arguments add to zero modulo $2\pi$. Thus $H(s)$ is the largest projection of the marked sums $\sum u_i v_i$ with $|u_i|=1$ and $\prod u_i=1$ in that direction. The maximum at $0\le s\le\pi$ may be restricted to deviations $t_i\ge0$ adding to $s$: shorten their absolute principal deviations until their sum is $s$. Equal deviations give $H(s)\ge L_0\cos(\pi/h)$. Since $m>L_0(1-\cos(\pi/h))$ here, every maximizing deviation is less than $\pi/2$. Strict concavity gives $t_i=\arcsin(\lambda_s/r_i)$, $\sum_i t_i=s$, and $H'(s)=-\lambda_s$. At positive $s$ an optimizer has no zero deviation, by transferring a small angle to it. Set
\[
 D=L_0-h^2/Q_0,\quad \delta_*=\arcsin[\pi/(mQ_0)],\quad
 C_*=1+\tfrac12\tan^2\delta_*,\qquad G(s)=H(s)-L_0\cos(s/h).
\]
The estimate
\[
 D\sin(s/h)/h\le G'(s)\le C_* D\sin(s/h)/h                         \tag{P0a}
\]
retains the radial variance. To verify it, put $x=s/h$, $\lambda=\lambda_s$. Concavity of sine gives $\lambda Q_0/h\le\sin x$; moreover $\lambda\le s/Q_0$ and $t_i\le\delta_*$. Its second derivative and the Lipschitz constant of arcsine on this interval give
\[
 0\le\sin x-\lambda Q_0/h
 \le\tfrac12(\lambda/m)\lambda^2\operatorname{Var}(1/r_i)/
                       \cos^2\delta_*.
\]
Here variance uses the uniform measure on the prefix. Since
$\sum(1/r_i-Q_0/h)^2\le (Q_0/h)^2D/m$, as follows from
$\sum(r_i-h/Q_0)^2/r_i=D$, the excess of $G'$ over its lower bound is at most $D\sin x\,\tan^2\delta_*/(2h)$. This proves (P0a).

At $B_0=\pi-\operatorname{dist}(h\arg z,2\pi\mathbb Z)$ the radius of $\Pi_h$ on this ray is $R_h=\cos(\pi/h)/\cos(B_0/h)$. Test the output in a direction with twist $\pi$, choosing the input coefficients to maximize the projection in its preimage direction, which has twist $B_0$ up to sign. The signed slack of a draw is the supporting value $H(\pi)$ minus its output projection. Its expectation is exactly $H(\pi)-rH(B_0)$. Since $r>R_h$, (P0a) integrated over $[B_0,\pi]$ and $[0,B_0]$ gives
\[
 H(\pi)-rH(B_0)\ \le C_*D(1-R_h)\le C_*D\,\pi^2/(2h^2).
\]
A preserving draw has nonnegative slack. If some angular change differs from $\arg z$ by at least $\tau$, one deviating argument in the output frame has magnitude at least $\tau-\delta_*$. Its slack is at least $d_A=m[1-\cos(\tau-\delta_*)]-\pi^2\bar r/(2h)$, where $\bar r=L_0/h$. For a leaking draw the sum of available radii falls by at least $r_h-r_{h+1}$. With $S=h^{-1}\sum r_i^2$, the assumed normalized squared gap $K/h$ implies slack at least $d_L=(KS/m-\pi^2\bar r)/(2h)$. Thus the probability of the union of the two bad events is bounded by $C_*\pi^2(D/h)/(2h\min(d_A,d_L))$, omitting $d_L$ at a full prefix.

Here $\tau\ge.78$, $h\ge193$, $m\ge.48$. The bound on $D/h$ below gives a $d_A$ ratio below $.05$, since $\delta_*<.035$, $C_*<1.001$. Indeed $D/h=\min_c h^{-1}\sum(r_i-c)^2/r_i\le(1-\sqrt m)^2$, by taking $c=\sqrt m$ and testing endpoints. If a fraction $b$ is at least $H>m$, then $S/m-\bar r\ge bH(H/m-1)$. Thus the $d_L$ ratio is at most $C_*\pi^2/K$ times
\[
 \frac{m(1-\sqrt m)}{1+\sqrt m}\frac{1-m}{bH(H-m)}.                \tag{P0b}
\]
The first factor decreases for $m\ge.48$ and the second increases up to $H$. For the three smallest denominators use $(m,b,H)\ge(.48,.7217,.965)$ componentwise. On the intervals starting at $.48,.75,.88$ and ending at $.75,.88,.94$, the products are below $.146,.110,.097$. For the remaining five types use respectively $(.55,.652,.990),(.52,.79,.976),(.61,.616,.9968),(.70,.499,.998),(.75,.428,.9984)$. The same intervals starting at the new minimum give upper products below $.133,.129,.068$ in all cases. For $m\ge.94$ use $\Delta/[(1-\pi^2/K)\bar r]<.077$ with $\Delta=D/h$. The factor $C_*\pi^2/K$ is less than one. Thus the bad probability is less than $.15$.

\paragraph{Winding differences of good permutations.}
Represent the doubly stochastic matrix by a random permutation of the labels. Each good permutation preserves the prefix, and each image's angle differs from the rotation angle $\theta=\arg z$ by at most $\tau$. The angle $\tau$ is specified in (P4). On a pair of independent good permutations $P,Q$ define
\[
 g=2\pi/h,\ a=4\sin^2(\pi/h),\ E=\|v_Q-v_P\|^2,\ 
 B=\langle(v_Q-v_P)^2\rangle,\ T=\|v_Q-e^{igk}v_P\|^2,\quad e=E/a .
\]
Here $k$ is the difference of the two integer windings, determined by lifting all good single-label angle differences near $\theta$. The relative angles $u_i$ satisfy $|u_i|\le2\tau<\pi$ and $\langle u_i\rangle=gk$. Put $w_i=|v_{Q(i)}v_{P(i)}|$. The chord energy bounds the sum of its radial energy and $\chi\langle w_i u_i^2\rangle$ whenever $\chi\le\operatorname{sinc}^2(u_i/2)$ on all edges. Here $\operatorname{sinc} t=\sin(t)/t$, with value one at zero. Also $\langle w_i\rangle\le1$ and $\langle1/w_i\rangle\le Z$ by Cauchy--Schwarz. In this weighted Hilbert space the projection $u\mapsto u-\langle u\rangle$ has squared operator norm $\langle w\rangle\langle 1/w\rangle$. Indeed this is the complementary projection to the rank-one map $u\mapsto\langle u\rangle\mathbf1$; orthogonally splitting off the span of $1,1/w$ proves the identity in two dimensions. Thus $T\le (Z/\chi)E$, by $4\sin^2((u_i-gk)/2)\le(u_i-gk)^2$. Also $e\ge\chi k^2/Z$. In the rows below $\chi=\operatorname{sinc}^2\tau>Z/2$. If $|k|=2$ and $e<2.5$, one edge gives $u_*:=\max|u_i|\le2\arcsin\sqrt{2.5\pi^2/(h\rho^2)}$. Using $h_0$ and $Z_1$ for each row gives $4\operatorname{sinc}^2(u_*/2)/Z>2.70$. This contradiction and $e>k^2/2$ at $|k|\ge3$ give
\[
  k^2-|k|\le4(e-|k|),\qquad T\le10a(e-|k|)\quad (|k|\ge2).
\]
For $k=0$ one has $T=E$. The one-dimensional projection identity in the constant-weight case is used simply as the inequality with constant one.

For $|k|=1$ conjugate if necessary to take $k=1$. Projection onto $1/w$ gives
\[
 \langle w u\rangle\ge g/Z-\sqrt{1-1/Z}\sqrt{E/\chi-g^2/Z}.
\]
Indeed the weight reciprocals have mean at most $Z$ and the right side decreases with this moment. The exact identity $T=E\cos g+a-2\sin g\langle w\sin u\rangle$, and
$|\sin u-u|\le u_* u^2/6$, therefore give
\[
 \frac Ta\le \Phi(e;u_*,\chi,Z):=
 e+1-1.9998/Z+2\sqrt{(1-1/Z)(e/\chi-1/Z)}
    +(\sin g)u_* e/(3\chi).                                  \tag{P3}
\]
Here $g\cot(g/2)>1.9998$. Below $e\le1.12$ use
$u_c=2\arcsin\sqrt{1.12\pi^2/(h_0\rho^2)}$, $\chi_c=\operatorname{sinc}^2(u_c/2)$, $h_0=\lceil An\rceil$. Otherwise use $u_*=2\tau$. In each range $\Phi(e)/(e-1)$ decreases on $e>1$ by direct differentiation.

The next two tables, collectively (P4), specify the prefix fraction, ordinary angular tolerance $\tau$, facet parameter $\delta$ and winding thresholds $y_0,t_0$. Here $[A,T]$ measures the rank $h/N$, whereas $\tau$ in these two
tables is an angle in radians. The separation matters, for a rank such
as $.80$ fixes neither the angle of an individual move nor its error.
The last three columns give bounds for the quantities following (P10).
Read, for example, the first row at $N=450$, $h=360$ in class $2$.
Its first margin gives $1-S>.394$, so the angular error with
coefficient $S$ can be absorbed. Of the bad permutations, one
that stays in the prefix but makes a large angular error has margin
$R_n>.291$. One that sends a label out has margin $R_l>.1164$.
Either pays more than the violated side permits. The same row covers
every rank in its indicated interval, under the same prefix hypotheses.
The second table gives an upper bound for $S$, which we compare with one.
{\small\setlength{\tabcolsep}{2.8pt}
\begin{center}\begin{tabular}{llrrrrrrrrr}
class&$[A,T]$&$\tau$&$\delta$&$y_0$&$t_0$&$C_T$&$C_B$&$1-S>$&$R_n>$&$R_l>$\\\hline
2&$[.7975,.85]$&.98&.13&.056&.167&12.5&37&.394&.291&.1164\\
2&$[.85,.94]$&.96&.15&.073&.227&13&29&.405&.857&.0496\\
2&$[.94,1]$&.88&.16&.061&.24&17.6&34&.309&.14&.0034\\
3&$[.7459,.85]$&1.12&.12&.079&.199&12.5&27&.077&.482&.0786\\
3&$[.85,.94]$&1&.14&.068&.231&16.2&31&.303&.285&.0788\\
3&$[.94,1]$&.96&.14&.048&.234&24.3&43&.105&.087&.0734\\
4&$[.7217,.85]$&1.12&.12&.077&.194&12.5&27&.101&.274&.078\\
4&$[.85,.94]$&1.02&.13&.066&.22&16.8&32&.083&.264&.0708\\
4&$[.94,1]$&.98&.13&.041&.233&27.9&50&.042&.195&.0742\\
2a&$[.85,.94]$&.96&.15&.074&.242&12.5&29&.381&.488&.0539\\
2a&$[.94,1]$&.9&.16&.06&.259&16.7&35&.284&.014&.0005\\
34a&$[.85,.94]$&1.1&.12&.077&.255&12.5&27&.096&.305&.0751\\
34a&$[.94,1]$&.96&.15&.052&.265&19.9&40&.157&.106&.0779\\
2b&$[.895,.94]$&.96&.14&.053&.218&12.5&39&.391&.738&.073\\
2b&$[.94,1]$&.92&.16&.065&.265&13.3&32&.386&.706&.0157\\
34b&$[.895,.94]$&1.1&.11&.056&.23&12.5&37&.103&.37&.0812\\
34b&$[.94,1]$&.96&.15&.052&.265&17&40&.237&.284&.0824\\
3c&$[.91,.94]$&1.1&.11&.045&.217&12.5&46&.079&.437&.0825\\
3c&$[.94,1]$&.98&.15&.05&.265&16.8&41&.21&.426&.083\\
\end{tabular}\end{center}}
{\small\setlength{\tabcolsep}{3.2pt}\begin{center}\begin{tabular}{lrrrrrrrrrr}
class&$[A,T]$&$\tau$&$\delta$&$y_0$&$t_0$&$C_T$&$C_B$&$S<$&$R_n>$&$R_l>$\\\hline
5&[.706,.94]&1&.16&.064&.219&17.1&33&.793&.473&.079\\
5&[.94,1]&.86&.16&.036&.234&30.7&57&.979&.19&.075\\
6--7&[.689,.92]&.98&.15&.065&.198&15.4&32&.682&.226&.083\\
6--7&[.92,1]&.86&.15&.037&.23&27.8&56&.841&.417&.079\\
8--18&[.652,.92]&.98&.15&.06&.224&16.1&35&.691&.518&.084\\
8--18&[.92,1]&.82&.15&.028&.227&34&73&.928&.4&.081\\
5--8 R&[.79,.92]&1&.15&.073&.242&13.1&29&.627&.396&.084\\
5--8 R&[.92,1]&.96&.15&.044&.248&25.2&47&.947&.426&.075\\
19--36&[.616,.92]&.96&.15&.048&.226&18.6&43&.738&.348&.083\\
19--36&[.92,1]&.78&.15&.024&.223&36.5&85&.887&.326&.08\\
$37\le d\le N/4$&[.499,1]&.88&.15&.035&.197&16.4&59&.63&.232&.092\\
$d>N/4$&[.428,1]&.98&.12&.03&.158&13.7&68&.581&.363&.086\\
\end{tabular}\end{center}}
Use $I=1$ when $|k|=1$ and $e<1+y_0$, and $I=0$ otherwise. Put $c=\langle v\rangle$, $M_2=\langle v^2\rangle$, and $D'=E-a|k|+1000Ia|c|^2$. The radial-excursion coefficient used below is $D_c$, rounded down to the fourth decimal place from
\[
 \min\left\{\frac{\operatorname{sinc}^2(\pi/[2(1+(1-\lambda)h_0)])}{4(1-\lambda+1/h_0)^2},\
 \chi_c\left[\frac{1+1/\rho}{1-t_0+\rho^{-2}(1-1/V_1)}\right]^2\right\}.
\]

\paragraph{The long winding cycle.}
For $I=1$ the sum of the reciprocal squared radii on a positive winding cycle, divided by $h$, is at least $\chi_c/e$ by Cauchy--Schwarz, in prefix normalization. Thus there is just one, its winding is one, and its omitted fraction $t'$ satisfies $t'\le V(Z-\chi_c/e)\le V_1(Z_1-\chi_c/(1+y_0))<\lambda$ in each row of (P4). In particular it meets a high label. Write $Q=1-t'$, $v=V-1$, and $u=h^{-1}\sum_{\rm cycle}|v_i|^{-2}$. Angular energy on the cycle is at least $a\chi_c/u$. The sublevel excursions of $f_i=(H\sqrt V-|v_i|)_+$ have at most $(1-\lambda)h$ positive consecutive labels. By odd reflection, an excursion of length $j$ costs at least $4\sin^2(\pi/(2(j+1)))\sum f_i^2$ in squared radial differences. Division by $a\le4(\pi/h)^2$ gives the first entry defining $D_c$. Both $h/(1+(1-\lambda)h)$ and its accompanying sinc factor increase with $h$.

Put $U_0=(1+Q)/V-1/V^2$. Expanding $|v_i|^{-2}$ at $V$ and using the off-cycle squared mass at most $Vt'$ yields
\[
 u\le U_0+\frac{(1+1/\rho)^2}{V^2h}
                   \sum_{\rm cycle}(\sqrt V-|v_i|)^2 .
\]
Indeed $r_i^{-2}=2/V-r_i^2/V^2+(V-r_i^2)^2/(V^2r_i^2)$. Hence the disjoint radial and angular energy terms give
\[
 e\ge\chi_c/u+
 D_c\left[\frac{V\sqrt{(u-U_0)_+}}{1+1/\rho}
                         -(1-H)\sqrt{VQ}\right]_+^2 .
\]
Also $u\le U=[Q+(Q+\rho^{-2})v]/(1+v)^2$ by the secant estimate $Z\le[(1+\rho^2)V-1]/(\rho^2V^2)$ and the off-cycle reciprocal mass $t'/V$. Suppose $t'\ge t_0$, and put $Q_0=1-t_0$. For $Q\le Q_0$ the right side as a function of $u\le U$ is decreasing since its derivative is at most $-\chi_c/U^2+D_cV^2/(1+1/\rho)^2\le0$. For this bound substitute $VU\le Q_0+\rho^{-2}(1-1/V_1)$ and the respective $D_c$ in (P4). Thus substitute $u=U$. With $q_0=(1-\rho)/(1+\rho)$ the first term inside the square brackets is exactly $\sqrt{q_0v}$. Both the angular and radial lower bounds now decrease as $Q$ increases, so use $Q_0$. Set
\[
 \omega=.38,\quad A=D_c[q_0-(1-H)\sqrt{q_0Q_0}(\omega+1/\omega)],
 \qquad A_0=D_c(1-H)\sqrt{q_0Q_0}\,\omega .
\]
The squared positive part is at least the first squared term minus twice the product. Use $2\sqrt{v(1+v)}\le\omega(1+v)+v/\omega$ to obtain $e\ge\chi_c/U+Av-A_0$. This exceeds $1+y_0$ if
\[
 [\chi_c+A(Q_0+\rho^{-2})]v^2+
 [2\chi_c+AQ_0-(1+y_0+A_0)(Q_0+\rho^{-2})]v+
 \chi_c-(1+y_0+A_0)Q_0>0 .                                 \tag{P5}
\]
Its minimum on $[0,V_1-1]$ is greater than $.00019$ on each row of (P4), by clamping the vertex to this interval and using (P2). Thus $t'<t_0$.

On this cycle of $l=Qh>162$ labels write the coordinates as
$\sum_j d_j e^{2\pi ijm/l}$ in cycle order $m$. Let $a_j=4\sin^2(\pi j/l)$, and use a prime to omit $0,1,-1$. Write $V_o,E_o$ for the squared mass and energy off the cycle, divided by $h$. Parseval gives
\[
 D'= E_o+Q\sum{}'(a_j-a_1)|d_j|^2+
           a_1(1-V_o-Q|d_0|^2)-a+1000a|c|^2.
\]
Here $V_o\le Vt'$, $Q|d_0|^2\le[1.007V(t')^2+(1+1/.007)|c|^2]/Q$, $a_1/a\le Q^{-2}$, and $a/a_1\le Q^2+.0003t'$ by integrating $1.99974<2s\cot s\le2$ on $0<s\le\pi/162$. The last error follows from $-Q\log Q\le1-Q$ and $Q>.73$. Since $a_j\ge3.9a_1$ on primed modes,
\[
 D'\ge E_o+.7Q\sum{}'a_j|d_j|^2+600a|c|^2
               +a_1 t' b(t'),\quad
 b(t')=2-t'-V_1-1.007V_1t'/(1-t')-.0003 .                 \tag{P6}
\]
In particular $D'\ge0$. Write $y=D'/a$. If $y\ge y_0$, (P3) at $1+y_0$ suffices by monotonicity of its numerator, and $|B-IaM_2|\le(2+y_0)a$. If $y<y_0$, then $t'<.10$. Indeed $t'b(t')-y_0((1-t')^2+.0003t')$ is concave for $.10\le t'\le t_0$, and its two endpoint values exceed $.00016$ on all the rows of (P4).

We give numerical bounds for the remaining modes in a form that uses the same formulas in every row. Set
\[
 c_t=1/b(.10),\quad v_o=V_1c_t,\quad c_0=(\sqrt{.10v_o}+1/\sqrt{600})^2/.90,
 \quad C_e=v_o+c_0+.367.
\]
Thus $t'\le c_ty$, $V_o\le v_o y$, $Q|d_0|^2\le c_0y$, $Q\sum'|d_j|^2\le.367y$. The last inequality follows from (P6) and $a_j\ge3.9a_1$. The complex moment therefore satisfies $|M_2-2Qd_1d_{-1}|\le C_e y$.

For the signed area $S=\langle\operatorname{Im}(\overline{v_P}v_Q)\rangle$ projection as in (P3) gives
\[
 S/g\ge s_0=1/Z_1-\sqrt{(1-1/Z_1)((1+y_0)/\chi_c-1/Z_1)}
 -u_c(2\pi/h_0)(1+y_0)/(6\chi_c).
\]
Write $[u,v]=\operatorname{Im}(\bar u v)$, the signed double area of the triangle $0,u,v$. The off-cycle contribution in these units is at most $t'E_o/a$ by the cycle sector inequality $|\sum [v_i,v_{i+1}]|\le\frac12\cot(\pi/j)\sum|v_{i+1}-v_i|^2$ on cycles of length $j\le t'h$. This inequality follows by expanding in discrete Fourier modes. The higher modes on the long cycle cost at most $Qy/(2\cdot.7)$, since $|\sin(2\pi j/l)|/a_j\le l/(4\pi)$ on the primed modes. Hence
\[
 Q(|d_1|^2-|d_{-1}|^2)\ge Q[s_0-(.10+1/1.4)y].
\]
Their combined squared mass is at least $1-C_e y$. Consequently $|d_1|^2\ge m_0=[1-C_e y_0+.90(s_0-(.10+1/1.4)y_0)]/2$. In all the rows of (P4) $m_0>.590$. The squared chord multiplier on the negative mode for $T$ is at most $4.47a$. Put $c_-=4.47/(4\cdot.90m_0)$ and $p_-=3.6/c_--1>0$. Its contribution is at most
\[
 3.6a|M_2|^2 +c_-(1+1/p_-)C_e^2y_0D' .
\]
This follows from $(A+B)^2\le(1+p_-)A^2+(1+1/p_-)B^2$. The positive mode has squared multiplier at most $1.25a(t')^2$, the constant mode has multiplier squared $a$, the primed multipliers are at most $3.1(a_j-a_1)$, and the off-cycle part is at most $2E_o+2aV_o$. For the primed bound use $(\sqrt{a_j}+\sqrt a)^2/(a_j-a_1)\le (1+1/\sqrt{3.9})^2/(1-1/3.9)<3.1$. The primed mass times $a_j-a_1$ is at most $D'$ by the Parseval identity and the positive remainder of (P6). The combined $D'$ coefficient is bounded by
\[
 C_{\rm near}=c_-(1+1/p_-)C_e^2y_0+1.25(.10)c_t+c_0+3.1+2+2v_o<14.41.
\]
In $B-aM_2$, the two first modes cost at most $a_1-a\le2.4at'$, and the others at most $Q\sum' (a_j+a)|d_j|^2+aQ|d_0|^2+E_o+aV_o$. The extra term $aQ\sum'|d_j|^2$ may be retained. This gives a coefficient at most $2.4c_t+1/.7+.367+c_0+1+v_o<12.55$. We have proved, with $G=2Z_1/\chi_c$,
\[
\begin{gathered}
D'\ge0,\qquad a(k^2-|k|)\le4D',\qquad |B-IaM_2|\le C_BD',\\
 T\le C_TD'+3.6 I a|M_2|^2,\qquad C_B=\lceil (2+y_0)/y_0\rceil .          \tag{P7}
\end{gathered}
\]
Here $C_T$ rounds upward to the next tenth the maximum of $12.5$, $C_{\rm near}$, $\Phi(1+y_0;u_c,\chi_c,Z)/y_0$ and $\Phi(1.12;2\tau,\operatorname{sinc}^2\tau,Z)/.12$. In these upper bounds one may replace $\sin g$ by $\sin(2\pi/h_0)$, and $Z$ by any local upper bound. At nonexceptional $|k|=1$ the complex bound uses $|B|\le E$.

\paragraph{From these estimates to one violated side of the regular polygon.}
Recall the normalized Hilbert norm on the prefix. Let $s$ be the mean fraction of prefix labels leaving it and $p$ the probability of a bad draw, $q=1-p$. The variance about $zv$ is $x-\Gamma$, where $\Gamma=1-\mathbb E\|v_P\|^2$. Matching incoming and missing labels gives $\Gamma\ge Ks/h$, $|(1-z)c|\le2\sqrt V s$ and a change of at most $2Vs$ in the expected second complex moment. Every bad nonleaking draw has error at least $4r\rho^2\sin^2(\tau/2)/h$ about $zv$. A leaking draw has at least one outgoing label. Therefore
\[
 p\le h^2\Gamma/K+[4r\rho^2\sin^2(\tau/2)]^{-1}h(x-\Gamma)
      \le h^2x/K .
\]
The last coefficient comparison holds since $h\ge193$, $K<15$, $r>.9999$, $\rho\ge.48$ and $\tau\ge.78$. The supporting-side estimate (P0a)--(P0b) also gives $p<.15$.

Write $qV_G=q\mathbb E_G\|v_P-\mu_G\|^2$, let $p_*=\mathbb E_G I$ on independent good pairs, and let $J=q\mathbb E_G D'/2$. Put $D=p\mathbb E_B\|v_P-\mu_B\|^2+pq\|\mu_B-\mu_G\|^2$, so $x=\Gamma+D+qV_G$. Let the mean of the good integer windings $k_P$ be $j+u$, $j\in\mathbb Z$, $0\le u<1$. Put $\eta=e^{igj}(1-u+ue^{ig})$, $b_G=\mathbb E_G e^{igk_P}$ and $W=a(\mathbb E_G|k|/2-u(1-u))$. Here $W\ge0$: the half mean absolute difference of an integer-valued random variable is $\sum F(j)(1-F(j))$, which is concave in its finite-support law. With the mean fixed it is minimized on a two-point law, where, in terms of the two distances from the mean, it is $AB/(A+B)$, minimized at adjacent integers. The integer-mean case is immediate. The second finite differences of $e^{igj}$ have magnitude $a$; summation from $j,j+1$ gives
\[
 |b_G-\eta|\le \tfrac a2[\operatorname{Var}_G k_P-u(1-u)]\le W/2+2J/q .
\]
The real variance identity and its complex analogue now give
\begin{align*}
 q(V_G-1+|\eta|^2)&=J+qW-500qap_*|c|^2,\\
 |(1-z^2-qap_*/2)M_2|&\le C_BJ+D+2Vs .
\end{align*}
Indeed the within-good complex variance equals the mean of $B/2$ on good pairs; the other complex variances have modulus at most $D$. Also $qap_*\le G(x-\Gamma)$ by $E\ge\chi_c a/Z$ on exceptional pairs. Unconditional complex variance gives $|1-z^2||M_2|\le x-\Gamma+2Vs$. Multiplying by $x-\Gamma$ and using $(x-\Gamma)\Gamma\le x^2/4$ therefore gives
\[
 (x-\Gamma)|M_2|\le x^2[1+Vh/(2K)]/|1-z^2|.
\]
With (P7) we conclude that
\[
 q\mathbb E_G T\le 2C_TJ+F(C_BJ+D+2Vs),                       \tag{P8}
\]
where one may take
\[
 F=\frac{3.6G(1+Vh/(2K))}{(.315N)(.315N-G/2)}\quad(d=2),\qquad
 F=\frac{3.6G(80)}{3.8N^2-40G}\quad(d\ge3).
\]
For the second use $|M_2|\le1$. The mean correction is bounded by
$500qap_*|c|^2\le500GV[L/(NK)]^2\Gamma$ using the same product inequality.

It remains to compare the average argument of the row means with $\arg\eta$. For $|\arg(W/U)|\le\gamma<\pi$, $|W|,|U|\ge m$ one has
\[
 \left|\log(W/U)-(W/U-1)\right|
 \le |W-U|^2\frac{\operatorname{arsinh}(\tan(\gamma/2))}
                            {2m^2\sin(\gamma/2)} .                    \tag{P9}
\]
Indeed $(z-1-\log z)/(z-1)^2=\int_0^\infty[(1+s)^2(z+s)]^{-1}ds$. Substitute $z=Re^{i\alpha}$ and $s=Rt$. After multiplication by $\min(1,R)^2$, the absolute integrand in $t$ is at most $(1+t)^{-2}|e^{i\gamma}+t|^{-1}$, since $\min(1,R)/(1+Rt)\le1/(1+t)$. Substitution $v=1/(1+t)$ and symmetry about $v=1/2$ evaluate this integral. The branch is the one containing the segment from $1$ to $z$.

\begin{lemma}\label{proj500}
Let $F(w)=\log w-(w-1)$ on the principal branch. If
$|\arg w|\le2.10$ and $|\alpha|\le .18$, then
\[
 \min(1,|w|)^2 |\operatorname{Im}(e^{i\alpha}F(w))|
                              \le .516\,|w-1|^2 .
\]
For $\alpha=0$ the factor $1/2$ suffices up to $|\arg w|\le2.4$.
\end{lemma}

\begin{proof}
Let $t=|\arg w|$, $s=\sin t$, $B=1-\cos t$, $A=t-s$ and
$d=\max(|w|,|w|^{-1})-1$. Dividing by the squared distance and normalizing by the smaller squared modulus leaves the denominator
$D=d^2+2(1+d)B$. For $|w|\ge1$ the real and signed imaginary numerators are $(1+d)B-[d-\log(1+d)]$ and $A-ds$; for $|w|\le1$ they are $B/(1+d)-[\log(1+d)-d/(1+d)]$ and $A+ds/(1+d)$. The two square brackets are between $0$ and $d^2/2$ by integration, so the normalized real part has magnitude at most $1/2$.

For either radius the absolute imaginary numerator is bounded by $A+ds$. Thus it is at most $pD$ whenever $2pB\ge A$ and
\[
            4p(2pB-A)\ge (s-2pB)_+^2 .
\]
If $s>2pB$ the last test is equivalent to
$4p(sB-A)\ge(1-4p^2)s^2$, by $s^2+B^2=2B$. For $p=1/2$ this holds up to $\pi/2$ since $sB-A=s(2-\cos t)-t$ has derivative $2\cos t(1-\cos t)\ge0$. Above $\pi/2$ we have $s\le B$. On $[0,2.4]$ one has $A\le B$: the derivative of $B-A$ changes from positive to negative at $\pi/2$, and its endpoint value at $2.4$ is positive.

For $p=.42$ consider $.60\le t\le2.10$. The quotient $A/B$ increases, since it is the weighted average of the increasing $A'/B'=\tan(t/2)$ from $0$ to $t$. At $2.10$, $A/(2B)<.412$. On $[.60,1.75]$ the function $(sB-A)/s^2$ has derivative $2\cos(t)A/s^3$, and its endpoint values exceed $.198$ and $.406$. These exceed $(1-4p^2)/(4p)<.176$. For $t\ge1.75$ one has $s/B=\cot(t/2)<.84$. We obtain a bound of $.42$ for the normalized imaginary part throughout the stated interval, and hence a projected bound $.42+.5\sin(.18)<.510$.

Finally the absolute normalized logarithmic remainder at angle $t$ is at most
\[
 K(t)=\operatorname{arsinh}(\tan(t/2))/(2\sin(t/2)).
\]
Indeed $(w-1-\log w)/(w-1)^2=\int_0^\infty[(1+v)^2(w+v)]^{-1}dv$.
After $v=|w|y$ and multiplication by the smaller squared modulus,
bound the absolute integrand by $(1+y)^{-2}|e^{it}+y|^{-1}$. Integration with $x=1/(1+y)$ gives $K(t)$, increasing with $t$ by the same integral. At $t\le.60$ one has $K(t)<.516$.
\end{proof}

Set $U_i=\mu_{G,i}/b_G$, $W_{Pi}=e^{-igk_P}v_{P(i)}$, $\xi_P=e^{igk_P}/b_G$. By Cauchy--Schwarz on weighted angular errors,
\[
 \mathbb E_G(gk_P-\theta)^2\le Zx/[qr\operatorname{sinc}^2(\tau/2)]=:v_\theta .
\]
Thus $|b_G|\ge b_*=1-v_\theta/2$, $|\arg b_G-\theta|\le\arctan(\sqrt{v_\theta}/b_*)$. Conditional variance gives $|\mu_{G,i}-zv_i|\le e_0=\sqrt{phx/q}$ and $|z_G-z|\le e_s=\sqrt{px/q}$. The moduli of $U_i,W_{Pi}$ are at least $m=\max(\rho-e_0-.0001,\rho\cos\tau)$; the second bound follows by projecting the good images on the ray of $zv_i$. The opening of the pair is at most
\[
 \gamma=2\tau+\arcsin(e_0/(r\rho))+\arctan(\sqrt{v_\theta}/b_*).
\]
The term $2\tau$ includes the compensating phase $gk_P-\theta$, whose absolute value is at most $\tau$. Write $K_s=\operatorname{arsinh}(\tan(\gamma/2))/(2\sin(\gamma/2))$. We use constants $C,A_w$ satisfying
\[
 \left|\operatorname{Im}\left\langle\log\frac{\mu_G}{v b_G}\right\rangle\right|
   \le C\mathbb E_GT+\tfrac12A_w\mathbb E_G D'.
\]
The ordinary estimate always permits $C=K_s/(2m^2b_*^3)$, $A_w=0$. For $\gamma\le2.10$ put $v=v_\theta$, $e_b=\arctan(\sqrt v/b_*)$, $\omega=.18$, $u_s=.063$, $\eta_s=.09$, and
\[
 q_s=1-v/u_s^2,\quad p_s=\min(1,v/(\omega-e_b)^2),\qquad
 k_s=1+\lfloor(\omega-e_b-u_s)/g\rfloor.
\]
If $k_s\ge3$, $q_s>0$, set $\chi=\operatorname{sinc}^2\tau$,
\[
 L_s=(Z_1/\chi)/(1-Z_1/(\chi k_s)),\quad
 H_s=\left(\sqrt{(1-q_s)/q_s}+\sqrt{b_*^{-2}-1}\right)^2.
\]
and put $k_b=\min(.516,K_s)$, $d_b=K_s-k_b$. Use
\[
 C=\frac{k_b/b_*^2+d_b(1+1/\eta_s)p_sH_s}{2b_*m^2},
 \qquad A_w=\frac{d_b(1+\eta_s)L_s}{b_*m^2q_s}.
\]
Otherwise retain the ordinary estimate. These expressions may be evaluated with lower $h_0,b_*$ and upper $v,e_b,\gamma$.

Indeed the linear part in (P9) vanishes after averaging with factors $\xi_P$, and the lifted logs $\log(W_{Pi}/v_i)$ sum to zero on the prefix for each $P$. All angles are in one branch, since $|\arg(W_{Pi}/v_i)|\le2\tau$, $|\arg(U_i/v_i)|\le\gamma-2\tau$ and $\gamma<\pi$. Write $D_P=\|W_P-U\|^2$ and $V_W=\mathbb E_GT/2$. By Cauchy--Schwarz and $\mathbb E_G|\xi_P-1|^2=|b_G|^{-2}-1$, one has $\mathbb E_G D_P\le b_*^{-2}V_W$. Lemma~\ref{proj500} bounds the projected remainder by $[k_b\mathbb E_G D_P+d_b\mathbb E_G(\mathbf1_{\{|\arg\xi_P|>\omega\}}D_P)]/(b_*m^2)$.

The second event has probability at most $p_s$. Let $V_s$ be the mean of $W_Q$ conditional on $|gk_Q-\theta|\le u_s$, an event of probability at least $q_s$. Cauchy--Schwarz about the ordinary mean shows $\|V_s-U\|^2\le H_s V_W$. Between these disjoint classes $|k_P-k_Q|\ge k_s$, so $T\le L_sD'$ by $T\le(Z_1/\chi)E$, $E/a\ge\chi k^2/Z_1$, $D'\ge E-a|k|$. The symmetric pair with $P,Q$ reversed has the same $D'$ and is also available. Using the square of a sum on $W_P-V_s+V_s-U$ and Jensen over $Q$ gives
\[
 \mathbb E_G(\mathbf1_{\{|\arg\xi_P|>\omega\}}D_P)
 \le (1+\eta_s)L_s\mathbb E_G D'/(2q_s)+(1+1/\eta_s)p_sH_sV_W .
\]
This proves both estimates. Both $|b_G|,|\eta|>.99$, and $|b_G-\eta|<.02$, since $g(j+u)$ differs from $\theta$ by at most $\sqrt{v_\theta}$ and the chord interpolation has error at most $g^2$. Their angular difference is at most twice their distance. Consequently take $d_0=C\mathbb E_GT+W+(4+A_w)J/q$ in Lemma~\ref{pay}.

\paragraph{The final strict inequality.}
Use the facet at the angle $\theta$, and write $S_h=2\tan(\pi/h)$, $c_h=\cos(\pi/h)$. The condition of Lemma~\ref{pay} on the logarithm at $z_G$ is met if $S_h/[2\rho^2(r-e_0/\rho-2e_s)^2]\le1$. If $S_h[C(2C_T+C_BF)+4+A_w]<1$, its good term satisfies
\[
 E_0/c_h\le \epsilon_c\Gamma+A(D+2Vs),\qquad
 \epsilon_c=500GV[L/(NK)]^2/c_h,\qquad A=(S_h/c_h)CF .
\]
We used (P8), the good variance identity, and $S_h<1$. As
$D\le(p/q)\mathbb E_B\|v_P-zv\|^2\le(2/q)(Y+10p/h^2)$, where $|z-e^{i\psi}|<3.16/h$, the violated side would require
\[
 (1-\epsilon_c-2VhA/K)\Gamma+(\delta-2A/q)Y<P_0p/h^2,
 \]
\[
 P_0=\pi^2+M(Z-1)+20A/q,\qquad
 M=\frac{9.872}{(1-\delta)c_h^2(r-e_s)^2}.                     \tag{P10}
\]
The formula bounds $h^2\tan^2(\pi/h)$ by $9.872$ and $2h^2(1-c_h)$ by $\pi^2$. A leaking draw costs at least $K/h^2$ in $\Gamma$, and a nonleaking bad draw costs at least $4\rho^2\sin^2((\tau-\pi/h)/2)/h$ in $Y$.

For completeness, the three final entries of each row of (P4) can be checked with the following substitutions. Put $K_+=K_0+11.85(Z_1-1)$, $p_0=\min(.15,X_*T^2/K_0)$, $q_0=1-p_0$ and $r_0=\sqrt{1-X_*/n^2}$. In $F$ use $N=n$, $h=Tn$, $V=V_1$, $K=K_0$; this upper bound decreases at larger orders because $h\le TN$. Take
\[
 e_0=\sqrt{p_0 X_*T/(nq_0)},\qquad
 e_s=\sqrt{p_0X_*/q_0}/n,\qquad
 v_\theta=\frac{Z_1X_*}{n^2q_0r_0\operatorname{sinc}^2(\tau/2)}.
\]
In the remaining trigonometric factors use $h=h_0$. With $S_h,c_h,A,C,\epsilon_c,M$ defined above, take
\begin{align*}
 S &=S_h[C(2C_T+C_BF)+4+A_w],\\
 R_{\rm nl}&=(\delta-2A/q_0)4\rho^2h_0\sin^2((\tau-\pi/h_0)/2)
              -\pi^2-M(Z_1-1)-20A/q_0,\\
 R_{\rm leak}&=K_0-\pi^2-20A/q_0-K_+\epsilon_c-2V_1(6.29)CF.
\end{align*}
On the prefix rows one obtains $4\operatorname{sinc}^2(u_*/2)/Z_1>2.75$, $\operatorname{sinc}^2\tau/Z_1>.529$, and $\lambda-V_1[Z_1-\chi_c/(1+y_0)]>.19$. All rows have $M<11.82$, $\gamma<2.34$ and $S_h/[2\rho^2(r-e_0/\rho-2e_s)^2]<.042$. In $\epsilon_c$ use $G=2Z_1/\chi_c$, $V_1$ and $K_0$. To obtain the last expression, retain the actual $K=K_0+11.85(Z-1)$ when subtracting the facet cost. Since $M<11.85$ this contributes at least $K_0$, and $hS_h/c_h<6.29$ on all rows. The signed margins in (P4) show that either type of bad draw costs more than $P_0/h^2$ in (P10). When there is no bad draw its strict inequality reads $0<0$. This proves (P1).

On returning to the physical scale, let $S_0$ be the mean squared radius on the prefix. The secant of $x^{-1}$ on $[\rho^2,1]$ gives
$Z-1\le(1-S_0)(S_0-\rho^2)/\rho^2$. The physical cap at any threshold $\rho$ no smaller than the one in (P4) is therefore
\[
 K_{\rm ph}(\rho,\lambda,H)=\max_{\rho^2+\lambda(H^2-\rho^2)\le S_0\le1}
 \left[K_0S_0+\frac{11.85}{\rho^2}S_0(1-S_0)(S_0-\rho^2)\right].    \tag{P11}
\]
If $H<\rho$ use $\rho^2$ as the lower endpoint instead. Here $\lambda$ is any known lower bound on the fraction at radius $H$ or above. To evaluate this cubic take both endpoints and its stationary points in the interval. When $\rho^2\ge11.85/(11.85+K_0)$ it is bounded by $K_0$ directly, since $S_0(S_0-\rho^2)\le1-\rho^2$. In particular the caps at $.59,.61,.70,.75$ for $d=18$, $19\le d\le36$, $37\le d\le N/4$ and $d>N/4$ are $10.73,10.53,10.06,10.10$ respectively.

\subsection{One supporting side at a radial rank}
The preceding estimate also applies when we know little about how often
a cut is crossed. We can improve it when the cut is usually preserved.
First count the average number of labels that leave. Then use the angle
of the multiplier to estimate the same loss in a chosen direction.
Return to the physical normalization $\max r_i=1$. For the first $h$ coordinates write $u_h$ for the probability of preservation and $q_h$ for the expected number of labels sent out. Then $q_h\ge1-u_h$. Summing the variance identity over this prefix gives
\[
 q_h a_h\le h^2x\,h^{-1}\sum_{i\le h}r_i^2\le X\tau^2,
 \qquad \tau=h/N.                                               \tag{R4}
\]
Indeed the radial loss $\sum_{i\le h}(r_i^2-\mathbb E r_{P(i)}^2)$ is at least $q_h(r_h^2-r_{h+1}^2)$ by matching the missing and incoming labels. Adding the nonnegative row variances gives $x\sum_{i\le h}r_i^2$.

The same marked-sum support $H(s)$ as in (P0a) gives a signed bound, with no choice of a branch for the arguments of individual draws. In a direction of twist $s$ choose input coefficients maximizing the support in its preimage direction. A preserving draw has support at most $H(s)$, while a draw sending $k$ labels out has available radius at most $L_0-k(r_h-r_{h+1})$. Write $D(s)=L_0-H(s)$. Thus
\[
 q_h(r_h-r_{h+1})\le (1-r)L_0-u_h D(s)+rD(s-h\theta).
\]
Put $\nu=\operatorname{dist}(h\theta/(2\pi),\mathbb Z)$. Choose $s$ having the sign of the shortest residue of $h\theta$ and magnitude $2\pi\nu+\pi v$, with $0\le v\le1-2\nu$. By (P0a)'s maximizing deviations,
\[
 (\operatorname{sinc}\delta_*)s/Q_0\le D'(s)=\lambda_s\le s/Q_0
       \qquad(0\le s\le\pi).
\]
Apply these two bounds respectively to $D(2\pi\nu+\pi v)-D(\pi v)$ and $D(\pi v)$. The same gap inequality with no signed terms follows by summing $r r_i\le\mathbb E r_{P(i)}$. Interpolation yields, for every $0\le c\le h^2/(L_0Q_0)$,
\[
 q_h\,2h(r_h-r_{h+1})\le \frac{L_0}{h}
   [J\tau^2-u_h\gamma\nu(\nu+v)+(1-u_h)B_1v^2],\quad
 J=\frac{2X}{1+r},\quad
 \gamma=4\pi^2 c\operatorname{sinc}\delta_*,\quad B_1=\pi^2 c.       \tag{R5}
\]
The bound $c\le4\rho/(1+\rho)^2$ suffices when the prefix radii are at least $\rho$, by averaging $r_i+\rho/r_i\le1+\rho$ and applying the arithmetic--geometric mean inequality. On these ranks $h\rho>100$ (use the individual input rows), so $\operatorname{sinc}\delta_*>.9998$. We may increase $J$ to $\eta X$, $\eta=1.00005$, since $r>.99998$.

For a radial interval $[a,b]$ put $T(x)=\min(b,\max(a,x))$ and
$\tilde a_h=h[T(r_h)^2-T(r_{h+1})^2]$. This gap is at most both $a_h$ and $2bh(r_h-r_{h+1})$. In particular (R5) applies to $q_h\tilde a_h$ with the extra factor $b$ in front of its bracket, since $L_0/h\le1$ and the bracket is nonnegative. If the run starts at the last radius at least $b$ and ends at the last at least $a$, the telescoping sum of the physical squared gaps is $\sum_h\tilde a_h/h=b^2-a^2$.

\paragraph{A cap among the earlier radial ranks.}
Let $h$ be a radial rank, $\tau=h/N$, and let $w_h=N[T(r_h)^2-T(r_{h+1})^2]$, where $T(r)=\min(b,\max(a,r))$ truncates the radii to $[a,b]$. Only the part of a gap inside this interval contributes. Write $u=1-(d-1)/n$, where $2\le d\le36$, $N\ge n$, and let $t\in(0,1)$ measure position in the Farey interval from the endpoint of denominator $d$. Put $v=\min(t,1-t)>0$. Equations (R3), (R5) give
\begin{align*}
 J&\le C v(1-uv/d)/d,\qquad C=4\pi^2\xi\eta/u^2,\\
 p w_h&\le(b/\tau)\{J\tau^2-(1-p)\gamma\nu(\nu+s)+p B_1s^2\},
       \qquad 0\le s\le1-2\nu,\qquad
 p w_h\le bJ\tau,\\
 \nu&=\operatorname{dist}(h\arg(z)/(2\pi),\mathbb Z),\qquad
 \gamma=4\pi^2 c\operatorname{sinc}\delta_*,\quad B_1=\pi^2 c.
\end{align*}
Here $p$ is the probability of leaking at the given cut, and $\xi,\eta\ge1$ bound the elementary finite-order errors. The signed bound follows first for the mean number of leaking labels in place of $p$ on the left and hence also for $p$. When $\tau\le u$, the denominator ratio $y=e/N$ of the other Farey endpoint satisfies $y\ge u\ge\tau$. The distance from a grid residue to a multiple of $d$ after the shift $\tau t/y$ is at least $\min(\tau t/y,1-\tau t/y)\ge\tau v$. Thus $\nu\ge\tau v/d$ \emph{without} an additive error. Since $\gamma\le4B_1$ and $\nu(1-\nu)$ increases for $\nu\in[0,1/2]$, the signed bound at $s=1-2\nu$ gives
\[
 p w_h\le b\left[\frac{pB_1}{\tau}
       +(C\tau-\gamma)\frac{v}{d}
       -(Cu-\gamma)\frac{\tau v^2}{d^2}\right].
\]
In the range $\tau\le\min(u,\gamma/C)$ the coefficient $Cu-\gamma$ is positive, so $p=0$ is impossible. Average the last estimate with $p w_h\le bC\tau v/d$, giving the former weight $C\tau/\gamma$. The linear term cancels and division by $p$ proves
\[
                  w_h\le bB_1 C/\gamma.                              \tag{R6}
\]
In particular no packet-boundary surcharge is needed where this individual bound suffices. At a given radial interval, a stronger choice of $c$ is available from any fraction $\lambda\ge1/2$ of the prefix whose radii are at least $H$, the other radii being in $[a,1]$. Namely
\[
 \frac{LQ}{h^2}\le 1+\lambda(1-\lambda)(1-a)^2/a+(1-H)^2/(4H).
\]
The product of mean radius and mean reciprocal radius is separately convex in each radius. At a maximum, the low coordinates can therefore be $a$ or $1$, the required high ones $H$ or $1$. In the pair expansion of $LQ/h^2-1$, pairs incident to a coordinate at $a$ contribute at most $\lambda(1-\lambda)(1-a)^2/a$ and pairs at $H,1$ contribute at most $(1-H)^2/(4H)$. One may take the reciprocal of this upper bound as $c$.

\subsection{Several consecutive cuts}\label{app:packet}
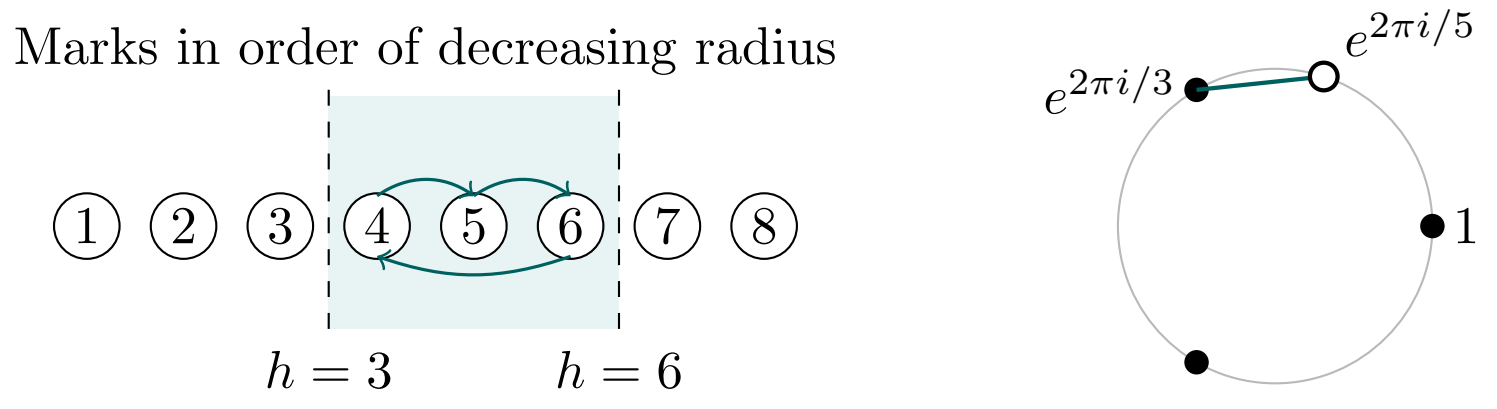
\begin{figure}[htbp]
\centering
\begin{tikzpicture}[font=\small]
\begin{scope}[x=.64cm,y=1cm]
\fill[teal!9] (3.5,-.68) rectangle (6.5,.86);
\foreach \x in {3.5,6.5} \draw[dashed] (\x,-.68)--(\x,.93);
\foreach \j in {1,...,8} \node[fill=white,draw,circle,inner sep=1.5pt] at (\j,0) {$\j$};
\draw[->,teal!75!black,semithick] (4,.2) to[bend left=35] (5,.2);
\draw[->,teal!75!black,semithick] (5,.2) to[bend left=35] (6,.2);
\draw[->,teal!75!black,semithick] (6,-.2) to[bend left=19] (4,-.2);
\node at (4.5,1.16) {Marks in order of decreasing radius};
\node[below] at (3.5,-.69) {$h=3$}; \node[below] at (6.5,-.69) {$h=6$};
\end{scope}
\begin{scope}[xshift=8.5cm]
\draw[gray!55] (0,0) circle (1.04);
\foreach \ang in {0,120,240} \fill (\ang:1.04) circle (2.3pt);
\draw[thick,teal!75!black] (72:1.04)--(120:1.04);
\draw[fill=white,thick] (72:1.04) circle (2.6pt);
\node[right] at (0:1.04) {$1$};
\node[above right] at (72:1.04) {$e^{2\pi i/5}$};
\node[left] at (120:1.04) {$e^{2\pi i/3}$};
\end{scope}
\end{tikzpicture}
\caption{A permutation that fixes both cuts has to pay for the
rotation error on the enclosed block. The three black points are
the eigenvalues of the illustrated three-cycle. The segment ends at
the nearest one to the chosen fifth root.}\label{fig:block}
\end{figure}

Two preserved cuts enclose a block. The point of grouping $d$ consecutive
cuts is that the block between any two of them has fewer than $d$ sites.
Consequently it cannot carry the root of unity at the chosen Farey
endpoint. We estimate the squared error of this block before we sum
the gaps.
For consecutive cuts starting and ending in the indicated radial range, start packets of $d$ consecutive
cuts at the first cut, allowing a shorter final packet. If a particular
permutation preserves $m_P$ cuts of a packet, the $(m_P-1)_+$ blocks
between consecutive preserved cuts are disjoint invariant blocks with
$s<d$ labels. These blocks are disjoint over all packets.
For example, take $d=5$ and suppose the cuts after labels $3$ and $6$
are preserved. Then the three labels $4,5,6$ map into themselves.
Even if the restriction is a three-cycle, it cannot rotate a nonzero
Fourier mode through one fifth of a turn. Figure~\ref{fig:block}
shows this particular restriction. The estimate also covers the
identity and a transposition on those same three labels.

On such a block the permutation is unitary, with eigenvalues among
roots of unity of orders at most $s$. Each of these differs from the
Farey endpoint of denominator $d$ by at least $1/(sd)$ turns, whereas
$z$ is at most $1/(de)$ turns from that endpoint. Hence the
block's squared error is at least $s\rho^2\delta_s^2$, where
\[
 \delta_s\ge 2\sqrt r\sin\bigl(\pi(1-s/e)/(sd)\bigr).
\]
If $sd<6$ the block size is one and $sd\delta_s>3.96$ by the
chord on $[0,\pi/2]$. For $6\le sd<10$ it is at most two and the
chord on $[0,\pi/6]$ gives $sd\delta_s>5.78$. For $sd\ge10$ the
chord on $[0,\pi/10]$ gives
$sd\delta_s\ge20\sqrt r\sin(\pi/10)(1-(d-1)/e)$.
In particular all block sizes with $sd\ge6$ use $c_*=5.78$ at $d\le18$, and $c_*=5.70$ at $19\le d\le36$,
since $e>N-d$; in $19\le d\le36$ use the chord at $\pi/19$, giving $38\sqrt r\sin(\pi/19)(1-35/415)>5.70$. Other constants $c_*$ at larger denominators follow below.
The total expected squared error in the
original normalization is $x\sum_i r_i^2\le X/N$. Consequently,
with $e_j=\mathbb E(m_P-1)_+$ in packet $j$,
\[
 \sum_{h\text{ in packet }j}u_h\le1+e_j,\qquad
 \sum_j e_j\le e_*=
 \frac{X}{N\rho^2}\max\left\{\frac{d^2}{3.96^2},
                       \frac{d^2(d-1)}{c_*^2}\right\},     \tag{R8}
\]
\paragraph{A rank-dependent baseline.}
The two bounds $a_h\le K$ and $a_h\le X\tau^2/(1-u_h)$ imply
$a_h\le X\tau^2+Ku_h$. This is the chord of the increasing convex
function up to its first meeting with $K$, followed by the constant
cap. For a run from hull rank $b$ ending at or before $N-1$, set
$\ell\le b/N$. In a packet starting at $c$ the contribution of its
$u_h$ terms to the squared radial span is at most $K(1+e_j)/c$.
The term $X\tau^2$ allows a gap even when the preservation probability is zero. This baseline contributes at most
$X\sum_{h=b}^{N-1}h/N^2\le X(1-\ell^2)/2$. If $b>d/2$,
convexity of $1/x$ gives
\[
 \sum_{j=0}^{\lfloor(N-1-b)/d\rfloor}\frac{1}{b+jd}
 \le\frac1d\log\frac{N-1+d/2}{b-d/2}.
\]
Thus the complete radial span is at most
\[
 \frac X2(1-\ell^2)+\frac Kd\log(1/\ell)+
 K\left[\frac{1/2-1/d}{N}+\frac{1}{2(N\ell-d/2)}\right]
 +\frac{K e_*}{N\ell}.                                   \tag{R9}
\]
We used $\log(1+x)\le x$ and $-\log(1-y)\le y/(1-y)$.

\subsection{Denominators of at least eighteen}\label{app:large-d}
Let a packet begin near a rank fraction $t$. Its unsigned gap estimate
has a baseline proportional to $Xt^2$. The preserving probability can
raise a few cuts above that baseline, up to the cap from (P11). The
block estimate (R8) controls the total excess probability. We compare
the largest total squared span permitted by these bounds with the
distance from a hull vertex down to the first smaller mark.
We first record the concentration inequality which will
also control the terminal fragment of a run. Suppose $K>2X$ and put $w=1-X/K$. For a cut with preservation
probability $u$ its excess over $X$ is at most
$\min(K-X,X/(1-u)-X)$. Up to $u=w$ the function is convex.
At fixed total preservation mass $s$, repeatedly transfer mass
between two entries until one is zero or the other is $w$.
The excess is at most
\[
 F(s)=\lfloor s/w\rfloor(K-X)+
 \frac{X}{1-s+w\lfloor s/w\rfloor}-X .
\]
This also holds for $s$ larger than the possible total capped mass
by monotonicity. We have $F(1)=K-2X+X/w$, and for $s\ge1$,
$F(s)\le F(1)+(K/w)(s-1)$. Indeed the line agrees with $F$ at
$1,2w$ and dominates the convex function between these points.
In each subsequent complete period, the line gains $K$ while $F$
gains $K-X$. From (R8), full packets and the terminal fragment
therefore give the valid pair
\[
            \bar K=X+F(1)/d,\qquad E_*=F(1)+(K/w)e_* .           \tag{R14}
\]
For $18\le d\le36$ retain the actual rank in the squared deficit.
For a packet starting at rank $c$ put
$W_-=X(c/N)^2$, $W_+=X\min((c+d-1)/N,1)^2$.
A cut with preservation probability $u_h$ has excess over its baseline
$X(h/N)^2$ bounded by
$\min(K-W_-,W_+u_h/(1-u_h))$.
Apply the concentration argument (R14) with deficit $W_+$ and
cap $K'=K+W_+-W_-$. Mass one costs at most
\[
    K-W_-+W_+^2/(K-W_-).
\]
Each unit of further mass costs at most
$(K+2Xd/n)^2/(K-X)$, since $K'/[1-W_+/K']=K'^2/(K-W_-)$.
Divide these costs by the packet's starting rank and sum the
baseline as in (R9). For a group $a\le d\le c$ let $u=1-(c-1)/n$, $D=2(9.86961)\xi(1-u/(2c))/u^2$, $\ell=.99945/[\xi(1.5/u-1/(2c))]$, $X_0=D/a$. Estimate (R8) by $e_*\le D\max(c/3.96^2,c(c-1)/5.70^2)/(n\rho^2)$. We may start packets at $\lceil\ell N\rceil$, counting zero for the gaps before the last hull-vertex rank. The preserved cut budget is still valid since all intermediate labels have radius at least $\rho$. The preceding estimate without the early-rank cap is
\begin{align*}
 U_*={}&\frac{X_0(1-\ell^2)}2+
 \left(K-X_0\ell^2+\frac{X_0^2}{K-X_0}\right)
 \left(\frac{\log(1/\ell)}a+\frac{1/2-1/a}{n}
 +\frac1{2(n\ell-c/2)}\right)
 \\ &+\frac{(K+2cX_0/n)^2e_*}{(K-X_0)n\ell}.
\end{align*}
In the harmonic sum keep $(\log(1/\ell)-1/n)/d$ together when replacing $d$ by its lower bound $a$. The numerator is positive in all six groups. It suffices to take $X=X_0$ in this bound and its proof: replacing all unsigned deficits by $X_0\tau^2$ only increases their bounds.

A packet starting at $\lceil\ell N\rceil+jd$ has every rank fraction at most $t_j=\ell+(j+1)d/N$. Use $\rho=.61$, $K=10.53$ from $19$ to $36$, and $\rho=.59$, $K=10.73$ at $18$. Put
\[
 K_s=\frac{9.86961\xi(1.00005)}{.9997u^2},\quad
 t_* =\min\left(\frac K{K_s},
 \frac{4\rho}{(1+\rho)^2}\frac{.9997u^2}{\xi(1.00005)}\right).
\]
If $t_j\le t_*$ then (R6), with the actual $u_d\ge u$, bounds the physical gap $a_h$ by $K_st_j$ throughout that packet. Here take the upper truncation $b=1$ and no lower truncation. The reciprocal-moment factor $4\rho/(1+\rho)^2$ uses only the actual prefix radii, all at least $\rho$. The decrease in its baseline-excess charge, divided by its starting rank, is at least
\[
 \frac{1}{N}\left(1-\left(\frac{X_0}{K_s\ell-X_0}\right)^2\right)
 (K-K_st_j)/t_j .
\]
This follows by differentiating $K-W_-+W_+^2/(K-W_-)$ in $K$ and using $W_+\le X_0$, $W_-\le X_0$. Additional preservation mass can still use the original coefficient, since $K_s\ell>2X_0$. Put $t_b=\ell+c/n<t_*$. Integrating the decreasing function $(K-K_st)/t$ up to $t_*$, and charging each interval to the grid point immediately on its left starting from $t_0=\ell+d/N$, gives a total saving of at least
\[
 D_* = \frac{1-(X_0/(K_s\ell-X_0))^2}{c}
      [K\log(t_*/t_b)-K_s(t_*-t_b)] .
\]
If the run ends early, fill the rest of the list with zero gaps and zero
preservation probabilities for the purpose of summing these upper bounds.
Such terms satisfy the cap, baseline and excess bounds for a slot, and
contribute zero to the actual span. The improved cap is used only on
packets wholly below its threshold. Direct substitutions give
\[
\begin{array}{c|rrrrrr}
 [a,c]&18&19\!:\!20&21\!:\!23&24\!:\!28&29\!:\!33&34\!:\!36\\\hline
 U_*-D_*<&.629&.619&.607&.613&.607&.596
\end{array}
\]
The required squared span exceeds $.6464$ at $18$ and $.6215$ on the other rows.

For $d\ge37$ and $s=d/N\le b=.065$ use (R9) and the cap $K=10.06$ at $\rho=.70$. Put
\[
 \ell=\frac{.99945(1-b)}{1.5\xi},\qquad D=\frac{2\pi^2\xi}{(1-b)^2},
 \qquad c_*=20(.9999)\sin(\pi/10)(1-b/(1-b)).
\]
Here $\tau\ge\ell$, $dX\le D$, and (R8) uses only the case in which the product of block size and $d$ is at least ten. By (R9) the span is at most
\[
 \frac{D(1-\ell^2)}{74}+\frac K{37}\log(1/\ell)
 +\frac{K}{2n}\left(1+\frac1{\ell-b/2}\right)
 +\frac{KD b^2}{\ell\rho^2 c_*^2}<.441.
\]
This is below $(.998)^2-.70^2$.

\paragraph{Mixed Farey ranks.}
Write $s=d/N$, $y=e/N$ and $m=s\lfloor1/s\rfloor$, and put
\[
 k=4\pi^2\xi,\qquad \eta=1.00005,\qquad
 \gamma=38.18,\quad B_1=9.56,\quad K=10.06,\quad C_0=k\eta/\gamma .
\]
Here $\xi=1.00020$ and $n=450$. For $37\le d\le N/4$ use $\rho=.70$ and $c=.9675$ in (R5). Above $N/4$ use $\rho=.75$, $c=.9785$, $\gamma=38.6$, $B_1=9.66$ and $K=10.10$. In either case $4\rho/(1+\rho)^2\ge c$, $4\pi^2c\operatorname{sinc}\delta_*>\gamma$ and $\pi^2c<B_1$. The squared gap is bounded by $K/h$. We write $\tau=h/N$ and $a_h=h(r_h^2-r_{h+1}^2)$.

After reflection, let $p/d<q/e$ be the adjacent endpoints. Write the pair consisting of the rank $h<N$ and the nearest integer to $h\arg(z)/(2\pi)$ as $i(d,p)+j(e,q)$. This is an integral combination since $dq-ep=1$. A rank is called pure if one coefficient vanishes and mixed if both are nonzero. The coefficients cannot both be positive since $d+e>N$. Thus the two endpoint errors have the same sign after multiplication by the respective coefficients. At fraction $t$ of this Farey interval measured from $p/d$, the distance $\nu$ to the closest integer therefore obeys
\[
 \nu=|i|t/e+|j|(1-t)/d\ \ge XD/k,\qquad D=|i|ms+|j|y^2 .
\]
The inequality follows from $dX/k\le\min(t/(my),(1-t)/y^2)$.

Call a rank exceptional if $\gamma\nu(1-\nu)\le J\tau^2$. Choose $e_0<1/2$ with $\eta X/\gamma\le e_0(1-e_0)$ and put $C=C_0/(1-e_0)$. At an exceptional rank $\nu\le e_0$, hence $D\le C\tau^2$. Else take the signed twist just above $[J\tau^2/(\gamma\nu)-\nu]_+$ and then pass to the limit. It lies below $1-2\nu$ by the strict reverse of the exceptional condition. If $\nu=1/2$ the only possible twist is zero and (R5) at twist zero gives the same conclusion. Since $q_h\ge1-u_h$, this gives
\[
 \frac{a_h}h\le\frac{J\tau}{N}
 +\frac{B_1}{N\tau}\left[C_0\tau^2/D-XD/k\right]_+^2.       \tag{F1}
\]
At $\nu=1/2$ and nonexceptionality, $J\tau^2-\gamma\nu^2<0$, so the putative case $u_h=1$ cannot occur. The same reasoning for general nonexceptional $\nu$ uses $J\tau^2-\gamma\nu(\nu+t_1)<0$ before taking the limit. On combining (F1) with the individual cap for an exceptional rank one obtains at every rank
\[
 \frac{a_h}h\le \frac{J\tau}{N}
                  +\frac{KC^2\tau^3}{ND^2}.                    \tag{F2}
\]

We give the sum of (F1) on $.065\le s\le.25$. On a subinterval $a\le s\le b$ put
\[
 l=\frac{.99945(1-b)}{1.5\xi},\quad m_0=y_0=1-b,\quad
 X_b=\frac{2\pi^2\xi}{(1-b)^2\max(37,na)},\quad
 e_0=\frac{1-\sqrt{1-4\eta X_b/\gamma}}2 .
\]
Thus $l\le\tau$, $X\le X_b$. There are two families of mixed coefficients at $h>0$. For $j=-v<0$, $D=m\tau+v y(m+y)$ and the possible ranks at fixed $v$ have spacing $s$ in units of $N$. For $i=-v<0$ use $D=y\tau+vs(m+y)$ and spacing $y$. All such ranks are nonexceptional. For the first family $m+y(m+y)\ge(1-s)+2(1-s)^2\ge1.875$; for the second $y+s(m+y)\ge(1-s)(1+2s)>1.0565$ for $.065\le s\le.25$. Since $C<1.0522$ in the table, both imply $D>C\tau^2$.

Use the triples $(c,M,H)=(a,m_0,2y_0^2)$ and $(y_0,y_0,2ay_0)$ for the two families. They bound below the spacing and the coefficients in $D\ge M\tau+vH$. The successive values of $\tau$ at a fixed $v$ are bounded above by $u=1-jc\ge l$, $j=0,1,\ldots$. The positive part in (F1), divided by $\sqrt{\tau}$, increases with $\tau$ wherever it is positive: its expression at the lower bound on $D$ is
$\sqrt{\tau}[C_0/(M+vH/\tau)-(X/k)(M+vH/\tau)]_+$.
For fixed $u$, its square is convex decreasing in $v$ on $v\ge1/2$, since $A/t-Bt$ is convex and its squared positive part is convex. The sum over $v=1,2,\ldots$ is bounded by its integral from $1/2$. With $d_0=Mu+H/2$ and $w=(d_0/u)\sqrt{x/(C_0k)}$ the resulting charge is
\[
 S(c,M,H;x)=B_1\sum_{u=1-jc\ge l}
 \frac{C_0^2u^3}{Hd_0}
 \begin{cases}1-8w/3+2w^2-w^4/3,&w\le1,\\0,&w>1.\end{cases} \tag{F3}
\]
Indeed the integral reduces by substitution to that of $(v^{-1}-v)^2$ from $w$ to $1$, with the limiting value at $w=0$ as displayed.

For the second family use a sharper enumeration. Fix its positive coefficient $j$ in front of $e$, writing $i=-v<0$. Then $\tau=jy-vs$ and $D=jy(m+y)-m\tau$. Put $H=y_0(m_0+y_0)$. For $j=1$, each $v\ge1$ has $\tau\le u_v=1-va$ and $D\ge d_v=y_0^2+vam_0$. Its total cost is at most
\[
 T_1(x)=B_1\sum_{u_v\ge l}u_v^{-1}
                          [C_0u_v^2/d_v-xd_v/k]_+^2.
\]
For every fixed $j\ge2$, order its ranks downwards, so the respective bounds are $u=1-ra\ge l$, $r=0,1,\ldots$, with $D\ge jH-u$. For each $u$ integrate the decreasing convex function of $j$ from $3/2$, exactly as for (F3). Thus the second family costs at most
\[
 T_2(x)=T_1(x)+B_1\sum_{u=1-ra\ge l}
 \frac{C_0^2u^3}{H(3H/2-u)}
 \bigl[1-8w/3+2w^2-w^4/3\bigr]_{w<1},
 \]
where
\[
 w=(3H/2-u)\sqrt{x/(C_0k)}/u .
\]
Here the bracket is set to zero at $w\ge1$. Its endpoint charge is convex as a function of $x$ by the integral representation.

Apart from the one possible rank $e$, charge only pure multiples $h=jd$ with
$j b\ge l$, $j a<1$. Put $u=\min(1,jb)$, $v=\max(l,ja)$ and recall $D=m\tau$. When $u<m_0/C$, (F1) bounds the excess after multiplication by $N$ by $B_1u(C_0/m_0-xm_0/k)^2$. For every other multiple use
\[
 \min\{K/v,\ KC/m_0,\ KC^2u/m_0^2\}.
\]
The first and last bounds follow from the cap and (F2); use the last when $\tau\le m_0/C$ and the first with its actual denominator $\tau$ otherwise to obtain the middle bound. Include $K/y_0$ for $h=e$ and call the sum of these charges $P(x)$.

On multiplying the total span by the fixed order $n\le N$, the preceding estimates give the bound
\[
 n\eta x(1-l^2)/2+P(x)+S(a,m_0,2y_0^2;x)+T_2(x),
 \qquad x=X\in[0,X_b].                                       \tag{F4}
\]
The baseline uses $\sum_{h\ge lN}^{N-1}h/N^2\le(1-l^2)/2$. The expression in (F4) is convex in $x$: each integrated charge is a squared positive part of an affine function before the integral is evaluated. Thus its maximum occurs at $0$ or $X_b$. The table gives upper bounds at both endpoints, and the respective summands at $X_b$. The index $j$ in (F3) runs only through $\lfloor(1-l)/c\rfloor$.
\[
\begin{array}{c|c|r|rrrr|r}
 [a,b]&j\text{ in }P&(F4)_0&
  \text{baseline}&P&S_1&T_2&(F4)_{X_b}\\\hline
{[.065,.075]}&9\text{--}15&136.35&87.05&81.32&6.39&32.89&207.64\\
{[.075,.085]}&8\text{--}13&123.50&90.13&71.37&6.21&30.48&198.17\\
{[.085,.1]}&6\text{--}11&118.99&91.85&71.31&5.91&27.88&196.94\\
{[.1,.15]}&4\text{--}9&132.13&92.83&79.04&7.24&31.68&210.77\\
{[.15,.25]}&2\text{--}6&137.15&87.80&75.10&9.96&39.13&211.97\\
\end{array}
\]
For instance, if $.15\le d/N\le .25$, the last row bounds the
two endpoint totals by $137.15$ and $211.97$. Convexity places the
entire graph of (F4) below the line joining its endpoint values.
Consequently every intermediate total is at most $211.97$.
The separate contributions are rounded up individually, so one uses
the total computed before those roundings. Across the five rows
the largest bound is below $212$, whereas the required span satisfies
$450[(.998)^2-.70^2]>227$.

For $s>1/4$ take $C=1.04$ in (F2); this is permitted since $X\le.470$. The squared span on $[.25,.5]$ or $(.5,1]$, multiplied by $N$, is bounded by the respective rows
\[
\begin{array}{c|cc|rrrr}
s &(c,M,H)_{j<0}&(c,M,H)_{i<0}&
\text{baseline}&\text{pure}&j<0&i<0\\\hline
{[.25,.5]}&(1/4,2/3,3/4)&(1/2,1/2,3/8)&57.35&55.60&23.58&50.70\\
(.5,1]&(1/2,1/2,1/2)&(1/2,1/2,1/2)&62.20&40.4&34.60&34.60
\end{array}
\]
For the first baseline use $\tau\ge l=.428$ and $s(m+y)y\ge9/32$: below $1/3$ this follows from $m\ge3s$, $s(1-s)(1+2s)\ge9/32$, and above it from $m=2s$, $s(1-s^2)\ge8/27$. For the second use $l=.461$ and $s(m+y)y\ge1/4$. In both cases $NJ\le k\eta/[s(m+y)y]$. Below $s=1/3$ charge the second and third multiples of $d$ and $1/y\le1.5$ in units of $K$; at $s=.25$ the fourth multiple is the full order and is absent. When $1/3<s<.428$ only the second multiple and $1/y$ remain. Up to $s=.5$, charging all possibilities costs at most $K(1.5/.428+2)<5.51K$. Above $.5$ the total pure cost is at most $4K$. Each mixed column is bounded from (F2) by
\[
 KC^2\sum_{u=1-jc\ge l}\frac{u^3}{H(Mu+H/2)}
\]
using convexity of $(Mu+vH)^{-2}$ once more. For the triples in the first row use $m\ge2/3$, $y\ge1/2$, $y(m+y)\ge.75$, $s(m+y)\ge3/8$; above $.5$ each entry exceeds the corresponding half. The totals are less than $187.23$ or $171.81$, whereas the required span is at least $450[(.9984)^2-.75^2]>195.43$.
\subsection{How many marks lie between two radii?}\label{app:counts}
For the remaining denominators we need more than a bound on the total
span. If consecutive cuts each allow little change in the squared radius,
a fixed radial interval takes many cuts to cross. The number of such
cuts gives a lower bound for the mass available to the angular test.
A packet consists of up to $d$ consecutive eligible cuts. We use this count on seven specified pairs of radii, some of which overlap. For any one interval $[a,b]$, apply its constants to gaps truncated to that interval, and write $w_h=N[T(r_h)^2-T(r_{h+1})^2]$. We give numbers $\mu,E_p$ for which any run of $k$ such cuts has $\sum w_h\le k\mu+E_p$. Thus the fraction of labels in $[x,y)$, for $a\le x<y\le b$, is at least
\[
 C(x,y;\mu,E_p)=\max\{0,[y^2-x^2-(\mu+E_p)/n]/\mu\}.           \tag{L5}
\]
To see the extra $\mu$ in (L5), suppose $m$ marks lie in $[x,y)$.
Starting from the last radius at least $y$, it takes $m+1$ gaps to
reach the first radius below $x$. A smaller radius exists by (P1),
and the hull reaches above $y$. After truncation the $m+1$ drops
telescope to $y^2-x^2$. Therefore
$N(y^2-x^2)\le(m+1)\mu+E_p$. Divide by $N\mu$ and use $N\ge n$
to obtain (L5). The term $E_p$ is paid once for the whole run.

Use the actual $d$, $u=1-(d-1)/n$, $X_0=X_d$ and let $(\ell,H)$ be the high population specified in the angular tables below. Use $e_*$ from (R8) at $\rho=a$, $N=n$, $X=X_0$ and $c_*=5.78$. Evaluate $K_{\rm ph}$ from (P11), rounding up to three decimal places. Partition the possible smallest rank fraction of a packet at $\ell,.76,.80,.85,.88,.90,.92,.94,.96,.98,u,1$, sorted and with internal points below $\ell$ omitted. For any resulting interval $[A,T]$ put $U=\min(1,T+(d-1)/n)$, $\lambda=\ell/U$, and
\begin{gather*}
 c=[1+\lambda(1-\lambda)(1-a)^2/a+(1-H)^2/(4H)]^{-1},\quad
 \gamma=4(9.86960)(.9998)c,\quad B_1=9.86961c,\\
 C_1=4(9.86961)\xi\eta/u^2,\quad y=bU\eta X_0,\quad
 k=K_{\rm ph}(a,\lambda,H)/A,\qquad \eta=1.00005.
\end{gather*}
If $U\le\min(u,\gamma/C_1)$ replace $k$ by $\min(k,bB_1C_1/\gamma)$, using (R6). For $d\ge6$ put $q=1-y/k>1/2$, $F_0=k-2y+y/q$. By (R14) each slot costs at most $y$ plus an excess bounded by $F_0$ per packet and $k/q$ per extra unit of preservation mass. Thus
\[
 \mu=\max_{[A,T]}(y+F_0/d),\qquad
 E_p=\max_{[A,T]}(y+F_0)-\mu+e_*\max_{[A,T]} k/q.             \tag{L6}
\]
The first slot of any terminal fragment already pays $\mu$ instead of $y$. For a group take the separate maxima over the allowed denominators.

For $2\le d\le5$, let $v_0=\min(t,1-t)$ for the Farey parameter $t$. Set $e=0$ if $U\le u$, and otherwise $e=d/(n-d)$. Then $\nu\ge\tau(v_0-e)_+/d$ by (R6) on the early ranks. On all other ranks the same grid shift $\tau t/y_F$ uses $y_F=e_F/N$ for the complementary denominator, with $1/y_F-1<d/(n-d)$. The distance to an integral multiple of $d$ is at least the smaller of $\tau t/y_F$ and $|1-\tau t/y_F|$; when $\tau t/y_F>1$ one has $1-t<d/(n-d)$. Put
\[
 M=\gamma[(1/2-e)/d]^2,\qquad
 D(U)=\max_{0\le v\le1/2}\left[\frac{C_1U-\gamma}{d}v
       -\frac{(C_1u-\gamma)U}{d^2}v^2\right]+\gamma e/d .       \tag{L7}
\]
Evaluate the quadratic at $(C_1U-\gamma)d/[2(C_1u-\gamma)U]$ clamped to $[0,1/2]$. For $p=1-u_h$ the unsigned twist in (R5) gives $p w_h\le b\tau[\eta X_0-(1-p)M]$: the bracket from (R3) with variable $v_0$ increases up to $1/2$, by the derivative lower bound $[4\pi^2(1-1/d)-\gamma/d]/d>0$. The signed twist $1-2\nu$ gives $pw_h\le b[pB_1/\tau+D(\tau)]$ as in (R6), replacing $\nu$ with $\tau v_0/d$ at cost at most $\gamma e/d$ after division by $\tau$. Here $D(\tau)$ is increasing. Set $m=bUM$, $b_0=bB_1/A$, $g=bD(U)$. Each cost is at most
\[
 f(p)=\min\{k,m+(y-m)/p,b_0+g/p\},\qquad f(0)=k.               \tag{L8}
\]
Throughout $f(1)=y$ and $g,y-m\ge0$. For $1\le j\le d$ the total of $j$ slots is bounded by $S_j$ defined as follows. Put
\[
 \mathcal Q=\{0,1,1-g/(k-b_0),1-(y-m)/(k-m),
                        1-(y-m-g)/(b_0-m)\}\cap[0,1],
\]
discarding undefined entries. Choose $j-1$ (possibly repeated) members $q_i$ of $\mathcal Q$ of sum $s\le1+e_*$ and set $q_j=\min(1,1+e_*-s)$. Let $S_j$ be the largest $\sum_{i=1}^j f(1-q_i)$ among these choices. Indeed the preserving masses have total at most $1+e_*$ by (R8); $f(1-q)$ increases and is convex between its switches. With any two interior masses, holding their sum fixed and varying one to a boundary of its current interval cannot decrease the objective if the better endpoint is chosen. Iteration proves the rule, including $q_i=1$ by the endpoint cap.

The verifier enumerates the entire finite set with repetitions allowed, including every terminal fragment. Now use $\mu=\max S_d/d$ and $E_p=\max_{[A,T],1\le j\le d}(S_j-j\mu)$. Full packets and their terminal fragments prove (L5). In all rows the separate resulting bounds $\mu,E_p$ are rounded upwards to three decimals.

\subsection{A convex test for a small denominator}\label{app:small-d}
A count of interior marks on its own gives no contradiction: a regular
polygon may contain many of them. We will construct a convex function
that takes advantage of both their radii and their angular positions.
At selected marks we give the function a positive slope in the direction of rotation.
The inequalities below compare the gain along the rotation with the loss
due to the contraction.
It remains to exclude $2\le d\le17$. After conjugation $z^d=Re^{i\varphi}$ with $\varphi>0$. The same configuration realizes this return. Consequently every convex function on the plane satisfies
\[
                       \sum_i f(z^dv_i)\le \sum_i f(v_i).          \tag{A1}
\]
We first choose the values of a convex function as a function of radius. At a source point, keep the physical plane fixed, lift each mark to its chosen height, and tilt a supporting plane in the forward tangential direction until it touches a second lifted mark. At that point its slope records a lower bound on the increase of the function under the rotation. A mark well inside the hull can obstruct a tilt only when its height is close to the plane. Varying the radius at which the outer slope begins bounds the total effect of such marks. Two ordinary convex functions supply additional slopes at the interior radii.

Write $u=1-(d-1)/n$ for a single denominator. In a group use its largest denominator to bound $u$ below. Put $\eta=1.00035$. The two endpoint estimates, now applied to $E=N(1-R)/(\pi R\varphi)$ and $F=N\varphi/(2\pi)$, give
\[
 E\le\frac\eta u\min\left(1-\frac{ut}{d},\frac{1-t}{t}\left(1-\frac{1-t}{d}\right)\right),
 \qquad F\le t/u.                                                \tag{A4}
\]
Indeed, (R3) bounds $N(-\log R)/(\pi\varphi)$ with factor $\xi$, and $\xi/R<\eta$. Also $\xi\varphi/(R\sin\varphi)<\eta$ and $\tan(\varphi/2)/(\varphi/2)<\eta$ throughout this range. In each of the cases below put $q=\min(t_{\max},1/2)$. The angular cost of a label of radius $r_i$ is expressed by $EJ_E+FJ_F$ after averaging over the labels. Here $J_F$ is the mean of the active tangential branch's radial slope times $r_i$, with zero on the ordinary radial branches. The angular heights increase this coefficient to $J_E$ on an ordinary angular branch, and $J_E$ includes the radial coefficient on every radial branch. The construction gives $1\ge J_E\ge J_F\ge\beta h_*$, where $\beta N$ hull vertices have radius at least $h_*$. In particular $J_F/u\ge\eta J_E/d$ on all the rows below.
When $t\le1/2$, the first affine estimate for the combined cost increases with $t$. Above $1/2$, the second bound for $E$ decreases, and that for $E+F$ is at most its midpoint value: the derivative after addition of $t$ to the second expression in parentheses is $(1-1/d)(1-t^{-2})$. Decomposing $EJ_E+FJ_F$ into its radial and angular parts gives
\[
 EJ_E+FJ_F\le E_*(J_E-J_F)+B J_F,\qquad
 E_*\ge\eta(1/u-q/d),\quad B\ge\eta(1/u-q/d+q/u).             \tag{A5}
\]
Likewise $N(1-R\cos\varphi)/(\pi R\sin\varphi)\le B$: the two terms are $(1-R)/(R\sin\varphi)$ and $\tan(\varphi/2)$. The same midpoint estimate applies with their respective finite factors bounded by $\eta$.

For $2\le d\le4$ use range I for $0<t\le.4$ and range II for $.4\le t<1$. Put $u=1-(d-1)/n$, $\eta=1.00035$, $\xi=1.00020$, $q=\min(.5,t_{\max})$ and
\[
 E_* = \eta(1/u-q/d),\quad B=E_*+\eta q/u,\quad
 \varphi_* =2\pi/(nu),\quad
 X_d^+= 2\pi^2\xi (1-u/(2d))/u^2.
\]
Round these four upper bounds upwards as in the table below. Put $\beta=1/B$, rounded down; every hull vertex has radius at least $h_*$, also rounded down.
{\small\setlength{\tabcolsep}{2.8pt}
\begin{center}\begin{tabular}{lcrrrrrr}
$d$&range&$\beta$&$h_*$&$E_*$&$B$&$X_d^+$&$\varphi_*$\\\hline
2&I&.83088&.97798&.80251&1.20354&14.885&.013994\\
2&II&.79758&.98233&.7525&1.25378&14.885&.013994\\
3&I&.78531&.98514&.87144&1.27337&16.615&.014025\\
3&II&.74599&.9871&.8381&1.3405&16.615&.014025\\
4&I&.76344&.98874&.90703&1.30986&17.525&.014057\\
4&II&.72172&.98985&.88203&1.38556&17.525&.014057\\
\end{tabular}\end{center}}
Indeed $\sec^2(2\pi/450)<\xi$, so (R3) applies and $dN^2(1-r^2)\le X_d^+$. In (A4)--(A5) the maximum of the combined affine cost occurs at $q$, by the same monotonicity on either side of $1/2$. Thus $BJ_F+E_*(J_E-J_F)$ bounds the angular and radial cost provided $1\ge J_E\ge J_F\ge\beta h_*$; these bounds on the coefficients follow from the construction below. Here $J_F/u\ge\eta J_E/d$. The same comparison gives $N(1-R\cos\varphi)/(R\sin\varphi)\le\pi B$, so the hull-square-gap bound (A8) gives the required hull fraction. For the radius use $E_0=\eta(1/u-t_{\min}/d)$. The first $d$ powers of $z$ applied to a site of radius one have gaps at most $2\pi/d+\varphi$. For any ray choose its preceding one of these powers, and advance by the nearest nonnegative number of $z^d$ rotations. At most $2\pi/(d\varphi)+3/2$ more rotations suffice, coming within $\varphi/2$. The initial site already has radius at least $R$. Integrating the logarithmic slope gives $-\log R\le\pi E_0\varphi/N$. The hull therefore contains the disk of radius
\[
 \cos(\varphi_*/2)\exp[-\pi E_0(2\pi/d+5\varphi_*/2)/n].
\]
The bound holds with the displayed $h_*$; for instance, use $\cos s\ge1-s^2/2$ and $e^{-x}\ge1-x+x^2/2-x^3/6$.

At $5\le d\le17$ use the entire Farey interval, taking $h_*=1-22/(d_0n)$ with $d_0$ the smallest denominator in the row. In all tables $X_d^+$ bounds $dN^2(1-r^2)$. On a group take the largest $d$ for $E_*,B,X_d^+,\varphi_*$. The last column is a chosen radius below the hull.
{\small\setlength{\tabcolsep}{3.2pt}\begin{center}\begin{tabular}{lrrrrrr}
$d$&$E_*$&$B$&$\beta$&$X_d^+$&$\varphi_*$&$U$\\\hline
5&.90929&1.41395&.70723&18.107&.014088&.9769\\
6&.92823&1.43403&.69733&18.526&.01412&.979\\
7&.94242&1.44935&.68996&18.852&.014152&.9809\\
8&.95364&1.46172&.68412&19.119&.014184&.9835\\
9&.96289&1.47211&.67929&19.348&.014216&.9841\\
10&.97075&1.48114&.67515&19.55&.014248&.9848\\
11--13&.98929&1.50317&.66526&20.06&.014346&.985\\
14--17&1.00781&1.52643&.65512&20.624&.014478&.9887\\
\end{tabular}\end{center}}

\paragraph{The gaps of the hull.}
Set $k=(1-R\cos\varphi)/(R\sin\varphi)$. Between
successive positively ordered hull vertices $r_i e^{i\theta_i}$ and
$r_{i+1}e^{i\theta_{i+1}}$, put $g_i=\theta_{i+1}-\theta_i\in(0,\pi)$.
The return image of the first vertex lies to the left of the edge, so
$r_i/r_{i+1}\le\cos g_i+k\sin g_i$. The logarithm of the right side
has initial value zero, initial derivative $k$ and second derivative at
most $-1$. Its argument is positive on the whole interval from $0$ to
$g_i$, since the endpoint values are positive and $g_i<\pi$.
Summing over the hull gives $\sum g_i^2\le4\pi k$. Extra angular marks
only subdivide gaps and hence reduce the sum of their squares.
Consequently, if we select $m$ labels including all the hull vertices and take the successive angular gaps within this set,
the normalized gap sizes $x_i=mg_i/(2\pi)$ have mean one and
\[
 m^{-1}\sum_i (x_i-1)^2\le mB/N-1.                           \tag{A8}
\]

\paragraph{Choosing the angular slopes.}
Suppose first that $m$ marks are equally spaced and that we want a
positive angular derivative at every mark. Put $a=\pi/m$. On the gap
of length $2a$ following a mark, use $H(x)=ax-x^2/2$, measuring $x$
from that mark. The values at the two ends are both zero, and the
slope falls from $a$ to $-a$. Repeating this arc gives a continuous
periodic height whose right derivative at every mark is $a$.
Here $H''\ge-1$ means that $H(	heta)+	heta^2/2$ is convex when the angle is unwrapped onto the real line. The slope jumps upwards at the joins. This is exactly the sign
allowed by $H''\ge-1$. The construction below makes room for unequal
gaps and for marks at which slope zero is enough.

Among $t$ marked angles around the circle, counting multiplicities,
designate $m$ as special. If $0<q\le1$ and
$tq^2\le m(2-q)^2$, put $a=\pi q/m$ and $C=4.935>\pi^2/2$.
There is a continuous periodic function $H$ with
\[
 0\le H\le Cq^2/m,\qquad H''\ge-1,\qquad
 H'_+\ge0\text{ at every mark},\qquad H'_+\ge a\text{ at special marks}.
                                                               \tag{A6}
\]
For its construction start at a special angle with slope $a$ and
decrease the slope at unit rate between marks. At a special mark reset
it to $a$; at an ordinary mark reset it to the maximum of its incoming
value and zero. Coincident marks can be merged, retaining a special
designation. The resulting periodic function $D$ has only upward
jumps. Its positive intervals occupy total length at most $ma$ and
have total integral at most $ma^2/2$. The negative portions are at most
$t$ triangles with total base at least $2\pi-ma$. Thus
$\int_0^{2\pi}D\le ma^2/2-(2\pi-ma)^2/(2t)\le0$. Take
$H(\theta)=\max_{u\in[\theta-2\pi,\theta]}\int_u^\theta D$.
The nonpositive period integral allows the same maximum over all
$u\le\theta$. On each gap $H$ follows the integral of $D$, staying
at zero once it returns there. This proves (A6), including its bound
by the positive integral. Between jumps the slope falls at rate at most
one, and all its jumps are upwards. Thus $H(\theta)+\theta^2/2$ is convex. In particular for $t_0\ge0$,
\[
 H(\theta+t_0)-H(\theta)\ge H'_+(\theta)t_0-t_0^2/2.
                                                               \tag{A7}
\]
Put $p(re^{i\theta})=r(1+H(\theta))$. This extends the angular function to the plane by making $p$ proportional to radius on every ray. We claim that $p$ is convex. Indeed its
tangential Hessian is $(1+H+H'')/r$ and its other Hessian entries
vanish in the orthonormal polar frame. Periodic smoothing proves the
same assertion distributionally, then uniform convergence gives
convexity off zero. On lines through zero the outgoing slope dominates
the incoming one since $1+H>0$.

\subsection{Weights on short arcs}
Assign a weight $w_i\in[0,1]$ to each mark. Suppose every arc of
length $0<L\le2\pi$, taken half-open, contains weight at most $M'$ and at most $T'$ marks. If
\begin{equation}
 a(M'+\sqrt{M'T'})\le L,\qquad a>0,                    \tag{A11}
\end{equation}
there is an angular function as in (A6)--(A7), now satisfying
$H'_+(\theta_i)\ge aw_i$ and $0\le H\le M'a^2/2$. To prove this form the periodic function $D$ of slope $-1$ which resets at each mark to the greater of the arriving slope and $aw_i$. Equivalently it is the maximum of the forward lines of slope $-1$ from all the periodic copies of the targets. At simultaneous marks only the largest target is needed. The total length of positive slopes is at most $a\sum_i w_i<2\pi$: integrating the arc weight over its initial location gives $L\sum_i w_i\le2\pi M'$ and $aM'<L$. In particular there is an open interval of negative slopes.

At a local minimum of a primitive of $D$, say $\alpha$, an upward reset has occurred, so its incoming slope is nonpositive. In the interval from $\alpha$ to $\alpha+L$ all positive slopes are supplied by marks inside that interval. Their total length is at most $aM'$, and their integrated area is at most $a^2\sum w_i^2/2\le a^2M'/2$. The negative portions have total length at least $L-aM'$ and are partitioned into at most $T'$ pieces, one after each marked point including the first at $\alpha$. Each such piece of length $b$ has negative integral of magnitude at least $b^2/2$. Thus the primitive returns to or below its initial value by $\alpha+L$.

The mean of $D$ cannot be positive. Otherwise start a primitive on a forward half-line at a point with strictly negative slope. It tends to infinity, so its last global minimum is attained in the interior, contrary to the preceding return property. If the mean is negative set $H(\theta)$ equal to the primitive at $\theta$ minus its preceding running minimum, and if the mean is zero subtract the global minimum of the periodic primitive. This is periodic. Every positive excursion starts at a local minimum and lasts at most $L$, so its height obeys the required bound. When the running minimum is constant the derivative follows $D$; on a new declining minimum it follows zero. The latter reflection occurs with an upward jump of the derivative. It follows that $H''\ge-1$ and that each marked point retains the required right derivative.

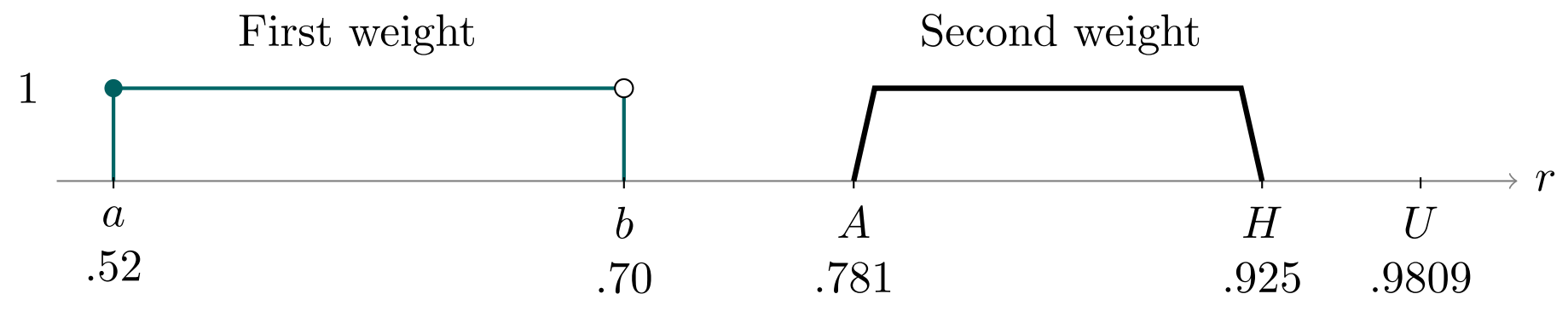
\begin{figure}[htbp]
\centering
\begin{tikzpicture}[font=\small,x=22cm,y=.72cm]
\draw[->,gray] (0,0)--(.515,0) node[right,black] {$r$};
\draw[thick,teal!80!black] (.02,0)--(.02,1)--(.20,1)--(.20,0);
\draw[very thick] (.281,0)--(.2884,1)--(.4176,1)--(.425,0);
\draw[fill=white] (.20,1) circle[radius=2pt]; \fill[teal!75!black] (.02,1) circle[radius=2pt];
\foreach \x/\name/\val in {.02/a/.52,.20/b/.70,.281/A/.781,.425/H/.925,.4809/U/.9809}{
 \draw (\x,.04)--(\x,-.08);
 \node[below,align=center] at (\x,-.09) {$\name$\\$\val$};}
\node at (.106,1.58) {First weight};
\node at (.354,1.58) {Second weight};
\node[left] at (.00,1) {$1$};
\end{tikzpicture}
\caption{Radial weights in case $7$R. A mark at $b$ is excluded
from the first count. For the second weight the shoulders both have
width $.0074$. The population variable $p$ counts \emph{all} marks
of radius below $U$, including sites whose weight here is zero.}
\label{fig:weights}
\end{figure}

\paragraph{Choice of the interior radii.}
Choose five interior radii $a<b<A<H<U$ below the hull radius
$h_*$. On $[a,b)$ give each source weight one for the first angular function, marking all sites below $A$. For the second angular function mark the sites in $[b,U)$ and use the trapezoidal weight on $[A,H]$ with shoulders of width $\delta$. Write $p=N^{-1}\#\{i:r_i<U\}$. The two pairs of positive slopes $s<v<S<V<1$ govern these angular functions and the radial functions between them.
Figure~\ref{fig:weights} shows the choice for $7$R. The first weight
records a count. The second increases gradually at the ends of its
support. Moving mass a short radial distance across either end then
changes its weight by a bounded multiple of that distance, which is
what the angular transport estimate (C1) will use.

In range I use the single row $d$I. In range II use $d$M if at least $.15N$ radii are below $.970$; since $U\ge.970$ on these rows, $p=N^{-1}\#\{i:r_i<U\}\ge.15$. Otherwise use $U=.965$, so that $p\le.15$ and at least $.85N$ radii are at least $.970$. Within this case use $d$R for $.105\le p\le.15$, and $d$S for $0\le p\le.105$, except that at $d=3$ the part $p\le.09$ uses $3$T. On $d$S and $3$T, at least $(1-p_{\max})N$ labels have radii at least $.965$. This provides classes $2b,34b,3c$ where needed. Thus some label has radius strictly below $a$ in every row.
{\small\setlength{\tabcolsep}{2.8pt}
\begin{center}\begin{tabular}{lrrrrrrrrrrr}
row&$a$&$b$&$A$&$H$&$U$&$s$&$v$&$S$&$V$&$\varepsilon$&$\delta$\\\hline
2I&.53&.605&.699&.89&.9662&.268&.55&.686&.955&.229&.013\\
2M&.53&.595&.696&.895&.9724&.268&.55&.682&.965&.199&.013\\
2R&.51&.59&.694&.89&.965&.268&.55&.677&.966&.199&.0129\\
2S&.5&.585&.685&.87&.965&.266&.526&.684&.961&.244&.0122\\
3I&.53&.623&.699&.895&.9695&.276&.582&.727&.955&.15&.0131\\
3M&.53&.639&.725&.924&.9759&.268&.561&.669&.943&.141&.0132\\
3R&.5&.592&.68&.894&.965&.264&.55&.687&.961&.149&.0137\\
3T&.48&.592&.68&.89&.965&.265&.551&.688&.951&.149&.0133\\
3S&.49&.592&.68&.89&.965&.265&.551&.688&.951&.149&.0133\\
4I&.54&.655&.736&.923&.9758&.31&.658&.745&.944&.133&.014\\
4M&.54&.668&.755&.925&.9728&.312&.664&.743&.951&.124&.0142\\
4R&.5&.622&.706&.908&.965&.299&.633&.734&.964&.135&.0143\\
4S&.49&.716&.745&.937&.965&.418&.778&.783&.926&.235&.0172\\
\end{tabular}\end{center}}
For an interval $[x,y]$ let $(K_{xy},e_{xy})$ be the (L5) packet cost per rank and the terminal excess. Its guaranteed fraction is
\[
 c(x,y)=\left[(y^2-x^2-(K_{xy}+e_{xy})/n)/K_{xy}\right]_+ .
\]
For (L5)--(L8) set $\xi=1.00020$, $\eta_1=1.00005$, and use $u$ and the unrounded (R3) bound $X_0=2(9.86961)\xi(1-u/(2d))/(du^2)$. The factors apply since $r>.99998$, $h\ge325$ and $m\ge.48$. The block constants $3.96,5.78$ in (R8) remain valid since $s\le3$ and $e>N-4$. The starting fraction $\ell$ and its high radius $H'$ are $(\beta,h_*)$ on I and M, $(.85,.970)$ on R, $(.895,.965)$ on S, $(.91,.965)$ on T. For $3$S use $.895$ also when $p>.09$. The seven pairs are rounded up to three decimals. Each is obtained from the same maximization over at most four slots in (L8); its switch set and the intervals for the possible starting rank are unchanged.
{\small\setlength{\tabcolsep}{2.8pt}
\begin{center}\begin{tabular}{lrrrr}
row&$[a,b]$&$[b,A]$&$[A,H]$&$[H,U]$\\\hline
2I&\(8.397,4.008\)&\(8.207,2.934\)&\(9.315,1.296\)&\(10.018,.34\)\\
2M&\(8.385,4.054\)&\(8.259,3.223\)&\(9.384,1.234\)&\(10.082,.276\)\\
2R&\(8.532,4.25\)&\(8.259,3.039\)&\(9.311,1.312\)&\(10.005,.353\)\\
2S&\(8.512,4.145\)&\(8.229,3.125\)&\(9.084,1.347\)&\(10.006,.352\)\\
3I&\(6.383,6.176\)&\(6.214,4.858\)&\(7.721,2.764\)&\(8.318,2.949\)\\
3M&\(6.449,6.129\)&\(6.286,4.666\)&\(7.964,2.847\)&\(8.37,2.965\)\\
3R&\(6.427,6.432\)&\(6.234,5.141\)&\(7.715,2.8\)&\(8.28,2.935\)\\
3T&\(6.48,6.503\)&\(6.217,4.919\)&\(7.676,2.839\)&\(8.28,2.935\)\\
3S&\(6.44,6.458\)&\(6.22,4.918\)&\(7.677,2.838\)&\(8.28,2.935\)\\
4I&\(5.168,7.22\)&\(5.097,5.756\)&\(6.06,4.359\)&\(6.375,4.199\)\\
4M&\(5.23,7.433\)&\(5.152,5.636\)&\(6.07,4.348\)&\(6.357,4.189\)\\
4R&\(5.236,7.833\)&\(5.042,6.039\)&\(5.967,4.479\)&\(6.308,4.159\)\\
4S&\(5.568,7.343\)&\(5.064,5.366\)&\(6.148,4.27\)&\(6.305,4.153\)\\
\end{tabular}\end{center}}
{\small\setlength{\tabcolsep}{2.8pt}
\begin{center}\begin{tabular}{lrrr}
row&$[a,A]$&$[A,U]$&$[A,H-\delta]$\\\hline
2I&\(8.741,3.664\)&\(10.113,.499\)&\(9.179,1.432\)\\
2M&\(8.758,3.68\)&\(10.196,.422\)&\(9.248,1.37\)\\
2R&\(8.912,3.87\)&\(10.095,.528\)&\(9.176,1.447\)\\
2S&\(8.914,3.744\)&\(10.076,.355\)&\(8.956,1.474\)\\
3I&\(6.652,5.908\)&\(8.364,2.994\)&\(7.608,2.852\)\\
3M&\(6.755,5.823\)&\(8.412,3.007\)&\(7.851,2.806\)\\
3R&\(6.749,6.111\)&\(8.328,2.983\)&\(7.597,2.919\)\\
3T&\(6.822,6.161\)&\(8.323,2.979\)&\(7.562,2.954\)\\
3S&\(6.775,6.123\)&\(8.324,2.98\)&\(7.563,2.953\)\\
4I&\(5.425,6.963\)&\(6.403,4.254\)&\(5.968,4.451\)\\
4M&\(5.499,7.164\)&\(6.383,4.242\)&\(5.976,4.441\)\\
4R&\(5.514,7.555\)&\(6.341,4.226\)&\(5.873,4.573\)\\
4S&\(5.67,7.242\)&\(6.331,4.206\)&\(6.035,4.383\)\\
\end{tabular}\end{center}}
For an example of the numerical reduction take row $3$M, interval
$[a,b]=[.53,.639]$. Its first rank window starts at $A_0=.74599$
and ends at $.76$, so its last cut has rank fraction at most
$U_0=.76+2/450$. At this window (L7) gives $g=0$. Hence (L8), written
in terms of the preservation probability $q$, is the lesser of a constant
cap $k$ and $m+(y-m)/(1-q)$ (the third term exceeds the cap).
Here $6.36598<k<6.36599$, $.53034<m<.53035$ and
$2.70543<y<2.70544$. The switches are $0$, $1$ and
$1-(y-m)/(k-m)$. The branch change lies near $.6273$, found by solving
$m+(y-m)/(1-q)=k$. The budget is $1+e_*<1.02515$. Of the three
slots at most one needs an interior
value. Placing the others at switches gives $S_3<13.21412$ in this
window. Repeating the calculation over the rank windows gives
$\mu<6.449$ and $E_p<6.129$. These are the two entries in column
$[a,b]$ of row $3$M. The calculation illustrates the need to take
the maximum over windows before using a coefficient in (L5).

In row $4$S the two disjoint counts $c(a,A)+c(A,U)>.10620$ already contradict $p\le.105$. We henceforth omit this row.

For $5\le d\le8$ take row M when $p\ge.21$ and row R when $p\le.21$, using the same $U$ from the Farey table for both. In row R at least $.79N$ marks have radius at least $U\ge.976$, providing the indicated prefix class. For $9\le d\le17$ there is one case per row. Use the first radius $a=.55$ at $5$M, $.57$ at $6$M or $7$M, $.59$ at $8$M and all higher denominators, and $.52$ on every R row from $5$ to $8$. The other radii and slopes are as follows.
{\small\setlength{\tabcolsep}{3.2pt}\begin{center}\begin{tabular}{lrrrrrrrrr}
$d$&$b$&$A$&$H$&$s$&$v$&$S$&$V$&$\varepsilon$&$\delta$\\\hline
5&.684&.771&.935&.33&.713&.769&.936&.121&.0141\\
6&.678&.758&.925&.283&.611&.717&.942&.129&.0101\\
7&.7&.781&.925&.254&.56&.679&.923&.13&.0074\\
8&.684&.772&.922&.256&.551&.677&.936&.162&.0082\\
9&.708&.789&.935&.251&.554&.67&.934&.138&.0071\\
10&.685&.772&.922&.256&.551&.682&.944&.138&.0077\\
11--13&.714&.797&.937&.25&.553&.665&.93&.134&.0064\\
14--17&.73&.8&.94&.24&.555&.686&.944&.112&.0063\\
\end{tabular}\end{center}}

To calculate the packet cost on these rows use $\ell=b_d$, $H'=1-22/(dn)$, or $(\ell,H')=(.79,U)$ on the R rows, for which every relevant cut contains this many high sites. Apply (L5)--(L8) separately at each denominator. For a grouped row use the maximum of each entry over its denominators. The first entry bounds cost per cut and the second bounds the excess at the last packet. The pairs are
{\small\setlength{\tabcolsep}{3.2pt}\begin{center}\begin{tabular}{lrrrr}
row&$[a,A]$&$[a,b]$&$[b,A]$&$[A,H]$\\\hline
5M&(4.623,7.861)&(4.374,8.11)&(4.363,6.363)&(5.028,5.612)\\
5R&(4.734,8.045)&(4.482,8.297)&(4.36,6.365)&(5.027,5.39)\\
6M&(3.933,10.397)&(3.715,10.474)&(3.797,8.857)&(4.272,8.272)\\
6R&(4.066,11.545)&(3.854,11.617)&(3.797,8.656)&(4.272,8.272)\\
7M&(3.488,11.41)&(3.289,11.483)&(3.36,9.441)&(3.72,8.825)\\
7R&(3.603,12.76)&(3.409,12.827)&(3.36,9.227)&(3.72,8.825)\\
8M&(3.048,11.997)&(2.855,12.069)&(2.969,10.315)&(3.295,9.482)\\
8R&(3.18,14.054)&(2.991,14.115)&(2.969,10.315)&(3.295,9.482)\\
9&(2.778,13.34)&(2.617,13.4)&(2.701,10.847)&(2.994,10.257)\\
10&(2.498,14.475)&(2.341,14.527)&(2.434,12.169)&(2.702,11.035)\\
11--13&(2.334,18.972)&(2.195,18.977)&(2.269,14.72)&(2.505,13.303)\\
14--17&(1.884,27.189)&(1.796,27.159)&(1.831,19.211)&(2.02,17.084)\\
\end{tabular}\end{center}}
{\small\setlength{\tabcolsep}{3.2pt}\begin{center}\begin{tabular}{lrrr}
row&$[A,U]$&$[H,U]$&$[A,H-\delta]$\\\hline
5M&(5.234,5.406)&(5.213,5.205)&(4.952,5.687)\\
5R&(5.233,5.184)&(5.213,5.205)&(4.951,5.466)\\
6M&(4.435,8.268)&(4.435,7.824)&(4.241,8.274)\\
6R&(4.435,8.268)&(4.435,7.824)&(4.241,8.274)\\
7M&(3.864,8.806)&(3.864,8.458)&(3.701,8.828)\\
7R&(3.864,8.806)&(3.864,8.458)&(3.701,8.828)\\
8M&(3.435,9.455)&(3.435,8.971)&(3.277,9.486)\\
8R&(3.435,9.455)&(3.435,8.971)&(3.277,9.486)\\
9&(3.094,10.23)&(3.094,9.444)&(2.979,10.262)\\
10&(2.819,11.002)&(2.819,10.234)&(2.688,11.039)\\
11--13&(2.587,13.284)&(2.587,12.128)&(2.494,13.305)\\
14--17&(2.086,17.08)&(2.086,15.053)&(2.012,17.085)\\
\end{tabular}\end{center}}
\subsubsection*{Column-weighted angular transport}
The angular error can be charged to the radius of the receiving mark. If $Pv=R' e^{i\theta}v$, choose $\theta_i=\arg v_i$ at each nonzero site and any angle at a zero site. Put
\[
 e_i=(Pr)_i-R'r_i=\sum_jP_{ij}r_j(1-\cos(\theta_j-\theta_i-\theta)),
 \qquad \sum_i e_i=(1-R')\sum_i r_i .
\]
Take $G(x)=\max\{c_-(L-x),0,c_+(x-H')\}$ and $w=(1-G)_+$. Suppose $w(r)\le C_0r$ and put $S_0=\sum r_i$, $S_H=\sum_{r_i>H'} r_i$. If an arc $I$ is at angular distance at least $2\epsilon$ from its rotation by $\theta$, let $I_+,I_-$ be the two rotated arcs enlarged by $\epsilon$, where $0<\epsilon\le\pi/2$. Their weights satisfy
\[
 2W(I)\le W(I_+)+W(I_-)
 +(1-R')\{2c_+S_H+S_0[c_-+\max(c_-,C_0/(1-\cos\epsilon))]\}.\tag{C1}
\]
In fact for any set $J$ of positive-weight rows,
$\sum_J(w-Pw)\le\sum_J(PG-G)\le\sum_{i\notin J}(G-PG)_+$.
On a row $r_i<L$ scalar Jensen bounds the last summand by $c_-e_i$. On $r_i>H'$ it is at most $c_+(1-R')r_i$. Elsewhere it is zero. For the forward comparison take $J$ within $I$, for the reverse one outside $I_-$. Each missed edge carries destination weight at most $C_0r_j$. If missed twice, both endpoints belong to $I$ and its angular error is at least $2\epsilon$. Since $1-\cos(2\epsilon)\ge2(1-\cos\epsilon)$, the row pays both misses from a single $e_i$. A row in neither selected set pays at most $2c_-e_i$ at low radius, one in exactly one set pays at most $(c_-+C_0/(1-\cos\epsilon))e_i$, and a row in both has no low-radius charge. This proves (C1). Approximate rotations may be used when each omitted angular difference is at least $\epsilon$.

Use $k$ root directions, as in the following table, pairing opposite exponents. The comparisons at exponents $j$ and $-j$ use the forward and reverse estimates of (C1) for the same power $A^{|j|}$. Their root rotations are separated by $\pi\gamma_0$ or more. Use every root for $d\le5$, the powers $0,\pm2$ for $d=6$, and at larger denominators compare the identity with the closest root to $-1$, averaging the two signs of its exponent. In the last comparison only separation from the identity matters. A centred exponent has size at most $d/2$. Its angular error is at most $\varphi_*/2$ and its mean absolute size in the comparison is at most $d\sigma$.
\[
\begin{array}{l|ccccc}
d&2\text{--}5&6&7&8,10&9\\\hline
k&d&3&2&2&2\\
\gamma_0&2/d&2/3&6/7&1&8/9\\
\sigma&\lfloor d^2/4\rfloor/d^2&2/9&3/14&1/4&2/9
\end{array}\qquad
\begin{array}{l|cc}
d&11\text{--}13&14\text{--}17\\\hline
k&2&2\\
\gamma_0&10/11&14/15\\
\sigma&1/4&1/4
\end{array}
\]
Indeed an odd denominator has an antipodal error of $\pi/d$ and an even denominator has none. In a two-arc comparison both enlarged rotated arcs avoid the original arc; their averaged weight is at most the mass outside it. With $k\ge3$ the enlarged root copies have disjoint interiors. Use half-open arcs $[\alpha,\alpha+L)$ on the circle, and enlarge a rotated arc by $\epsilon$ at each end with the same convention. Touching arcs count a mark at their common endpoint once. A missed endpoint has angular error at least $\epsilon$, which suffices in (C1). In (A11) the reset at the final endpoint has no effect on the integral. For the arc of length $L=\pi\gamma_0-2\varepsilon$, apply (C1) with $\epsilon=\varepsilon-\varphi_*/2$ and put $e_*=\sigma X_d^+/n$, $c_\epsilon=\epsilon^2/2-\epsilon^4/24$. A doubly charged edge has error at least $2\varepsilon-\varphi_*/2\ge2\epsilon$. The following upper bounds for the normalized radial masses will be used:
\begin{align*}
 \overline S&=1-c(a,b)(1-b)-c(b,A)(1-A)-c(A,H)(1-H)-c(H,U)(1-U),\\
 \overline S_H&=1-c(a,A)-c(A,H-\delta),\qquad
 \overline S_U=\min(\overline S,1-c(a,A)-c(A,U)).
\end{align*}
These bound $S_0/N,S_H/N$ at thresholds $H-\delta,U$ as indicated. For the second gauge's trapezoid use $L=A+\delta$, $H'=H-\delta$, $c_-=c_+=1/\delta$, $C_0=1/(A+\delta)$. With the lower core weight $M_2$ defined below, the fraction of a core's total weight on one tested arc is at most
\[
 \kappa=\frac1k+\frac{e_*}{M_2}\left[
 .501\frac{\overline S_H}{\delta}+
 .251\overline S\left(\delta^{-1}+\max(\delta^{-1},[(A+\delta)c_\epsilon]^{-1})\right)\right].
\]
For the unweighted count on $[b,U]$ use $w(x)=\min\{x/b,1,(L'-x)_+/(L'-U)\}$. All those sites have weight one; the total weight is at most the fraction below $L'$ times $N$. In the count of marks entering (A15) use the following upper bound, in units of $n/N$:
\[
 T(L',b')=n(1-b')/k+e_*\left[
 .501\overline S_U/(L'-U)+.251\overline S\{1/b+1/(b c_\epsilon)\}\right].
\]
Take $T=T(h_*,\beta)$, also taking the minimum with $T(.970,.85)$ on R, S and T rows with $d\le4$. The constants follow from $1-r^k\le k(1-r^2)/(1+r)$ and mean absolute exponent $d\sigma$, assigning half of the paired cost to each sign. This weighted estimate allows $L'-U$ anywhere above zero.

\paragraph{The two interior angular targets.}
For a count $c(x,y)$ use (L5) with the pair for exactly that interval. For a subinterval within a moment integral, keep the pair for the original interval. Put
\[
 \underline h_1=c(a,b),\quad J=c(b,A),\quad o=c(A,U),\quad
 P_1=\underline h_1,\quad P_2=J+o,\quad V_2=J+c(H,U),
\]
\[
 \underline h_2=
 [H^2-A^2-\delta(H+A)-(K_{AH}+e_{AH})/n]/K_{AH},\qquad M_j=n\underline h_j .
\]
For $a'\le x\le y$ and a cost pair $(k,e)$ put
\[
 m(a';x,y;k,e)=\{(x-a')(y^2-x^2)+(y-x)^2(x+2(y-x)/3)
                  -(k+e)(y-a')/n\}/k .
\]
Take $D_1=m(a;a,b;K_{ab},e_{ab})$ and take $D_2$ to be the average of
$m(A;A+t\delta,H-t\delta;K_{AH},e_{AH})$ on $0\le t\le1$. It is the Simpson average of the values at $0,1/2,1$. These are lower bounds for $A_j=N^{-1}\sum_i w_i(r_i-a_j)$ with $(a_1,a_2)=(a,A)$. Indeed $(x-a_j)c(x,y)$ plus the integral of $c(t,y)$ from $x$ to $y$ bounds the excess moment on $[x,y)$. Use the superlevel intervals of the trapezoid and integrate their counts to get $\underline h_2$ and their moments to get $D_2$. Thus the actual weight fractions $h_j$ satisfy $h_j\ge\underline h_j$.

Intersect the stated interval for $p$ with $[c(a,A)+o,1-\beta]$. In rows $2$S and $4$R divide the intersection into two equal parts. Each resulting interval will be written $[p_-,p_+]$.
Put
\[
 \lambda_1=\frac{(A-b)(v-s)}{CAs},\quad
 \lambda_2=\min\left(\frac{(A-b)(S-v)}{CbS},\frac{(U-H)(V-S)}{CUS}\right).
\]
The maximum number of marks available to the second angular function
also gives
\[
 T_0=\min(T,n(p_+-P_1)).
\]
Set $\gamma=\gamma_0-2\varepsilon/\pi$, $c_1=E_*Cb/n$, $c_2=E_*CU\kappa/n$. The constants $q_-,q_+,q_2$ in the table are used as written. The columns for $p_\pm$ are approximate, with exact endpoints as just defined. For the displayed tangent parameter $y>0$ they satisfy
\[
 q_\pm\le \min\left(1,\sqrt{\lambda_1 M_1},
 \frac2{1+y}-\frac{p_\pm-o-P_1y^2}{P_1y(1+y)^2}\right),\qquad
 q_2\le\min\left(1,\sqrt{\lambda_2 M_2/\kappa},
 \frac{\gamma}{\kappa+\sqrt{\kappa T_0/M_2}}\right).
\]
Interpolate $q_1(p)$ linearly between $q_-$ and $q_+$. The affine expression in $p$ is the tangent lower bound to the convex function $2/(1+\sqrt{(p-o)/P_1})$. This gives feasibility throughout the interval.

There are at most $N(p-o)$ first-function marks. Apply (A6) with target slope $\pi Q_1/(Nh_1)$ on the special ones. For the second function use (A11) with target slope $\pi Q_2 w_i/(Nh_2)$. Its arc weight is at most $\kappa Nh_2$ and its arc count at most $NT_0/n$. Ignoring for a moment an additional height restriction, the targets permitted are
\[
 Q_1(h)=\min\{\sqrt{\lambda_1nh},2\sqrt h/(\sqrt{p-o}+\sqrt h)\},\qquad
 Q_2(h)=\min\{1,\sqrt{\lambda_2nh/\kappa},\gamma/(\kappa+\sqrt{\kappa T_0/(nh)})\}. \tag{A15}
\]
At $q=1$ the condition of (A6) forces $t=m$, and the same construction applies. At actual core size $h$ also impose the bound $q_j\sqrt{h/\underline h_j}$, where $q_1=q_1(p)$. The resulting angular heights are bounded by $C\min(\lambda_1,q_-^2/M_1)$ and $C\min(\lambda_2,\kappa q_2^2/M_2)$, respectively. More precisely the height at core $j$ is bounded by $C\kappa_j Q^2/(Nh)$ with $\kappa_1=1$, $\kappa_2=\kappa$ and the actual restricted target $Q$.

These choices give an angular contribution at least $q_j(a_j+D_j/\underline h_j)-c_jq_j^2$ after charging the angular height and replacing the radial moment by its lower bound. To see the monotonicity on the two increasing branches of (A15), write $A'=a_j-2c_j$ and $D'=D_j$. The charged expression is $Q(h)(A'+D'/h)+c_jQ(h)(2-Q(h))$.
The second term increases on $0\le Q\le1$. It suffices to show that the first increases. On a square-root branch this follows from $A'\ge3D'/\underline h_j$. For the second branch of $Q_1$, differentiating in $\sqrt h$ gives the sign of $A'\sqrt{p-o}-D'\sqrt{p-o}/h-2D'/\sqrt h\ge0$ since $h\le p-o$. On the last branch of $Q_2$, use $m=nh$, $v'=1/\sqrt m$. The derivative in $v'$ has the sign of $2nD'\kappa v'+nD'\sqrt{\kappa T_0}(v')^2-A'\sqrt{\kappa T_0}\le0$, using $T_0\ge\kappa M_2$ and $m\ge M_2$. The constant branch supplies at least $a_2-c_2$, sufficient since $q_2(a_2+D_2/\underline h_2)\le a_2-c_2$. The added square-root branch agrees with $q_j$ at the lower endpoint and obeys the same monotonicity. Finally at $h=\underline h_j$ the expression increases with the target up to one. All the used tests, including $3D_j/\underline h_j\le a_j-2c_j$, are direct from the displayed parameters.

\paragraph{Fitting the functions together.}
Draw a radial line $l_1(x)=s(x-a)$ and a second line $l_2$ of slope $S$ such that $l_2(b)-l_1(b)=-(A-b)(S-v)$. Give $l_H$ slope $V$ and make it meet $l_2$ at $H$. Give the outer line slope one and make it meet $l_H$ at a radius $b_o\in(U,h_*]$. Take the maximum of zero, the two outer radial lines and the two inner functions $l_j(|x|)+s_j|x|H_j(\arg x)$ built with (A6), (A11). It is convex by (A7). Below $b$ the second angular branch lies beneath $l_1$ by the definition of $\lambda_2$; at $A$ the first branch is below $l_2$ by the definition of $\lambda_1$ and stays there since $s(1+C\lambda_1)\le v<S$. At $U$ the second branch is below $l_H$ and stays there by $S(1+C\lambda_2)\le V$. Choose a supporting branch at each mark, giving each core its own branch wherever tied. Thus each angular branch only needs its chosen marks. Above $U$ the two outer radial lines supply all the source heights.
\subsubsection*{Averaging a radial breakpoint near the hull}
Let $f_0$ be the maximum of the two angular branches and the radial branches just drawn. Suppose $U<b_o\le h_*$, so that above $b_o$ the outer line is active, and between $U$ and $b_o$ the active branch is the radial line $l_H$, of slope $V$. Let $F$ be the corresponding radial maximum, using the same lines without their gauges. At a site $i$ with $r_i\ge b_o$ we have $F(r_i)=f_0(v_i)$ and $F$ has subgradient $1$. Increase the tangential slope of the supporting plane there as far as the other marks permit. Its final value is
\[
 T_i=\min_{0<\alpha_{ij}<\pi}
 \left[\tan(\alpha_{ij}/2)+
 \frac{F(r_j)-F(r_i)-(r_j-r_i)}{r_j\sin\alpha_{ij}}\right],
 \quad\alpha_{ij}=\theta_j-\theta_i\pmod {2\pi}.
\]
The minimum ranges over nonzero forward sites. Such a site exists by the imaginary part of the row equation at the return. At a nonforward site the tangential projection is nonpositive, so the plane lies below its height automatically. Add all these planes to the maximum defining $f_0$. They lie below $F(r_j)\le f_0(v_j)$ at every site. The resulting function is convex, leaves each source height fixed, and gives at each outer site the gain
$r_i[R\sin\varphi\,T_i-(1-R\cos\varphi)]$. Between $U$ and $b_o$ the radial evaluation costs at most $(1-R)Vr_i$. Below $U$ the two inner evaluations continue to hold. Thus the normalized gain uses (A5), with $J_F\ge\beta h_*$ and $J_E\le1$. At every label with $r_i\ge U$ the initial cost allowance $B$ saves at least $B(1-r_i)$; whenever $r_i<b_o$ it saves an additional $(B-E_*V)r_i$.

Let $\mu$ be the fraction of actual hull vertices and $M=1-p_0$, where $p_0$ is any guaranteed lower bound on $p$. Then $\beta\le\mu\le M$. Fix a hull interval $u_0\le\mu\le u_1$, a number $D>1$ and a width $w>0$. Take $b_o$ uniform on $[h_*-w,h_*]$. Write $T_0=\pi D/(N\mu)$ and assume
\[
 h_*-w-\frac{(\pi D/(n u_0))^2}{2(1-V)}\ \ge U.                 \tag{C2}
\]
For the hull gap angles $g_i$, the variables $N\mu g_i/(2\pi)$ have mean one and variance at most $B\mu-1$ by (A8). The sum of the clipped hull half-gaps divided by $\pi$ is at least one minus
\[
 \tfrac12\{\sqrt{(D-1)^2+B\mu-1}-(D-1)\}.
\]
This follows from $2(x-D)_+=|x-D|+x-D$ and Cauchy--Schwarz. Inserting the other sites of radius at least $b_o$ only increases the sum of these clipped gaps.

A site below $b_o$ can reduce the clipped target at at most one active angular gap, the one that contains it. If its radius is $r_j$, the radial height excess is at least $(1-V)(b_o-r_j)$. Writing $x=\tan(\alpha_{ij}/2)$ and using $\sin\alpha_{ij}\le2x$ bounds the amount lost at that gap by
\[
 \big[T_0-\sqrt{2(1-V)(b_o-r_j)/r_j}\big]_+ .
\]
Indeed the angular cost is at least $x+(1-V)(b_o-r_j)/(2r_jx)$. This argument also applies when the identity of the preceding active site varies with $b_o$. Condition (C2) ensures that a site responsible for loss has $r_j\ge U$. Integrating the positive part over $b_o$ costs at most $r_jT_0^3/[6(1-V)w]$. If $x_j=\Pr(r_j\ge b_o)$, the additional bound $T_0(1-x_j)$ is immediate.

To keep track of radial weights, replace each active coefficient $r_i$ by $\min(r_i,h_*)$. Its deficiency from $h_*$ costs at most $T_0(h_*-r_i)_+$ when active. For a nonhull site the mean of $(h_*-r_i)_+\mathbf 1_{r_i\ge b_o}$ is at most $w/4$, by maximizing $x(1-x/w)$ over $0\le x\le w$. Sites above $h_*$ incur zero loss. Put
\[
 s_\mu=D/\mu,\quad g=(B-E_*V)U,\qquad
 a_\mu=\frac{\pi^2h_*D^3}{6(1-V)w n^2\mu^3}.
\]
After division by $\pi$ and rescaling each per-label contribution by $N$, a site can steal at most $\min(a_\mu,s_\mu(1-x_j))$, while its additional saving is at least $g(1-x_j)$. The net is at most $a_\mu(1-g/s_\mu)_+$. Indeed for two nonnegative numbers $k,y$ the positive part of $\min(k,s_\mu y)-gy$ is at most $k(1-g/s_\mu)_+$. Each site below $h_*$ also has the initial saving $B(1-h_*)$. Consequently the outer angular contribution together with the savings above $U$ is at least
\begin{align*}
 O=(1-\varphi_*^2/6)\bigg[h_*&
 -\frac{h_*}{2}\left\{\sqrt{(D-1)^2+B\mu-1}-(D-1)\right\}\\
 &-(M-\mu)\Big[a_\mu(1-g\mu/D)_++D w/(4\mu)-B(1-h_*)\Big]_+\bigg].\tag{C3}
\end{align*}
The bracket is positive in every application below. The lower sine factor bounds $\sin\varphi/\varphi$, and the unscaled nonnegative savings may also be multiplied by it. Since the lower bound is the average of the estimates for $b_o$, some choice attains it. Repeated angles are allowed; take zero gap for all but the last label at a repeated angle.

For a two-endpoint evaluation of (C3), majorize $\sqrt{(D-1)^2+B\mu-1}$ by its tangent at $t\in[u_0,u_1]$. One may use the upper bound
$(Z^2+(D-1)^2+B\mu-1)/(2Z)$ with any $Z>0$, taking $Z^2=(D-1)^2+(Bt-1)_+$.
The function $\mu^{-3}(1-g\mu/D)_+$ is convex decreasing, by direct differentiation on its positive interval. The positive part in (C3) is thus convex decreasing. Multiplication by $M-\mu$ preserves convexity on $[u_0,u_1]$. The largest estimated loss occurs at an endpoint. In a row take $p_0$ to be the lower endpoint of its feasible range before the optional angular split. Divide $[\beta,M]$ into the number of equal intervals indicated, inserting an extra division at $.85$ where requested, and read $D$ in increasing interval order. The attached letters L, M, R indicate the left endpoint, midpoint, or right endpoint for the tangent. Put
\[
 w=\min\left\{h_*-U-\frac{(\pi D/(nu_0))^2}{2(1-V)},\ 
 \frac{\pi D\sqrt{2h_*/3}}{n((u_0+u_1)/2)\sqrt{1-V}}\right\}.
\]
The first entry enforces (C2). Every indicated width is positive.
The table supplies a choice of plane for each range of hull populations.
Its final column is a guaranteed contribution from the outer marks.
One may choose a different tangent on each range because only one
value of the hull population occurs in a given configuration.
{\small\setlength{\tabcolsep}{2.8pt}
\begin{center}\begin{tabular}{lcccr}
row&equal intervals&cut at $.85$&$D$ and tangent&$O>$\\\hline
2I&1&--&2R&.9474\\
2M&1&--&1.8R&.9626\\
2R&1&Y&2R/2.2M&.9577\\
2S&3&Y&1.8R/1.8R/2.2R/2.4L&.9513\\
3I&1&--&2.2L&.9468\\
3M&1&--&2R&.9537\\
3R&1&Y&2R/2.4R&.9487\\
3T&1&Y&2.2R/2.6L&.9509\\
3S&1&Y&2.2R/2.6L&.9517\\
4I&1&--&2.2L&.9477\\
4M&2&--&2R/2.4L&.9546\\
4R&1&Y&2R/2.6L&.9425\\
\end{tabular}\end{center}}

On the rows from $5$ onward the following table gives $D,w$ themselves. Two entries indicate the consecutive halves of the hull interval. In this table take $Z=D-1$ in the tangent majorant and replace $B\mu-1$ by $(B\mu-1)_+$ there. This bounds the hull loss by $h_*(B\mu-1)_+/[4(D-1)]$ and is again convex. Direct substitution at the interval ends gives the listed $O$ rounded down.
{\small\setlength{\tabcolsep}{3.2pt}\begin{center}\begin{tabular}{lrrr}
row&$D$&$w$&$O$\\\hline
5M&2.2&.00963&.9639\\
5R&2, 2.4&.01027, .0098&.9455\\
6M&2.1&.00904&.9604\\
6R&2.1&.00904&.9389\\
7M&2.2&.00889&.9618\\
7R&2, 2.3&.00945, .00922&.9536\\
8M&2&.00713&.9554\\
8R&2&.00713&.9481\\
9&1.9, 2.2&.00757, .00725&.9512\\
10&1.9&.00686&.9463\\
11--13&2&.0074&.952\\
14--17&1.8&.00452&.953\\
\end{tabular}\end{center}}

\subsection{The two endpoint inequalities}\label{app:gains}
Fix a row, so its radii, angular functions and coefficients are chosen.
The variable that remains is the fraction $p$ of labels inside $U$.
The following calculation gathers the contributions of the inner marks,
the two radial intervals, and the outer supporting planes into one
function $\Delta(p)$. Its concavity will reduce a whole interval of
populations to two endpoint comparisons.
In the first radial bin write $p_1=N^{-1}\#\{i:r_i<b\}$, in the second write $p_2=p-p_1$. The two angular branches have gains bounded from below at each active mark by
\[
     s_jr_i[R\varphi H'_{j,+}-R\varphi^2/2-(1-R)(1+H_j)] .
\]
This follows by evaluating the same branch at the rotated point and applying (A7). A radial branch costs at most $(1-R)Vr_i$ between $U$ and $b_o$. The outer tilted plane pays the cost proved in (C3). Since $s(1+C\lambda_1)\le v$, $S(1+C\lambda_2)\le V$, the actual means in (A5) satisfy $J_E\le1$, $J_F\ge\beta h_*$ as claimed. Divide by $\pi R\varphi$; the credits at the two inner cores are $s_jQ_j(a_j+A_j/h_j)$.
\begin{align*}
 W_1&=Bs+E_*sC\min(\lambda_1,q_-^2/M_1),\\
 W_2&=\max\{E_*V,BS+E_*SC\min(\lambda_2,\kappa q_2^2/M_2)\},\\
 z_1&=B-aW_1,\quad z_2=B-UW_2,\quad z_*=\min(z_1,z_2),\\
 b'_1&=B,\qquad b'_2=\min(B,q_2/(p_+-P_1-V_2)),\\
 L_1&=a(W_1-Bs),\quad L_2=U(W_2-BS)+b'_2 S(U-A).
\end{align*}
Use $b'_2=B$ if its denominator is nonpositive. The tabulated targets have $q_+\ge B(p_+-P_2)$, as required in the first core.
For $P^-(x,y)=\{(y-x)^2(x+(y-x)/3)-(K_{xy}+e_{xy})(y-x)/n\}/K_{xy}$ set
\[
 Q_*=W_2[J(U-A)+P^-(b,A)+\max(0,P^-(H,U))].
\]
In the final line use $(a_1,a_2)=(a,A)$, $(s_1,s_2)=(s,S)$, $V_1=0$. The outer planes add no height to the source marks. 
Here $W_1,W_2$ are upper bounds on the cost per unit radius from (A5) in the two bins. In the second bin the first branch costs at most $Bv r_i<BSr_i$. Initially allot $B$ to every label. In bin one the $h_1N$ special marks use their true radii and every other mark has radius below $a$. In bin two allot $UW_2$ to each mark. Active core weights retain a further nonnegative saving against these allotments. In bin $j$ this amounts to
\[
 u_j(W_j-Bs_j)h_j + Bs_j[(u_j-a_j)h_j-A_j]
       -E_*C u_js_j\kappa_j Q_j^2/n,
 \qquad (u_1,u_2)=(b,U).
\]
For the first bin subtract $(b-a)W_1h_1$ because its initial allotment already used $a$ in place of $b$. The two disjoint intervals $[b,A)$ and $[H,U)$ in bin two supply $Q_*$ by integrating the counts of $y-r_i$. Since $p_j\ge\max(P_j,h_j+V_j)$ and $Q_j/h_j\ge b'_j$ (use $Q_1\ge q_+$ and $h_1\le p_+-P_2$ for the first core, and $Q_2\ge q_2$ and $h_2\le p_+-P_1-V_2$ for the second), replace $B$ by $b'_j$ in the nonnegative square bracket. The remaining coefficient of $A_j$ is nonnegative, which allows its replacement by $D_j$. The coefficients $L_j$ are nonnegative. Formula (C3) includes the credits for all the outer labels. Consequently the normalized signed gain is bounded below by
\begin{align*}
 \Delta(p)=Q_*+O-B+z_*p+\sum_{j=1}^2\{&
 (z_j-z_*)\max(P_j,\underline h_j+V_j)+L_j\underline h_j\\
 &+s_j[q_j(a_j+D_j/\underline h_j)-c_jq_j^2-b'_jD_j]\}.
\end{align*}
On each listed interval $q_1$ is affine in $p$, and $q_2$ and the
other coefficients are fixed. The only quadratic dependence is
$-s_1c_1q_1(p)^2$, with $s_1c_1\ge0$. This proves concavity.
Write $p=(1-\alpha)p_-+\alpha p_+$, $0\le\alpha\le1$. Then
\[
 \Delta(p)\ge(1-\alpha)\Delta(p_-)+\alpha\Delta(p_+).
\]
Here we need a \emph{lower} bound. Positive values at both ends place
the whole graph above zero. The last two columns in the following
tables therefore decide every population in the interval, even as
$N$ grows. The columns headed $p_-,p_+$ display rounded values for
orientation; the calculation uses their defining expressions.
{\small\setlength{\tabcolsep}{2.8pt}
\begin{center}\begin{tabular}{lrrrrrrrr}
row&$p_-$&$p_+$&$y$&$q_-$&$q_+$&$q_2$&$\Delta_-$&$\Delta_+$\\\hline
2I&.06227&.16912&2.9029&.29739&.28251&.37719&.028&.0674\\
2M&.15&.20242&5.7503&.27481&.27481&.39852&.0366&.0579\\
2R&.105&.15&2.4345&.31691&.25127&.39664&.0171&.0266\\
2Sa&.06499&.08499&1.2396&.31301&.255&.41903&.0025&.0026\\
2Sb&.08499&.105&1.5856&.31301&.24088&.41903&.0106&.0087\\
3I&.07797&.21469&2.0871&.36958&.20438&.44824&.0353&.0592\\
3M&.15&.25401&4.1116&.42707&.42412&.42968&.0015&.0482\\
3R&.105&.15&2.8421&.37814&.37814&.51917&.0048&.0244\\
3T&.08309&.09&2.315&.42339&.42339&.51595&.003&.006\\
3S&.09&.105&2.8121&.40195&.40195&.51589&.0032&.0097\\
4I&.10144&.23656&1.8605&.48944&.38241&.44034&.0588&.0866\\
4M&.15&.27828&2.9561&.53543&.4951&.40296&.0027&.0435\\
4Ra&.105&.1275&.9506&.49961&.43012&.52051&.0152&.0127\\
4Rb&.1275&.15&1.1388&.49961&.38763&.52051&.024&.0147\\
\end{tabular}\end{center}}
A displayed packet cost is a rational upper bound on its exact value. For the remaining square roots a lower or upper bound at the twelfth decimal suffices by squaring, using $3.141592<\pi<3.141593$; the indicated strict margins are rounded downward.

{\small\setlength{\tabcolsep}{3.2pt}\begin{center}\begin{tabular}{lrrrrrrrr}
row&$p_-$&$p_+$&$y$&$q_-$&$q_+$&$q_2$&$\Delta_-$&$\Delta_+$\\\hline
5M&.21&.29277&3.2112&.57464&.52844&.39755&.0334&.0528\\
5R&.1268&.21&2.1868&.67115&.66904&.39816&.0078&.0362\\
6M&.21&.30267&1.7706&.55678&.36357&.48278&.034&.0395\\
6R&.14647&.21&1.0984&.66956&.5348&.48297&.0108&.0152\\
7M&.21&.31004&1.744&.6769&.5322&.43112&.009&.0336\\
7R&.16802&.21&1.1646&.78316&.74312&.43145&.0036&.0174\\
8M&.21&.31588&1.9252&.60272&.47613&.50894&.0181&.0484\\
8R&.19009&.21&1.1384&.79906&.78693&.50923&.0434&.0511\\
9&.18813&.32071&1.9386&.71252&.6146&.48286&.0102&.0597\\
10&.20586&.32485&1.4302&.65044&.39205&.56757&.0383&.0539\\
11--13&.21859&.33474&1.5557&.77533&.60834&.50521&.0264&.0593\\
14--17&.26203&.34488&1.4612&.83806&.71104&.53845&.0681&.089\\
\end{tabular}\end{center}}

For example consider row $7$R. At least $.79N$ sites are at $U=.9809$ or above, so $a=.52$. The core weights divided by $N$ are bounded below by $.05383396$ and $.05515102$. The corresponding lower moments exceed $.00417764$ and $.00357356$. The exact counts give $p_-=.16801517\ldots$, $p_+=.21$. From the weighted arc estimate $\kappa<.5478670$, $T<70.48609$, and $T_0<70.27472$ after imposing the total population. Bisect $[.68996,1-p_-]$, the hull interval. Substitution of the two tabulated pairs $(D,w)$ gives $O>.95367$ before rounding down to $.9536$. The targets $.78316,.74312,.43145$ then give $\Delta(p_-)>.00367$ and $\Delta(p_+)>.01749$.

The input parameters in row $7$R are exact decimals. For instance its
two packet pairs at $[a,A]$ and $[A,U]$ are $(3.603,12.76)$ and
$(3.864,8.806)$. Thus the lower endpoint can already be obtained on
paper from two fractions:
\[
 p_-=
 \frac{.781^2-.52^2-(3.603+12.76)/450}{3.603}
 +\frac{.9809^2-.781^2-(3.864+8.806)/450}{3.864}.
\]
The contributions are $.08415175\ldots$ and $.08386341\ldots$.
They add because the intervals $[a,A)$ and $[A,U)$ are disjoint.
Thus more than $.16801517$ of the marks lie below $U$. The definition
of case $7$R gives $p_+=.21$. Substitution in the
two endpoint expressions gives the figures above. This also explains
why the rounded entry $.16802$ in the table cannot serve as the
definition of the interval.

Every endpoint gain is strictly positive. Concavity now gives $\sum_i f(z^dv_i)>\sum_i f(v_i)$ for some choice of the outer breakpoint and its supporting planes. This contradicts (A1) and excludes all the remaining denominators. The regular shift construction in the introduction proves the reverse inclusion, completing the proof of Theorem~\ref{main}.

%% file: smallorders.tex
\section{The eigenvalue regions at orders five and six}
\label{app:loworders}

The number of marks now equals the matrix order. This gives another
use for the polygons $K_c(v)$ of Section~\ref{sec:object}. One can
group continuously varying configurations according to which marks
are vertices of their convex hull (the \emph{extreme} marks) and where
the remaining marks lie. Within each of the
finitely many cases used below, a suitable choice of coefficients $c$
gives an inequality throughout that case. A \emph{certificate} here
supplies coefficients for an identity that proves such an inequality.
The verifier reconstructs the identity and checks the required signs
over a whole interval of proposed eigenvalues. The accompanying proofs
specify the identities and their inputs.

\subsection{Why five is exceptional}
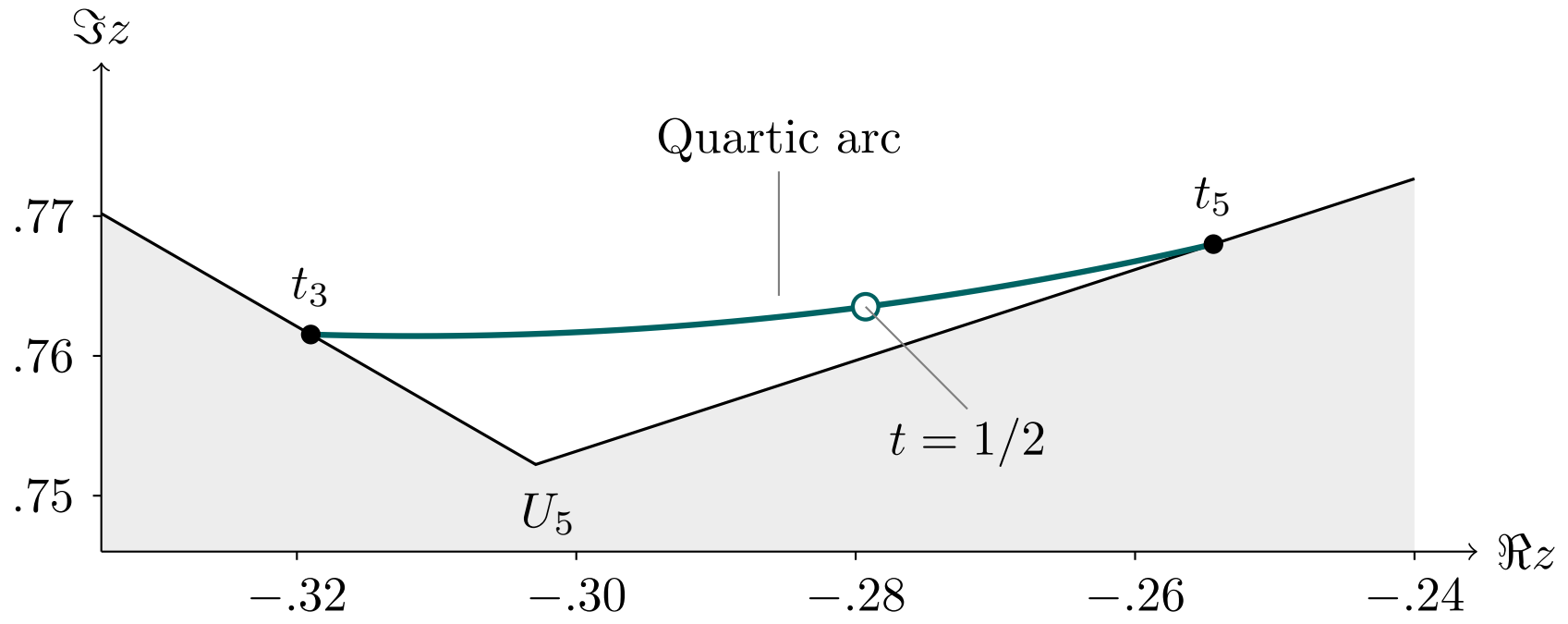
\begin{figure}[htbp]
\centering
\begin{tikzpicture}[x=1cm,y=1cm,font=\small]
\fill[gray!14] (0,0)--(0.00000,2.41853) -- (3.10967,0.62316) -- (9.40000,2.66701)--(9.4,0)--cycle;
\draw[black,semithick] (0.00000,2.41853) -- (3.10967,0.62316) -- (9.40000,2.66701);
\draw[very thick,teal!78!black] plot coordinates {(1.49933,1.55289) (1.62650,1.54943) (1.75411,1.54661) (1.88214,1.54445) (2.01060,1.54295) (2.13947,1.54211) (2.26875,1.54194) (2.39843,1.54244) (2.52851,1.54363) (2.65897,1.54550) (2.78982,1.54807) (2.92104,1.55133) (3.05263,1.55530) (3.18457,1.55997) (3.31687,1.56536) (3.44950,1.57146) (3.58248,1.57829) (3.71578,1.58585) (3.84940,1.59413) (3.98332,1.60316) (4.11755,1.61292) (4.25208,1.62343) (4.38689,1.63469) (4.52197,1.64670) (4.65732,1.65946) (4.79292,1.67299) (4.92878,1.68728) (5.06487,1.70233) (5.20120,1.71816) (5.33774,1.73475) (5.47450,1.75212) (5.61145,1.77027) (5.74860,1.78920) (5.88593,1.80891) (6.02343,1.82941) (6.16109,1.85069) (6.29890,1.87276) (6.43686,1.89562) (6.57495,1.91927) (6.71316,1.94371) (6.85148,1.96894) (6.98990,1.99497) (7.12842,2.02179) (7.26702,2.04941) (7.40569,2.07782) (7.54442,2.10702) (7.68320,2.13702) (7.82202,2.16782) (7.96087,2.19941)};
\draw[->] (0,0)--(9.85,0) node[right] {$\Re z$};
\draw[->] (0,0)--(0,3.5) node[above] {$\Im z$};
\foreach \x/\v in {1.4/{-.32},3.4/{-.30},5.4/{-.28},7.4/{-.26},9.4/{-.24}}{
 \draw (\x,0)--(\x,-.06) node[below] {$\v$};}
\foreach \y/\v in {.4/.75,1.4/.76,2.4/.77}{
 \draw (0,\y)--(-.06,\y) node[left] {$\v$};}
\node at (3.2,.27) {$U_5$};
\node at (4.85,2.94) {Quartic arc};
\draw[gray] (4.85,2.72)--(4.85,1.83);
\fill (1.49933,1.55289) circle (2pt) node[above=2pt] {$t_3$};
\fill (7.96087,2.19941) circle (2pt) node[above=2pt] {$t_5$};
\filldraw[draw=teal!75!black,fill=white,thick] (5.47074,1.75164) circle (2.6pt);
\draw[gray] (5.47074,1.75164)--(6.2,1.02);
\node[below,fill=gray!14,inner sep=1pt] at (6.2,1.0) {$t=1/2$};
\end{tikzpicture}
\caption{An enlarged part of the upper boundary at order five.
The two straight sides belong to $\Pi_3$ and $\Pi_5$. The curve
plots the upper root of $q_t$ between the two certified crossings.
The marked root at $t=1/2$ lies strictly above the polygonal boundary.
Pointwise root evaluations draw the curve. The certificate bounds it over
the entire parameter interval.}\label{fig:quartic}
\end{figure}

Take the following doubly stochastic matrix, $0\le t\le1$:
\[
 A(t)=\begin{pmatrix}
 0&1&0&0&0\\0&0&1&0&0\\t&0&0&1-t&0\\
 1-t&0&0&0&t\\0&0&0&t&1-t
 \end{pmatrix}.
\]
Expanding the determinant gives $(z-1)q_t(z)$, where
\[
 q_t(z)=z^4+tz^3+t(1-t)z^2-t^2z+3t(1-t)-1.
\]
For $.468\le t\le.552$, write $z(t)$ for its root with imaginary
part above $1/2$. Moving $t$ carries this root through a side of the
pentagon and later into the triangle. Let $t_5$ and $t_3$ be the
respective crossings, with $t_5\approx .47052750$ and
$t_3\approx .54900128$. Figure~\ref{fig:quartic} shows why a drawing of
the whole unit disk can easily conceal the extra eigenvalues.

The marked point in Figure~\ref{fig:quartic} can be checked separately. At $t=1/2$, the quartic multiplied
by four is $4z^4+2z^3+z^2-z-1$. Writing $z=x+iy$ and dividing its
imaginary equation by $y$ gives
\[
 y^2=(16x^3+6x^2+2x-1)/(16x+2).
\]
Substitution into the real equation gives a polynomial in $x$ with
opposite signs at $-.279293$ and $-.279292$. Over this interval
the positive square root lies between $.7634$ and $.7637$.
The entire rectangle lies above an edge of each of the triangle,
square and pentagon. The short program \texttt{five\_witness.py}
checks these signs with fractions. This already proves the strict
inclusion $U_5\subsetneq D_5$, known from~\cite{five}.

\subsection{From the marks to an outer bound}

The matrices $A(t)$ attain the whole quartic arc, and their real entries
give its conjugate as well. How do we fill the rays to these arcs, and
how can we rule out points farther along them? We use the radial
contraction from the proof of Theorem~3 in~\cite{pm}. Write $J_n$ for the matrix with every
entry $1/n$. A nonreal eigenvalue $z$ has a centred marking $u$.
Thus $sA+(1-s)J_n$ realizes $sz$ on the \emph{same} marks for
$0<s<1$, and every entry is strictly positive. We call such a matrix
a \emph{positive realization} of $sz$ on $u$. Consequently, if a given
type of marking has no positive realization at a point of a ray,
no realization of that type lies farther out. A test point may lie
strictly inside the proposed region. The real axis is already settled:
a row at a mark of largest modulus gives $|z|\le1$, and $\Pi_2=[-1,1]$.

For the outer inclusion we may assume the marks are distinct and no
three lie on a line. This is what \emph{general position} means here. Indeed,
an exterior point remains exterior after a sufficiently small radial
contraction, so take $A>0$. The two columns of
$X=[\Re u,\Im u]$ sum to zero individually. They are independent:
otherwise $u$ would be a complex multiple of a real vector, and the
equation for the real matrix $A$ would force $z$ to be real.
For any nearby matrix $X'$ with centred columns the real matrix
\[
 T=I+(X'-X)(X^TX)^{-1}X^T
\]
satisfies $TX=X'$, $T\mathbf1=\mathbf1$ and
$\mathbf1^TT=\mathbf1^T$. For $X'$ sufficiently close to $X$, the
matrix $T$ is close to the identity and invertible. Thus $TAT^{-1}$
still has all row and column sums one, has positive entries and has
the same eigenvalue on the perturbed marks.
Avoiding coincidences and collinear triples requires avoiding only
finitely many proper polynomial equalities.

Here is an elementary bound common to both orders. Suppose that
multiplication by $z=x+iy$, with $y>0$, preserves a polygon with
$m$ vertices and with zero in its interior. Write the vertices
counterclockwise as $r_j e^{i\theta_j}$, and put
$\delta_j=\theta_{j+1}-\theta_j\in(0,\pi)$, taking
$\theta_{m+1}=\theta_1+2\pi$ at the last edge. The image of vertex $j$ lies to the left of edge
$j,j+1$. After rotating that vertex to the positive real axis,
the signed-area inequality becomes
\[
 \frac{r_j}{r_{j+1}}
 \le \cos\delta_j+c\sin\delta_j,\qquad c=\frac{1-x}{y}.
\]
The right side is positive. On its positive interval inside $(0,\pi)$,
the logarithm of $\cos\delta+c\sin\delta$ is concave
(its second derivative is $-(1+c^2)/(\cos\delta+c\sin\delta)^2$).
Multiply around the polygon and use $\sum\delta_j=2\pi$ to obtain
$1\le\cos(2\pi/m)+c\sin(2\pi/m)$. In particular,
\begin{equation}\label{small:sector}
 y\le(1-x)\cot(\pi/m).
\end{equation}
If the images of the vertices lie strictly inside, every edge
inequality is strict. For a triangle there is also the bound
$x\ge-1/2$. The barycentric coordinates of the three images form
a row-stochastic matrix with eigenvalues $1,z,\bar z$, whose trace
$1+2x$ is nonnegative. These two inequalities put $z$ in $\Pi_3$.

For the remaining cases, put
\[
 P_k(u)=\operatorname{conv}\Big\{\sum_{j\in S}u_j:|S|=k\Big\}.
\]
This is an assignment polygon with $k$ unit coefficients. A
permutation decomposition gives $zP_k\subseteq P_k$. For a positive
matrix every transported generator lies strictly inside whenever
$P_k$ has interior. To see this, first subtract a sufficiently small
positive multiple of the average of all permutations, which is $J_n$.
After rescaling the remainder, the decomposition in the proof of
Proposition~\ref{prop:assign} gives a representation of the original matrix in
which every permutation has positive weight. In the resulting expression
for a generator's image, every generator occurs.

\paragraph{Five marks.}
Write $\rho_m=e^{2\pi i/m}$. The point
$L_{m,e}(s)=(1-s)\rho_m^e+s\rho_m^{e+1}$ runs along one side of
$\Pi_m$ as $s$ goes from zero to one. A five-vertex marking
uses the singleton and pair-sum polygons $P_1,P_2$. Inequality
\eqref{small:sector} bounds the side $[1,\rho_5]$. On the other
upper sides, namely $L_{5,1}([0,1])$ and $L_{5,2}([0,1/2])$,
the certificate forces the five marks to form a regular pentagon
at every nonreal test point. At these points a vertex is transported
to the pentagon boundary, which a positive matrix cannot do.
The ray argument puts all such multipliers in $\Pi_5$.

A quadrilateral marking has an additional constraint. If $u_i$ is an
outer vertex, let $c_i(w)$ be the least weight of $u_i$ among all
convex combinations of the five marks giving $w$. This quantity,
called a \emph{deletion cap}, vanishes exactly when the other marks can
supply $w$ by themselves. Row $j$ provides
one such combination for $zu_j$, so $c_i(zu_j)\le a_{ji}$. The
column sum gives $\sum_j c_i(zu_j)\le1$. The two diagonals of
the quadrilateral split it into four triangles. Relabel the outer
marks $0,1,2,3$ so that the fifth mark, labelled $4$, is in the
triangle next to side $01$. In that case $P_2$ has successive vertices
$01,14,12,23,03,04$, where $ij=u_i+u_j$. These oriented edges and the
inequalities for $c_i$ give the polynomial constraints on this case.

The first quadrant follows from \eqref{small:sector} with $m=4$.
The following list covers the second quadrant in increasing order
of argument. Here $E=L_{5,1}$, $V=L_{5,2}$ and
\[
 S(s)=-s+i(1-s),\quad T(s)=-s+i(1+s)/\sqrt3,\quad
 U(s)=-1/2+is,\quad L(s)=-1+s+is.
\]
All decimal endpoints in the table are exact rationals.

\begin{center}\small
\begin{tabular}{@{}lll@{}}
\toprule
Test point & Parameter interval & Consequence for a quadrilateral marking\\
\midrule
$S(s)$ & $[0,.2215]$ & A square with a central mark is forced.\\
$E(s)$ & $[.4745,.5021]$ & No such marking exists.\\
$z(t)$ & $[.468,.552]$ & The marking of $A(t)$ is forced.\\
$T(s)$ & $[.32,.5]$ & No such marking exists.\\
$U(s)$ & $[.6875,\sqrt3/2]$ & No such marking exists.\\
$E(s)$ & $[.7237,1]$ & No such marking exists.\\
$V(s)$ & $[0,.337]$ & No such marking exists.\\
$L(s)$ & $[0,.192]$ & A square is forced at nonreal points.\\
\bottomrule
\end{tabular}
\end{center}

The square cases force a transported mark onto a square edge. In
the quartic case two of the marks have to be transported onto other
vertices. Up to scaling they are $1,z,z^2$ in succession, as the
first two rows of $A(t)$ show. None of these cases admits a positive
realization. The angles of consecutive tests either overlap strictly
or agree at the cubic or quintic root of unity. For example, the
determinant of the vectors from zero to the start of the first $E$
window and the end of the $S$ window exceeds $4.8\cdot10^{-5}$.
Its sign places the start of the next window before the end of the
preceding one. The same checks cover the other joins, and confirm
that the short continuations of $z(t)$ beyond its two crossings
lie in $\Pi_5$ or $\Pi_3$. The straight test points are attained by the
cycle construction of Section~\ref{sec:object}, and $A(t)$ attains the curved
ones. The matrices giving attainment may have different marks from
those excluded in the table. The ray argument, followed by conjugation,
proves
\[
 D_5=U_5\ \cup\
 \{s z(t),\,s\overline{z(t)}:0\le s\le1,\ t_5\le t\le t_3\}.
\]

\subsection{What changes with the sixth mark}

The matrix $A(t)$ can be extended by one fixed state. Its exceptional
eigenvalue at $t=1/2$ still occurs, but now the entire rectangle
found above lies in $\Pi_6$, since $|x|<1/2$ and $|y|<\sqrt3/2$.
In fact $D_6=U_6=\Pi_4\cup\Pi_5\cup\Pi_6$. We describe how the
tests now separate the four possible hull sizes.

If all six marks are extreme, the vertices of $P_k$, $k=1,2,3$,
are the sums of $k$ consecutive marks in their cyclic order. A
supporting direction selects the $k$ largest projected coordinates,
which are consecutive around the boundary. In particular $P_3$ has
at most six vertices and $P_3=-P_3$ by complementation. Both $z$ and
$-z$ preserve it. Applying \eqref{small:sector} and conjugation gives
$|y|\le\sqrt3(1-|x|)$. The only possible escape from $\Pi_6$ in the
upper half-plane is therefore across its horizontal side. The
six-vertex certificate on that side forces a regular hexagon except
at the left unit endpoint. A positive realization is excluded by the
contact at a side, and a positive matrix cannot realize a nonreal unit
root since all its singleton images lie strictly inside. Radial
contraction gives the required inclusion for this hull type. Three
extreme marks already give $\Pi_3$ by the triangle argument.

For the two remaining hulls, write $\mathcal K_h$ for the region of
ordinary row-stochastic eigenvalues of order $h$. An invariant
$h$-gon places $z$ in $\mathcal K_h$, by using barycentric coordinates
for its $h$ vertex images. The upper boundary of $\mathcal K_h$ runs from $1$ to $-1$.
Its intervening unit points have arguments $2\pi f$, where $f$
runs through the reduced fractions in $(0,1/2)$ with denominator
at most $h$, in increasing order. Between adjacent
fractions having denominators $q<s$, put $m=\lfloor h/q\rfloor$.
The Karpelevi\v c theorem selects the branch between the corresponding
unit roots of
\[
 z^s(z^q-t)^m=(1-t)^m z^{mq},\qquad 0\le t\le1.
\]
The region is filled radially up to the selected branch. Theorem~1.1
and Sections~2--6 of~\cite{vg} give the proof; the complete paper is
included in the supplement as \texttt{stochastic-regions.pdf}.
For example, the last arc has $(q,s,m)=(2,3,2)$ for $h=4$ and
$(2,5,2)$ for $h=5$. Thus its equations are respectively
$(z^2-t)^2=(1-t)^2z$ and $z(z^2-t)^2=(1-t)^2$.
No doubly stochastic conclusion enters this bound.

For four extreme marks the only rays left by this ordinary bound
are those between the rays through $i$ and $-1/6+i\sqrt3/2$.
On every other ray the ordinary bound already lies in $U_6$.
Indeed, the first arc is the segment $[1,i]$. On the next arc,
from $i$ to $e^{2\pi i/3}$, dividing the polynomial $z^4-tz-(1-t)$ by $z-1$
and taking the imaginary part gives $y^2=3x^2+2x+1$. Hence
$y^2-3/4=3(x+1/2)(x+1/6)\le0$ between $x=-1/2$ and $x=-1/6$.
On the final arc put $b=\sqrt{-2x-1}$. Its equation gives
$y^2-3(x+1)^2=-2b(b-1)^2\le0$. This accounts for the entire
left sector. For the rays still to test, the $U_6$ boundary is first
a square edge, then a horizontal hexagon edge.

Define an assignment polygon by placing coefficients
$1,\rho_h^{-1},\ldots,\rho_h^{-(h-1)}$ on $h$ of the six marks,
and taking the convex hull over all ordered choices of distinct labels.
A cyclic shift rotates this polygon by $\rho_h$. For the eight
possible arrangements of two inner marks in a quadrilateral (the signs
of the oriented triangle areas record which side of each line a mark
occupies),
the certificates force this polygon for $h=4$ to be a nonzero
complex multiple of $\Pi_4$ at every point of
\[
 L_{4,1}([0,.135]),\qquad L_{6,1}([.633,.668]).
\]
Indeed, with $H=\sqrt3/2$, the square piece ends at parameter
$1-H<.135$. The hexagon piece runs from $1/2+(1-H)>.633$
to $2/3<.668$. After scaling the assignment polygon to $\Pi_4$,
its vertex $1$ is sent to $z$. On the square piece $z$ is on a
square edge. On the required part of the hexagon piece, namely
$1/2+(1-H)\le s\le2/3$, it is on or outside the square.
Positivity would put it strictly inside. This excludes a positive realization at
each required point.

For five extreme marks the fifth-order ordinary region has upper
unit vertices $1,\rho_5,i,\rho_3,\rho_5^2,-1$. Its first arc is a
pentagon side. Its shared arc from $i$ to $\rho_3$ enters the hexagon
at $-1/6+i\sqrt3/2$. On the next arc, $z^5-tz^2-(1-t)=0$,
the only additional crossing of the line $y=\sqrt3(x+1)$ is at
$x=(-11+\sqrt{13})/12$. The part from $\rho_3$ to this crossing
is in the hexagon. On the final arc the crossing abscissa lies
between $-.842$ and $-.841$, and the part from there to $-1$ is
in the hexagon. The proof in \texttt{order6/proof.pdf} derives these
crossings by eliminating $t$ on the line and counting real roots.

At each vertex of the outer pentagon, its two incident sides
and the diagonal joining its neighbours enclose a triangle, called
an ear. The
inner mark lies in zero, one or two ears, giving three cell types
up to cyclic relabelling. Here a \emph{cell} comprises markings with the
same signs for the oriented areas of all labelled triples. The following intervals exclude positive
realizations on the remaining rays. The initial and final points
are joined by the full segment indicated; a test segment can be
shorter in radius than the outer boundary on its ray.

\begin{center}\small
\begin{tabular}{@{}lll@{}}
\toprule
Test point & Parameter interval & Number of ears containing the inner mark\\
\midrule
$L_{5,1}(s)$ & $[0,.18]$, $[.83,1]$ & $0,1,2$\\
$L_{5,1}(s)$ & $[.17,.20]$ & $0$\\
$L_{4,0}(s)$ & $[.90,1]$ & $0$\\
$L_{4,0}(s)$ & $[.89,1]$ & $1,2$\\
$L_{4,1}(s)$ & $[0,.14]$ & $0,1,2$\\
$L_{6,1}(s)$ & $[.635,.68]$ & $0,1,2$\\
$L_{6,2}(s)$ & $[.23,.25]$, $[.63,.69]$ & $0,1,2$\\
$L_{5,2}(s)$ & $[0,.23]$ & $0,1,2$\\
\bottomrule
\end{tabular}
\end{center}

The first hexagon crossing on its left sloping side has side
parameter $(5-\sqrt{13})/6\in(.23,.25)$; the final one has parameter
in $(.682,.684)$. Thus the tested portions of $L_{6,2}$ reach both
crossings. The other consecutive windows either share a unit endpoint
or have a positive overlap, as checked by the determinants of their
endpoint vectors. Each tested point lies in $U_6$. Radial contraction
of an exterior realization would give a positive realization on the same
marks at one of these points. This completes the outer inclusion.

\subsection{What the certificates check}

For an oriented edge with endpoints $a(u),b(u)$ of a marked polygon
and a generator $c(u)$ of that polygon, the assertion that $zc(u)$
lies to the left of the edge reads
\[
 H(u;z)=\Im\bigl(\overline{b(u)-a(u)}(zc(u)-a(u))\bigr)\ge0.
\]
The determinant in this formula is twice the signed area of the triangle
with vertices $a,b,zc$. All three generators are linear in the marks,
so $H=w^*Bw$ in independent coordinates $w$ of the centred marking,
where $^*$ denotes conjugate transpose and $B^*=B$. To multiply such
inequalities we need products of coordinates. Collect them in
$V(w)=(w_iw_j)_{i\le j}$ and
$M(w)=(w_i\overline{w_j})_{i,j}$. An identity of the form
\[
 0=V(w)^*RV(w)+M(w)^*QM(w)+\sum_j H_j(w;z)\,w^*G_jw
\]
with positive semidefinite matrices $R,Q,G_j$ expresses zero as a sum
of nonnegative terms on the cell. Indeed, a positive semidefinite
matrix is a Gram matrix $S^*S$, so its quadratic form at $x$ is
$\|Sx\|^2$. If $R$ is positive definite, this square vanishes only
on the zero vector. The identity then
forces $V(w)=0$ and hence $w=0$, which excludes all eigenmarkings.
If $\ker R$ is the line through $V(w_*)$, comparing a nonzero square
coordinate and its cross products forces $w$ to be a complex multiple
of $w_*$. This is the mechanism behind the regular polygons and the
quartic marking. Other identities include the cap inequalities and
products of two valid areas.

At order six we also have to determine when the assignment polygon is
exactly a regular polygon. Let
$a=h^{-1}\sum_{j<h}\rho_h^{-j}u_j$, numbering the marks from zero.
Cyclic assignments put the vertices of $ha\Pi_h$ in the polygon.
To prove the reverse inclusion, fix the normal to one side of
$ha\Pi_h$. Write $C_{ij}$ for the projection in this direction
of $\rho_h^{-i}u_j$ for $i<h$, and zero for the remaining rows.
If numbers $p_i,q_j$ satisfy $C_{ij}\le p_i+q_j$, every assignment
has projected value at most $\sum_i p_i+\sum_j q_j$, independently
of its permutation. The unused amount $p_i+q_j-C_{ij}$ is the
\emph{slack} of entry $ij$. The tests bound these slacks for explicit
$p_i,q_j$ whose sum is attained by the two cyclic assignments at
the ends of the selected side. Rotation supplies the other sides.
For a nonzero marking, the certificates give strictly positive slacks
away from those assignments. If $a=0$, the side normal used in the
tests and all the corresponding slacks would vanish. Thus $a\ne0$.

Bounds on an interval are obtained, for example, by expanding a matrix
polynomial in the Bernstein basis. Rescale the interval to $[0,1]$.
The degree-$m$ basis functions $\binom mr t^r(1-t)^{m-r}$,
$0\le r\le m$, are nonnegative and sum to one by the binomial theorem.
Thus positivity of the coefficient matrices
implies positivity throughout the interval. The implementations
reconstruct the polynomial identities, certify any corrections to
the proposed coefficients and check the required matrix signs.
The two accompanying proofs, \texttt{order5/proof.pdf} and
\texttt{order6/proof.pdf}, specify the forms for each marking type,
the interval joins and the exact algebraic and interval tests. The
order-five driver checks fifteen stages and the order-six driver
checks 87 certificates together with the cell and angular coverage.